\documentclass{article}
\usepackage[margin=1.2 in]{geometry}  
\usepackage{graphicx}
\usepackage{amssymb,amsmath}
\usepackage{epstopdf}
\usepackage{array}
\usepackage{amsthm}
\usepackage[shortlabels,inline]{enumitem} 
\usepackage{color} 
\usepackage[square,numbers,sort&compress]{natbib} 
\usepackage{hyperref}

\usepackage{algorithm}
\usepackage[noend]{algpseudocode}

\numberwithin{equation}{section}

\usepackage{tikz}
\usetikzlibrary{intersections}
\usetikzlibrary{arrows.meta} 
\usetikzlibrary{arrows}
\usetikzlibrary{decorations.markings}
\usetikzlibrary{decorations.text}
\usetikzlibrary{shapes}
\usetikzlibrary{calc}
\tikzset{
    disc/.style={fill=black!40, draw=black,thick}
}
\tikzset{
    discg/.style={fill=#1, draw=black,thick}
}
\tikzset{
    vertex/.style={circle, draw, fill=black!50, inner sep=0pt, minimum width=5pt}
}
\tikzset{
     mystar/.style={star,star points=5,star point ratio=2,fill=blue!20,draw,inner sep=1.5pt}
}
\def\centerarc[#1](#2)(#3:#4:#5)
    { \draw[#1] ($(#2)+({#5*cos(#3)},{#5*sin(#3)})$) arc (#3:#4:#5); }

\usepackage[charter]{mathdesign}
\usepackage[mathcal]{eucal}

\theoremstyle{definition}

\definecolor{jgGreen}{rgb}{0.0, 0.5, 0.0}

\newcommand{\dd}[2]{\frac{d #1}{d #2}}

\newcommand{\pp}[2]{\frac{\partial #1}{\partial #2}}

\newcommand{\R}{\mathbb{R}}
\newcommand{\B}{\mathbb{B}}
\newcommand{\e}[1]{^{(#1)}}
\newcommand{\grad}{\nabla}

\newcommand{\ieq}{I_{\rm eq}}
\newcommand{\iin}{I_{\rm in}}
\newcommand{\jeq}{J_{\rm eq}}
\newcommand{\jin}{J_{\rm in}}
\newcommand{\ngain}[1]{\mathcal N_{\rm gain}(#1)}
\newcommand{\nlose}[1]{\mathcal N_{\rm lose}(#1)}
\newcommand{\nlosesig}[1]{\mathcal N_{\rm lose}^{\sigma\bdy}(#1)}

\newcommand{\lose}{_{\rm lose}}
\newcommand{\gain}{_{\rm gain}}
\newcommand{\same}{_{\rm same}}
\newcommand{\bdy}{_{\rm bdy}}
\newcommand{\opt}{_{\rm opt}}
\newcommand{\sigtan}{\sigma_{\rm tan}}
\newcommand{\egain}{^{\rm gain}}
\newcommand{\elose}{^{\rm lose}}
\newcommand{\esame}{^{\rm same}}

\newcommand{\rev}{_{\rm rev}}
\newcommand{\nfcns}{n_{\rm fcns}}
\newcommand{\nmanif}{\texttt{nmanif}}

\title{Simulating sticky particles: A Monte Carlo method to sample a stratification}
\author{Miranda Holmes-Cerfon\footnote{Courant Institute of Mathematical Sciences, New York University, USA. Email: holmes@cims.nyu.edu.}}

\begin{document}

\maketitle


\begin{abstract}
Many problems in materials science and biology involve particles interacting with strong, short-ranged bonds, that can break and form on experimental timescales. Treating such bonds as constraints can significantly speed up sampling their equilibrium distribution, and there are several methods to sample probability distributions subject to fixed constraints. We introduce a Monte Carlo method to handle the case when constraints can break and form. More generally, the method samples a probability distribution on a stratification: a collection of manifolds of different dimensions, where the lower-dimensional manifolds lie on the boundaries of the higher-dimensional manifolds. We show several applications of the method in polymer physics, self-assembly of colloids, and volume calculation in high dimensions. 
\end{abstract}


\section{Introduction}

Simulating interacting particles is slow: a mere handful of colloidal particles with attractive pairwise interactions can take days or even weeks to collect accurate statistics. 
The trouble with colloidal particles, and related systems like $C_{60}$ molecules or atoms with covalent interactions,
is they form strong but short-ranged bonds \cite{Manoharan:2015ko,Girifalco:1992gg,Ryckaert:1977wr,Andersen:1983wg,Barth:1995vi} 
causing particles to jiggle rapidly around their current configuration as if they are attached by stiff springs that are gently plucked. 
Resolving these vibrating springs in a molecular dynamics or Monte Carlo simulation requires tiny time or space steps, steps that are usually much smaller than the scales of the interesting rearrangements, where springs break, form, or change their relative positions. 


A natural idea for allowing larger steps is to freeze the vibrations, by adding distance constraints between bonded particles: the springs become rods, which can't change their lengths. 
By evolving the system only in directions that preserve the distance constraints, one can take much larger steps, and accelerate the sampling significantly. Such ideas are the origin of the popular Shake and Rattle algorithms and their variants in molecular dynamics \cite{Ryckaert:1977wr,Andersen:1983wg}, and have led to several Monte-Carlo methods that sample a probability distribution subject to holonomic constraints \cite{VandenEijnden:2006gp,Ciccotti:2007fv,Lelievre:2012id,Byrne:2013ed,Zappa:2018jy,Lelievre:2019be}. 
Still, these methods can only be used in limited situations because they require the constraints, or bonds, to be fixed throughout the entire simulation. None allow constraints to be added or dropped, so none can simulate particles with stiff bonds that can break and form. 

This paper will introduce a Monte Carlo method to simulate particles with constraints -- distance constraints, angle constraints, or other kinds of scalar constraints -- that can break and form over the course of the simulation. A natural application is to study \emph{sticky} particles, i.e. those that interact only when their surfaces are exactly in contact. Sticky particles are a model for particles with very short-ranged attractive interactions, such as DNA-coated colloids, which have diameters on the order of 1$\mu$m but interact attractively over ranges around $20$nm or less \cite{HolmesCerfon:2013jw,Perry:2015ku,Wang:2015ep}. DNA-coated colloids are widely studied for their potential use as building blocks for new materials \cite{Rogers:2016bd}, but their short-ranged interactions make them a challenge to simulate. Current constraint-based Monte Carlo techniques cannot be applied because the bonds between the particles must break and form as they equilibrate.

The method we introduce samples more general systems than sticky particles: it can sample a probability distribution for any system restricted by holonomic constraints that can be added or removed. In mathematical language, it samples a probability distribution defined on a \emph{stratification}: a union of open manifolds of different dimensions, where the lower-dimensional manifolds lie on the boundaries of higher-dimensional manifolds in a nice enough way \cite{goresky1988stratified,Goresky:2012jg}. 
For example, a cube is a stratification, consisting of one three-dimensional manifold (its interior), six two-dimensional manifolds (its faces), twelve one-dimensional manifolds (its edges), and eight zero-dimensional manifolds (its vertices.) 
For our sampling method to work, each manifold in the stratification must be the level set of a set of scalar functions, the \emph{constraints} (such as bond-distance constraints), and the manifolds must be connected to each other by adding or removing constraints. 
That is, starting from a particular manifold, if a new constraint is added to the set of constraints which define it (e.g. a bond is formed), then one obtains a new manifold with one dimension lower than the starting manifold, and if a constraint is removed from this set (e.g. a bond is broken), one obtains a manifold with one dimension higher. We consider a stratification formed by taking the union of manifolds connected by adding or removing constraints, and construct a Monte Carlo algorithm to sample a probability distribution defined on this stratification.   
We call our algorithm the Stratification Sampler.

Our method is a natural extension of a well-known sampling algorithm, a random walk Metropolis scheme on a manifold 
\cite{Zappa:2018jy,Lelievre:2012id}. 
As in many of the constraint-based methods referenced above, the method sometimes proposes moves that preserve the current set of constraints. However, it additionally attempts to add or remove constraints: it removes each constraint with constant rate, but proposes to add a constraint only when the configuration is close enough to the boundary corresponding to that constraint. For sticky particles, this means the method proposes to break each bond with a constant rate, and it proposes to form a new bond between particles only when the surfaces of these particles are close enough together.  We show this choice of proposal gives 100\% acceptance probability for a stratification consisting of two flat manifolds (such as a plane and a line) defined by affine constraints. For general manifolds, we still obtain a high acceptance probability, so the proposals don't waste computation proposing moves that are likely to be rejected.

Overall, our method is a form of reversible jump Monte Carlo (RJMC), a method that is commonly used in Bayesian computations to determine the number of parameters in a model  \cite{Green:1995dq}. However, RJMC as it is typically presented requires a parameterization of each of the manifolds that it samples, whereas our method works even without a parameterization.

Our method builds on others that accelerate simulations of systems of particles with strong, short-ranged attractive interactions. Two methods that are highly effective for fluids with finite interaction ranges include Aggregation-Volume-Bias Monte Carlo, which transports single particles into the interaction range of other particles \cite{Chen:2000kr}, and the Geometric Cluster Algorithm of Liu \& Liujten \cite{Liu:2004eb}, which reflects clusters of particles about a pivot point; the latter extends the Swedsen and Wang lattice cluster algorithm to off-lattice systems. 
A third method to simulate fluids but when the particles are perfectly sticky, proceeds by choosing one sphere at random and adding or subtracting up to 3 bonds at a time, choosing randomly from all the geometrically possible ways to do so \cite{Seaton:1986fj,Seaton:1987cn,Kranendonk:1988,Miller:2004ba}. 
Moving beyond fluids, Virtual Move Monte Carlo moves particles with finite-ranged  interactions as a cluster, by displacing one particle initially,  moving other particles if their bond energies changed significantly, and iterating \cite{Whitelam:2007ci}. This method is effective for isolated structures as well as fluids, and is frequently used to simulate DNA \cite{Doye:2013ij}. Finally, a method to sample cross-linked networks of polymers proceeds by adding or removing cross-linked bonds, then regrowing segments of polymers in between the cross-linked sites, treating all interactions as if they were sticky \cite{DeGernier:2014iu,Oyarzun:2018hz}. 

What distinguishes the algorithm to be presented here, is first, its generality -- it can handle any kind of constraints, not just those representing stiff pairwise interactions, but also other kinds of constraints that may arise in biology, machine learning, and geometry calculations. Second, it allows the interactions to have infinitesimal range, a feature of only a small number of the previous studies. Finally, it proposes moves that are physically realistic, moving a system along its internal degrees of freedom in steps that can become arbitrarily small. Hence, we expect this method could be adapted to be a consistent discretization of a set of dynamical equations. Of the above methods to simulate strongly interacting particles, only Virtual Move Monte Carlo comes closest to proposing physically realistic moves, however because this method moves groups of particles rigidly, it cannot propose every motion that could occur dynamically.

Here is an overview of the paper.  Section \ref{sec2:stratification} provides the setup necessary to describe the sampler, by introducing the concept of a mathematical stratification and the notation to describe probability distributions on it. We start with the special case of sticky spheres in Section \ref{sec2:stickysetup}, and move to the more general setup with nearly arbitrary smooth constraints in Section \ref{sec2:generalsetup}. 
Section \ref{sec2:overview} gives an overview of our sampling algorithm, the Stratification Sampler.  
Section \ref{sec2:examples} gives a collection of examples that both illustrate different kinds of stratifications and show how our sampler performs for different systems. It also compares our sampler's data to that from more conventional Brownian dynamics simulations, and shows our sampler can give a good approximation to a system that is not perfectly sticky. 
Section \ref{sec2:algorithm} gives the details of the algorithm that are necessary to implement it, and Section \ref{sec2:flat} proves our proposals  lead to 100\% acceptance probability for a stratification consisting of two flat manifolds defined by affine constraints with no inequalities. 
Section \ref{sec:conclusion} discusses some additional applications of the Stratification sampler as well as ideas for extending and improving it. 
Details on how to implement the algorithm are contained in an Appendix.


\section{What is a stratification?}\label{sec2:stratification}

In this section we introduce the setup and notation required to describe our sampling algorithm. We describe the kind of configuration space we wish to sample, a mathematical stratification, as well as the probability distribution defined on it. We start by explaining these objects for a specific example, sticky spheres, in Section \ref{sec2:stickysetup}, and then we explain how these objects generalize to systems with other kinds of constraints, in Section \ref{sec2:generalsetup}. 

\subsection{Example: Sticky Spherical Particles} \label{sec2:stickysetup}

\begin{figure}
\def\cb{black!80}
\def\ca{black!50}
\def\cc{black!20}
\def\cd{black!60}
\begin{minipage}{0.33\textwidth}
\begin{tikzpicture}[scale=1]
    \draw [<->,very thick] (0,4) node (yaxis) [above] {$|y_1-y_2|$}
        |- (5,0) node (xaxis) [below,pos=0.94] {$|x_1-x_2|$};
 \draw[thick,\cd] (3,0) arc (0:90:3) ;
\node at (1.4,0.5) {$q_{12}(x)<0$};
\node at (3.7,3.55) {$q_{12}(x)>0$};
\node[\cd] at (3.2,2) {$q_{12}(x)=0$};
\begin{scope}[shift={(1.1,1.5)},scale=0.7]
\filldraw[discg=\ca] (0.7,0) circle (0.5) node{2};
\filldraw[discg=\cb] (0,0) circle (0.5) node {1};
\end{scope}
\begin{scope}[shift={(1.2,3.5)},scale=0.7]
\filldraw[discg=\ca] (1.3,0) circle (0.5) node{2};
\filldraw[discg=\cb] (0,0) circle (0.5) node {1};
\end{scope}
\end{tikzpicture}
\end{minipage}
\begin{minipage}{0.65\textwidth}
\centering
\begin{tikzpicture}[scale=0.8]
\begin{scope}[shift={(-0.5,0)}]
\filldraw[discg=\ca] ({cos(60)+0.3},{sin(60)+0.3}) circle (0.5) node{2};
\filldraw[discg=\cb] (0,0) circle (0.5) node {1};
\filldraw[discg=\cc] ({1+1.5},{0+0.3}) circle (0.5) node {3};
\end{scope}
\node[align=center] at (0.5,-1.3) {$M_\emptyset$ \\ no discs in contact};
\begin{scope}[shift={(5.5,0)}]
\begin{scope}[shift={(-0.5,0)}]
\filldraw[discg=\ca] ({cos(60)},{sin(60)}) circle (0.5) node{2};
\filldraw[discg=\cb] (0,0) circle (0.5) node {1};
\filldraw[discg=\cc] ({1+1.5},{0+0.3}) circle (0.5) node {3};
\end{scope}
\node[align=center] at (0.5,-1.3) {$M_{12}$ \\ discs 1-2 in contact};
\end{scope}
\begin{scope}[shift={(0,-4)}]
\begin{scope}[shift={(-0.25,0)}]
\filldraw[discg=\ca] ({cos(60)},{sin(60)}) circle (0.5) node{2};
\filldraw[discg=\cb] (0,0) circle (0.5) node {1};
\filldraw[discg=\cc] ({cos(60)+cos(-30)},{sin(60)+sin(-30)}) circle (0.5) node {3};
\end{scope}
\node[align=center] at (0.5,-1.3) {$M_{12,23}$ \\discs 1-2, 2-3 in contact};
\end{scope}
\begin{scope}[shift={(5.5,-4)}]
\filldraw[discg=\ca] ({cos(60)},{sin(60)}) circle (0.5) node{2};
\filldraw[discg=\cb] (0,0) circle (0.5) node {1};
\filldraw[discg=\cc] ({1},{0}) circle (0.5) node {3};
\node[align=center] at (0.5,-1.3) {$M_{12,23,31}$\\ all pairs of discs in contact};
\end{scope}
\end{tikzpicture}
\end{minipage}
\caption{Left: illustration of the distance function $q_{12}(x)=(|x_1-x_2|^2+|y_1-y_2|^2)^{1/2}-\sigma$ for two discs with positions $x=(x_1,y_1,x_2,y_2)$ and diameter $\sigma$. When $q_{12}(x)=0$, the system lies on a lower-dimensional subset of its configuration space, a grey semicircle above. Right: some examples of configurations from selected manifolds in the configuration space for three sticky discs. 
The configuration space also contains manifolds $M_{23},M_{31},M_{12,31},M_{23,31}$, not shown.
}\label{fig:setupexample}
\end{figure}
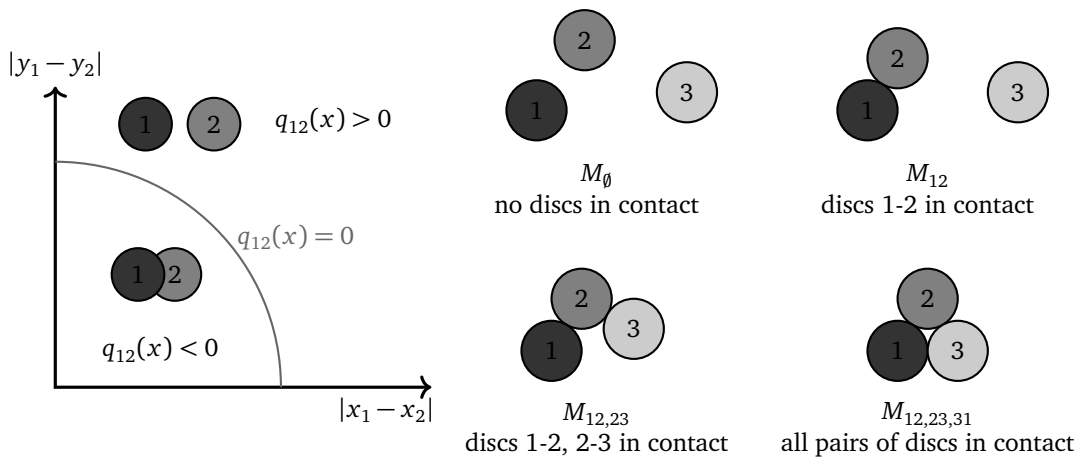

Consider a system of $N$ spheres in $\R^d$ with radii $r_1,r_2,\ldots, r_N$, at temperature $T$. The spheres are \emph{sticky} when the pair potential $U_{ij}(r)$ between spheres $i,j$ is such that
\begin{equation}\label{stickypotential}
e^{-U_{ij}(r) / k_BT} = \kappa_{ij}\delta(r-\sigma_{ij}) + H(r - \sigma_{ij}) \,,
\end{equation}
where  $\sigma_{ij}=r_i+r_j$. 
Here $\delta$ is a Dirac-delta function, 
which captures the stickiness of the interaction: it causes spheres' centers to spend a non-zero amount of time at a distance of exactly $\sigma_{ij}$, so the spheres' surfaces stick to each other. The parameter $\kappa_{ij}$, which we call the \emph{sticky parameter}, measures the stickiness of interaction $ij$ at temperature $T$, and is larger if spheres spend more time exactly in contact \cite{Baxter:1968dh,HolmesCerfon:2013jw,Perry:2015ku}.\footnote{
The way $\kappa_{ij}$ is defined here gives it units of length. To obtain a dimensionless sticky parameter we must write the argument of the delta function as $\delta(r/\sigma_{ij}-1)$. However, doing so would require carrying around factors of $\sigma_{ij}$ in our later calculations, so we instead work with a dimensional $\kappa_{ij}$ for notational convenience. 
} The function $H$ is a Heaviside function, which equals 1 if its argument is positive and equals 0 otherwise. This function captures both the hard-core interaction that prevents spheres from overlapping, as well as the density for spheres $i,j$ to not be in contact.

The sticky interaction potential \eqref{stickypotential} is a model for particles that interact with a strong and short-ranged attractive interaction. One such interaction is the square-well potential with large depth and small range, which Baxter and others used in their early studies of sticky particles \cite{Baxter:1968dh,Gazzillo:2004jc,Stell:1991va}.
Another family of potentials that can lead to sticky interactions, a family we will return to in our examples, is the Morse potential:
\begin{equation}
U^{\rm morse}(r) = E(1-e^{-\rho(r-\sigma)})^2 - E .
\end{equation}
This potential has an attractive well at distance $\sigma$ with depth $E$  and range proportional to $\rho^{-1}$. 
When  $E\to\infty, \rho\to \infty$ so the potential becomes arbitrarily deep and narrow, this pair potential approaches a sticky potential \eqref{stickypotential}, provided the Boltzmann factor for the particles to be bound approaches a constant: $\int_0^{\sigma+\epsilon} e^{-U^{\rm morse}(r)/k_BT} dr \to \kappa < \infty$ \cite{HolmesCerfon:2013jw,BouRabee:2020gg}.\footnote{
A dimensionless sticky parameter would be defined as $\kappa = \lim_{\rho,E\to\infty}\sigma^{-1}\int_0^{\sigma+\epsilon} e^{-U^{\rm morse}(r)/k_BT} dr$. 
}
 Here $\epsilon$ is a cutoff beyond which the potential is sufficiently close to zero; a reasonable value for the Morse potential is $\epsilon \approx 2.5/\rho$ \cite{Trubiano:2019tl}. The constant $\kappa$ that arises in this limit is the sticky parameter for the interaction. It is proportional to the second virial coefficient, which is used to characterize pair potentials in a cluster expansion \cite{Noro:2000bi,Miller:2004ba,Trubiano:2019tl}.

The configuration space for a system of sticky spheres can be divided into subspaces depending on how many spheres are in contact and which particular spheres are in contact. For concreteness, consider a system of $N=3$ discs in dimension $d=2$, with common diameter $\sigma$ and common sticky parameter $\kappa$. Let $x_1,x_2,x_3\in\R^2$ represent the centers of the discs, and let $x=(x_1,x_2,x_3)\in\R^6$ represent the configuration of the whole system. 
To prevent the system from evaporating, we also assume the discs are confined to a finite space $\B\subset\R^6$; for example we may take $x\in\B=[0,L]^{6}$ so the discs are in a box with length $L$. 
To describe which discs are in contact, it is convenient to introduce the distance functions 
\begin{equation}\label{distancefunction}
q_{12}(x) = |x_1-x_2|-\sigma, \qquad q_{23}(x) = |x_2-x_3|-\sigma, \qquad
q_{31}(x) = |x_3-x_1|-\sigma\,.
\end{equation}
Each function $q_{ij}$ measures the distance between the surfaces of discs $i,j$. When $q_{ij}(x) = 0$, the surfaces of discs $i,j$ are exactly in contact. 

One configuration the discs can adopt has no discs in contact. Let the set of all such configurations be $M_\emptyset$. 
This set can be characterized in terms of the distance functions as 
\[
M_\emptyset =  \{x\in \R^{6}\cap \B: q_{12}(x)>0,\; q_{23}(x)>0,\; q_{31}(x)>0\}.
\]
This set is a 6-dimensional subset of $\R^6$. 
Because of the Heaviside function in the sticky potential \eqref{stickypotential}, there is a finite probability to find the system in $M_\emptyset$. 

Another kind of configuration the discs can adopt has exactly one pair of discs in contact. Let us call these sets of configurations $M_{12}$, $M_{23}$, $M_{31}$ depending on which pair is in contact. These sets are given by 
\begin{align*}
M_{12}&= \{x\in \R^{6}\cap \B: q_{12}(x)=0,\; q_{23}(x)>0,\; q_{31}(x)>0\},  \\
M_{23}&= \{x\in \R^{6}\cap \B: q_{23}(x)=0,\; q_{12}(x)>0,\; q_{31}(x)>0\},  \\
M_{31}&= \{x\in \R^{6}\cap \B: q_{31}(x)=0,\; q_{12}(x)>0,\; q_{23}(x)>0\}.
\end{align*}
That is, one of the distance functions equals 0 exactly, while the other distance functions are inequalities. Fixing a distance removes one degree of freedom from the system, as illustrated schematically in Figure \ref{fig:setupexample}, so the sets $M_{12}$, $M_{23}$, $M_{31}$ are 5-dimensional subsets of $\R^6$. In fact, they each form part of the boundary of $M_\emptyset$, since if a configuration is in $M_\emptyset$, then by moving discs  $i,j$ so they become arbitrarily close without touching, the configuration approaches $M_{ij}$. Mathematically, these sets are 5-dimensional \emph{manifolds}, which means they could be mapped out locally using only 5 variables instead of 6 variables, if we were clever enough to find such a parameterization (we won't need to.) 

Because of the delta-function in the sticky potential \eqref{stickypotential}, there is   finite probability to find the system in any one of the manifolds $M_{12}$, $M_{23}$, $M_{31}$. That is, there is finite probability to find the system on a lower-dimensional subspace than the ambient space, $\R^6$. This property makes the sticky potential more challenging to sample than smoother interaction potentials, since any sampler must be able to propose moves that jump exactly to these lower-dimensional subspaces. 

Continuing, we may form lower-dimensional subspaces by considering more discs in contact. When two pairs of discs are in contact we obtain sets 
\begin{align*}
M_{12,23}&= \{x\in \R^{6}\cap \B: q_{12}(x)=0,\; q_{23}(x)=0,\; q_{31}(x)>0\},  \\
M_{23,31}&= \{x\in \R^{6}\cap \B: q_{23}(x)=0,\; q_{31}(x)=0,\; q_{12}(x)>0\},  \\
M_{31,12}&= \{x\in \R^{6}\cap \B: q_{31}(x)=0,\; q_{12}(x)=0,\; q_{23}(x)>0\}.
\end{align*}
Each of these sets is given by fixing two distances, so is a 4-dimensional manifold. Notice that each set is part of the boundary of a 5-dimensional manifold in the configuration space; for example $M_{12,23}$ is on the boundary of both $M_{12}$ and $M_{23}$. 
When all discs are in contact we obtain
\[
M_{12,23, 31} = \{x\in \R^{6}\cap \B: q_{12}(x)=0,\; q_{23}(x)=0,\; q_{31}(x)=0\}\,,
\]
which is a 3-dimensional manifold; it is the set of all states obtained by rotating and translating the triangular configuration of discs. As before, there is a finite probability to find the system in any of the 4- or 3-dimensional manifolds in the configuration space. 

The full configuration space for the 3 sticky discs is
\begin{equation}\label{Sdiscs1}
\mathcal S = M_\emptyset \cup M_{12} \cup M_{23} \cup M_{31} \cup M_{12,23} \cup M_{23,31} \cup M_{31,12} \cup M_{12,23, 31}\,.
\end{equation}
That is, it is a union of manifolds, each obtained by considering some collection of spheres in contact. 
This set $\mathcal S$ is a particular example of a \emph{stratification}, a union of manifolds of different dimensions, where the lower-dimensional manifolds lie on the boundaries of the higher-dimensional manifolds. 

We may write $\mathcal S$ more concisely by defining the set of possible pairs in contact to be 
\begin{equation}\label{Isticky}
\mathcal I = \{\emptyset,\{12\},\{23\},\{31\},\{12,23\},\{23,31\},\{31,12\},\{12,23,31\} \},
\end{equation}
and letting $I$ index some element in $\mathcal I$, for example $I=\{23\}$, or $I=\{12,23,31\}$. Then, from \eqref{Sdiscs1}, we have
\begin{equation}\label{Sdiscs2}
\mathcal S = \bigcup_{I\in \mathcal I} M_I\,.
\end{equation}

A system of sticky discs in equilibrium has finite probability to be found in any of the manifolds $M_I$ that are part of its configuration space $\mathcal S$. To write down this probability specifically we consider all pairwise interactions of the form \eqref{stickypotential}, and use the distance functions \eqref{distancefunction}, to obtain the Boltzmann density on $\mathcal S$ as 
\begin{equation}\label{rhos1}
\rho(x) = Z^{-1}\prod_{i<j} \kappa\delta(q_{ij}(x)) + H(q_{ij}(x)) \,.
\end{equation}
Here $Z$ is a normalization constant so the density integrates to 1. By expanding the product, we may write $\rho(x)$ as a sum of densities on each individual manifold as  
\begin{equation}\label{rhos2}
\rho(x) = Z^{-1}\sum_{I\in \mathcal I} \kappa^{|I|}H_{I^c}(x)\delta_{M_I}(x).
\end{equation}
Here $\delta_{M_I}(x) = \prod_{ij\in I} \delta(q_{ij}(x))$ is a singular density on $M_I$, and $H_{I^c}(x) = \prod_{ij\notin I}H(q_{ij}(x))$ enforces the non-overlap condition for pairs that are not exactly in contact, i.e. it equals 1 if no pairs are overlapping, and it equals 0 otherwise. 
The quantity $|I|$ is the size of $I$, equal to the number of pairs of discs in contact (the number of bonds.)

Expression \eqref{rhos2} shows that the probability to find the system with a given number of bonds $m$ is proportional to $\kappa^m$. That is, the manifolds may be thought of as constant-energy surfaces: on a given manifold $M_I$ with $|I|=m$, the probability density is constant (with respect to the delta-function on $M_I$), so the energy is constant on this surface. If the system then gains a bond, it moves to a manifold with density proportional to $\kappa^{m+1}$, effectively changing the energy by a discrete amount measured by $\kappa$. 

Yet another way to write the Boltzmann distribution, one that will prove necessary to implement our algorithm, is to write it in terms of the natural surface measure $\mu_I$ on each manifold $M_I$. This measure, also called the Hausdorff measure, is used to measure the volume of a manifold and to calculate integrals over it. To convert between delta-functions and the natural surface measure, we use the coarea formula \cite{Federer:2014uh,Morgan:2016tv}, which says that $\prod_{ij\in I}\delta(q_{ij}(x))dx = |Q_I^TQ_I|^{-1/2}\mu_I(dx)$, where 
\begin{equation}\label{Qsticky}
Q_I(x) = (\grad q_{ij}(x))_{ij\in I}
\end{equation}
 is the matrix whose columns are the gradients of the distance functions. For example, $Q_{12,23} = (\grad q_{12}\;\;\grad q_{23})$ is a $6\times2$ matrix. 
Using the coarea formula we may write \eqref{rhos2} as \cite{Ciccotti:2007fv}
\begin{equation}\label{rhos3}
\rho(x) = Z^{-1}\sum_{I\in \mathcal I} \kappa^{|I|}H_{I^c}(x)|Q_I^TQ_I|^{-1/2}\mu_{I}(dx).
\end{equation}
Physically, the matrix $Q_I^TQ_I$ arises when, instead of fixing the distances, we treat them as stiff springs, and construct an energy as $E(x)=\sum_{ij\in I} \frac{1}{2}kq_{ij}^2(x)$ for some spring constant $k$. Then on $M_I$, the Hessian of the energy is $\grad\grad E(x)=kQ_I(x)Q_I^T(x)$. The nonzero eigenvalues of $Q_IQ_I^T$ are the same as the nonzero eigenvalues of $Q_I^TQ_I$, so the logarithm of $|Q_I^TQ_I|^{-1/2}$ is proportional to the vibrational entropy. 

The fact that $\prod_{ij\in I}\delta(q_{ij}(x))dx \neq \mu_I(dx)$ stems from the difference between so-called ``hard'' constraints and ``soft'' constraints, that has been much discussed in the statistical mechanics literature (e.g. \cite{Fixman:1974dd}.) Hard constraints are those that are imposed directly on the equations of motion, and in this context would give rise to the equilibrium probability $\mu_I(dx)$ with no prefactor. Soft constraints arise by imposing constraints with increasingly stiff springs, which in the limit of infinite stiffness don't have the same equilibrium probability as hard constraints.

\subsection{General Mathematical Setup}\label{sec2:generalsetup}

In this section we introduce notation to handle more general stratifications and their associated probability distributions, not just those that come from distance  constraints between particles. For example, we may wish to consider fixed bond angles between particles, bonds between particles of different shapes such as elliptical particles, or, we may have an entirely different application that doesn't involve particles, but still involves sets of variables whose values may sometimes be collectively frozen. The notation in this section will be abstract, and we refer the reader back to Section \ref{sec2:stickysetup} and forward to Section \ref{sec2:examples} for concrete examples.

Let $x\in \R^n$, and let 
\[
\mathcal Q=\{q_1(x), q_2(x), \ldots, q_{|\mathcal Q|}(x)\}
\]
be a collection of continuously differentiable scalar functions. These functions will be either \emph{constraints}, when $q_i(x)=0$, or \emph{inequalities}, when $q_i(x)>0$. In Section \ref{sec2:stickysetup}, these functions were the distance functions \eqref{distancefunction}; in this section they can represent angles between bonds, distance functions between non-spherical surfaces, or other kinds of functions.

We wish to form manifolds by choosing some of the functions in $\mathcal Q$ to be constraints, and some to be inequalities. To this end, let  $\ieq\subset \mathcal Q$ be a particular set of constraints, and let $\iin\subset \mathcal Q$ be a particular set of inequalities. To make sense we must have $\ieq\cap\iin = \emptyset$. Write $I = (\ieq,\iin)$, and call this pair of constraints and inequalities the \emph{labels}. 
Define a manifold from this particular choice of labels as 
\begin{equation}
M_I = \{x\in\R^n : q_i(x) = 0 \text{ for } q_i\in \ieq, \;q_j(x) > 0 \text{ for } q_j\in \iin\}\,.
\end{equation}
We have changed notation slightly from Section \eqref{sec2:stickysetup}: we now include the inequalities explicitly in the labels $I$ for each manifold. In Section \eqref{sec2:stickysetup} we did not need to specify the inequalities, because any function in $\mathcal Q$ that was not a constraint was automatically an inequality, however in general we may wish to ignore some function $q_i$, and not use it as either a constraint or an inequality. 

If we let $m_I=|\ieq|$ be the number of constraints defining $M_I$, then generically, $M_I$ will be an open manifold of dimension $d_I = n - m_I$. This statement is not always true for specific, non-generic cases, but we discuss conditions when we are sure it will hold at the end of the section. 

%
%

To form a stratification, define a collection of labels
\begin{equation}
\mathcal I = \{ I\e{1}, I\e{2},\ldots, I\e{|\mathcal I|}\}. 
\end{equation}
Recall this collection was defined for sticky particles in \eqref{Isticky}. 
Form the union of all manifolds with these labels as  
\begin{equation}
\mathcal S = \bigcup_{I\in \mathcal I} M_I\,.
\end{equation}
Then set $\mathcal S$ is a \emph{stratification}, given some additional mild but technical conditions on the relationship between manifolds 
\footnote{
Loosely, if two manifolds $X,Y$ with $\text{dim}(Y)<\text{dim}(X)$ are in the stratification, and $Y \subset \overline X$, i.e. $Y$ is in the closure of $X$, then at every point $y\in Y$, the stratification near $y$ has to look locally like a cone, and, furthermore, the topology of this local picture is the same for all $y\in Y$. 
Specifically, a stratification is usually assumed to satisfy \emph{Whitney's condition B} at all $y\in Y$, which is as follows \cite{goresky1988stratified}. Given $X,Y$ as above, and suppose that (i) sequence $x_1,x_2,\ldots \in X$  converges to $y$, (ii) sequence $y_1,y_2,\ldots \in Y$  converges to $y$, (iii) the secant lines $l_i = \overline{x_iy_i}$ converge to a limiting line $l$ in some local coordinate system near $y$, (iv) the tangent planes $T_{x_i}X$ converge to a limiting plane $\tau$. Then $l \subset \tau$.  
} \cite{goresky1988stratified},  \cite{Goresky:2012jg}. 
We assume these conditions hold, and refer to $\mathcal S$ as a stratification hereafter. 

Now we wish to construct a probability distribution on the stratification $\mathcal S$. 
We construct this distribution from the surface measures on each manifold, in a similar way to \eqref{rhos3}. 
Let $\mu_I(dx)$ be the natural surface measure on $M_I$,  the $d_I$-dimensional Hausdorff measure; it is obtained from the Euclidean measure in the ambient space $\R^n$ by restricting it to $M_I$. 
Suppose we have a collection of scalar functions $\{f_I(x)\}_{I\in \mathcal I}$, where $f_I$ is defined on manifold $M_I$. 
For example, for sticky spheres we might choose $f_I = \kappa^{|I|}|Q_I^TQ_I|^{-1/2}$ as in \eqref{rhos3}. 
We construct a probability distribution from these functions as 
\begin{equation}\label{rho}
\rho(dx) = Z^{-1}\sum_{I\in \mathcal I} f_I(x)\mu_I(dx)\,,
\end{equation}
where $Z = \sum_{I\in \mathcal I} \int_{M_I}f_I(x)\mu_I(dx)$ is the normalizing constant, which we assume is finite. 
The measure $\rho$ is the probability distribution that we wish to sample. 

Our sampling algorithm will rely on manifolds in the stratification being sufficiently connected to each other in a way we now describe. Recall that in Section \ref{sec2:stickysetup}, every manifold (except the highest-dimensional manifold $M_\emptyset$) was on the boundary of a higher-dimensional manifold, because every lower-dimensional manifold was obtained from a higher-dimensional one by adding a distance constraint.  

With this in mind, let $\nlose{I}$ be the set of labels in $\mathcal I$ that can be obtained from a particular label $I$ by adding a single constraint: 
\begin{equation}
\nlose{I} = \{J\in \mathcal I : \;\; \jeq = \ieq \cup \{q\}, \; q\notin \jeq \}.
\end{equation}
We call a label $J$ in $\nlose{I}$ a \emph{Lose neighbour of $I$}, because generically, in moving from $M_I$ to $M_J$, we lose a dimension: $d_J=d_I-1$. 
Similarly, let $\ngain{I}$ be the set of labels in $\mathcal I$ that can be obtained from a particular label $I$ by removing a single constraint: 
\begin{equation}
\ngain{I} = \{J\in \mathcal I : \;\; \jeq = \ieq - \{q\} \}.
\end{equation}
We call a label $J$ in $\ngain{I}$ a \emph{Gain neighbour of $I$}, because generically, in moving from $M_I$ to $M_J$, we gain a dimension: $d_J=d_I+1$. 
We distinguish between two types of Gain neighbours. 
Sometimes a constraint that is removed becomes an inequality, so that $\jeq = \ieq - \{q\}$ with $q\in\jin$. We call such a $J$  a \emph{one-sided Gain neighbour} of $I$. 
Other times, the removed constraint is simply forgotten, so that $\jeq = \ieq - \{q\}$ with $q\notin\jin$. We call such a $J$ a \emph{two-sided} Gain neighbour of $I$.

Our hope is that if $J$ is a Gain neighbour of $I$ (which implies that $I$ is a Lose neighbour of $J$), then $M_I$ is part of the boundary of (the closure of) $M_J$: $M_I\subset\partial \overline M_J$.\footnote{For a two-sided Gain neighbour it may be that $M_I\subset\overline M_J$; that is, $M_I$ is part of the closure of $M_J$, but not necessarily part of its boundary. 
In the two-sided case we will still refer to $M_I$ as being part of the boundary of $M_J$, because it may be treated as if it were a boundary in our algorithm.} This is a property our algorithm will rely on to pass between manifolds. If it doesn't always hold, our algorithm will still work, provided this property holds for enough manifolds that the sampler can move between all manifolds in the stratification by following Gain and Lose neighbours.

With the setup complete, we now return to the delicate issue of what assumptions are needed to ensure that $M_I$ is a $d_I$-dimensional manifold, and hence, to ensure the probability distribution \eqref{rho} is well-defined. 
There are pathological cases where $M_I$ is \emph{not} a manifold of the correct dimension; for example the manifold $M_I = \{(x,y):x-y=0, \; 2x-2y=0\}$ is a line (1-dimensional), not a point (0-dimensional), whereas the manifold $M_I = \{(x,y,z): z -(x^2+y^2)=0, \; z +(x^2+y^2)=0\}$ is a point, not a one-dimensional curve. For less trivial examples involving sticky spheres, see \cite{HolmesCerfon:2016wa,Kallus:2017hi}. 

To avoid these situations we must make an assumption on the constraints, that we call the \emph{regularity assumption}. Form the $n\times m_I$ matrix whose columns are the gradients of the constraints:
\begin{equation}\label{Q}
Q_I(x) = \left( \grad q_i(x) \right)_{i\in \ieq} .
\end{equation}
Recall this matrix was defined for sticky particles in \eqref{Qsticky}. 
The regularity assumption is that for each $I\in \mathcal I$, the columns of $Q_I$  are linearly independent everywhere in $M_I$. When this assumption holds, the implicit function theorem implies that each $M_I$ is an open manifold of dimension $d_I = n - m_I$, and therefore each $\mu_I$ and hence the full probability distribution $\rho(x)$ are well-defined. This assumption will also be necessary later to implement our sampling algorithm. 
The regularity assumption doesn't always hold in applications, and we explain  what goes wrong with our algorithm when it doesn't hold in the conclusion.

\section{Overview of the Stratification Sampler}\label{sec2:overview}

In this section we give an overview our algorithm, the Stratification Sampler, for 
 generating a Markov chain $X_1,X_2,\ldots \in \mathcal S$ with stationary distribution $\rho$ as in \eqref{rho}. We present only the details needed to understand the basic ideas of the algorithm and the parameters in the subsequent examples, leaving the details necessary to implement the algorithm for Section \ref{sec2:algorithm}.

The algorithm keeps track of both a point and the current manifold as a pair $(X_k,I_k)$, with $X_k\in M_{I_k}$. 
Suppose the current point is $(X_k,I_k)=(x,I)$. 
The next step  $(X_{k+1},I_{k+1})$ of the Markov chain is generated by first constructing a proposal move $(y,J)$, and then accepting or rejecting it using a Metropolis-Hastings step. 
 
 The proposal move is constructed in two steps:
 \begin{enumerate}[(i),noitemsep]
 \item Choose a label $J$ with probability $\lambda_{IJ}(x)$;
 \item Choose a point $y \in M_J$ with probability density $h_{IJ}(y|x)$ with respect to $\mu_J(dx)$. 
 \end{enumerate}
 Here $\lambda_{IJ}(x)$ is the probability of proposing a label $J$, given the current label is $I$ and current point is $x$, and must be such that $\sum_{J\in \mathcal I}\lambda_{IJ}(x) \leq 1$ for each $(x,I)$. 
 The function $h_{IJ}(y|x)$ is the probability density of proposing a point $y\in M_J$, given the current label is $I$, the current point is $x$, and the proposed label is $J$. 
 Sometimes the method fails to produce a proposal point $y$, in which case the proposal is set to $(x,I)$. Therefore it is usually the case that $\int_{M_J}h_{IJ}(y|x)\mu_J(dy) < 1$. 
 
Given a proposal $(y,J)$, it is accepted according to a Metropolis-Hastings rule with acceptance probability
\begin{equation}\label{acc}
a(y,J|x,I) = \min\left( 1, \frac{f_J(y)\lambda_{JI}(y)h_{JI}(x|y)}{f_I(x)\lambda_{IJ}(x)h_{IJ}(y|x)} \right)\,.
\end{equation}
That is, the algorithm generates $U\sim \mbox{Unif}([0,1])$, and if $U<a(y,J|x,I)$ the proposal is accepted and the next step is $X_{k+1} = y$, $I_{k+1} = J$. Otherwise, the proposal is rejected and the next step is $X_{k+1} = x$, $I_{k+1} = I$. 
If the proposal itself was $(x,I)$, it is accepted automatically. 
Importantly, as with most Monte-Carlo samplers, this method can be implemented without knowing the overall normalization constant $Z$, since $Z$ cancels out in the acceptance ratio above. 

We show in Appendix \ref{sec:DB} that if $X_k \sim \rho$, then $X_{k+1}\sim \rho$. 
Therefore $\rho$ is the stationary distribution for the Markov chain above. 
We expect that for reasonable choices of proposals the algorithm is ergodic, so that for any initial condition, $X_n\to \rho$ as $n\to \infty$. 

The acceptance rule \eqref{acc} will look familiar to anyone who has worked with replica exchange methods, Reversible Jump Monte Carlo, or other Monte Carlo methods that switch between a discrete set of continuous spaces. 
Yet, implementing this acceptance rule for a stratification requires several additional considerations. 
One, is that the proposal densities in \eqref{acc} are calculated with respect to reference measures on different manifolds, so one must be careful to account for the Jacobians of the maps from one manifold to another.  
Another is that there must be a unique map between a starting point $(x,I)$ and a proposal point $(y,J)$, a nontrivial consideration when jumping between manifolds of different dimensions. 
Finally, finding proposal moves that ensure the acceptance probability \eqref{acc} is neither too high nor too low is nontrivial. 
The average acceptance probability should be large enough that computation isn't wasted generating proposal moves, but if it is too large, successive samples are highly correlated, so many samples must be generated to obtain approximately independent draws from $\rho$.




We now present a choice of label and point proposals that we have found to be effective. We introduce these proposals here, instead of in the more detailed Section \ref{sec2:algorithm}, so the reader can understand the parameters $\lambda\gain, \lambda\lose, \sigma, \sigma\bdy, \sigtan$ that we refer to in our examples. A reader not interested in these parameters can skip to the examples in Section \ref{sec2:examples}. 

While these proposals may seem unmotivated, we show in Section \ref{sec2:flat} that they give an acceptance probability of $a(y,J|x,I) = 1$ for two manifolds defined by affine constraints with $f(x) = cst$ and no inequalities. Therefore, when the step sizes for the proposals are small, we expect to have high acceptance probabilities for curved manifolds as well.

\subsection{Label proposal}\label{sec2:label}
To choose a label $J$, we consider three different types of proposals:
\begin{itemize}[noitemsep]
\item \emph{Same}: choose $J=I$.  \\
This retains the current labels, and takes a step on the current manifold. 
\item \emph{Gain}: choose $J\in \ngain{I}$.  \\
This moves to a manifold of one higher dimension, obtained by dropping a single constraint. 
\item \emph{Lose}: choose $J\in \nlose{I}$. \\
This moves to a manifold of one lower dimension, on the boundary of the current one, obtained by adding a single constraint. 
\end{itemize}
The overall probabilities of Same, Gain, and Lose moves are $\lambda\same,\lambda\gain,\lambda\lose$ respectively. In practice we set $\lambda\gain,\lambda\lose$, and then calculate $\lambda\same=1-\lambda\gain-\lambda\same$. If a particular manifold $I$ has no Gain neighbours or no Lose neighbours, then the probability $\lambda\same$ must be adjusted for that manifold. 

We must also decide how to choose the particular label $J$ for Gain and Lose moves. For a Gain move, we choose $J$ uniformly among all the possible Gain neighbours, $J\sim\mbox{Unif}(\ngain{I})$. For a Lose move, choosing labels uniformly doesn't work well; briefly, this is because most boundaries are far away from the current point, so the sampler usually fails to find a point on a randomly chosen boundary, as we discuss in more detail in the example in Section \ref{ex:isotropic}. Therefore, to help the sampler find boundaries and propose smaller, more physically realistic moves, we only propose moves to boundaries that are estimated to be sufficiently close to the current point.
Letting $\nlosesig{I,x}$ be the set of nearby Lose neighbours to a point $(x,I)$, we choose the label uniformly from this set, $J\sim\mbox{Unif}(\nlosesig{I,x})$. 
 
The set of nearby Lose neighbours is the set of neighbours estimated to lie within a distance  $\sigma\bdy$ from the current point:
\begin{equation}\label{Nlosesig}
\nlosesig{I,x} = \{J\in \nlose{I}: h_{\rm opt}(x,I,q_{I\to J}) < \sigma\bdy\}\,.
\end{equation}
The function $h_{\rm opt}(x,I,q_{I\to J})$ is the estimated distance from the current point $(x,I)$ to the closest point on a Lose neighbour $J$ formed by adding constraint $q_{I\to J}$. This distance is estimated by linearizing each possible additional constraint $q_{I\to J}$ as explained in Section \ref{sec2:algorithm}. 

Proposing labels requires three parameters: $\lambda\gain,\lambda\lose,\sigma\bdy$. 
In some cases it is efficient to relate these parameters. Specifically, if it is known that for manifolds $M_I, M_J$, the function we are sampling is $f(x) \approx c_I$ on  $M_I$ and $f(x)\approx c_J$ on $M_J$, where $c_I,c_J$ are constants, then we should instead choose $\lambda\gain = (c_J/c_I)\sigma\bdy\lambda\lose$to depend on $I,J$, see Section \ref{sec:01d}. We usually set $\lambda\gain = \sigma\bdy\lambda\lose$ in our implementation.


\subsection{Point proposal}\label{sec2:point}

Once we have chosen a label $J$, we must propose a point $y\in M_J$. Here is an overview of how we do this for each of the three move types. 
In addition to the parameter $\sigma\bdy$ introduced earlier, these moves depend on two new parameters, $\sigma$ and $\sigtan$. 

Our moves depend on taking steps in the tangent space to each manifold in the stratification.  The tangent space to manifold $M_I$ at $x$, which we write as $\mathcal T_{x,I}$, is the set of vectors $v$ satisfying $Q_I^T(x)v = 0$, where $Q_I$ was defined in \eqref{Q}. It is the linear approximation to a manifold $M_I$ near $x$, i.e. it approximates a surface locally as a plane. The set of vectors perpendicular to $\mathcal T_{x,I}$ is the normal space, $\mathcal N_{x,I} = \mathcal T_{x,I}^\perp$. This space is spanned by the columns of $Q_I$. 
Under the regularity assumption, the dimension of $\mathcal T_{x,I}$ is $d_I$, the dimension of the manifold, and the dimension of $\mathcal N_{x,I}$ is $m_I$, the number of constraints. 


To represent a given tangent space we construct a matrix 
$T_{x,I}\in \R^{n\times d_I}$ whose columns form an orthonormal basis of $\mathcal T_{x,I}$. The matrix $T_{x,I}$ can be computed using standard linear algebra operations, such as from the QR decomposition of $Q_I$ as described in Appendix \ref{sec:summary}, Algorithm \ref{alg:tangent}. 



\begin{figure}
\begin{minipage}{0.47\textwidth}
Gain move\\
\end{minipage}\hfill
\begin{minipage}{0.47\textwidth}
Lose move\\
\end{minipage}
\begin{minipage}{0.47\textwidth}
\begin{tikzpicture}[scale=1]
\def\ra{15};    
\def\raa{15.07}; 
\def\rb{17.5};  
\def\rbb{17.45}; 
\def\vn{1.5};  
\def\vt{1};     
\def\th{7};   
\centerarc[ultra thick](0,{-\ra})(80:100:\ra);
\centerarc[dotted](0,{-\ra})(80:100:\rb);
\node[vertex,label=below:{x}] (x) at (0,0) {};
\node[vertex,label=right:{y}] (y) at (\vt,\vn) {};
\draw[->,thick] (x) -- node[midway,right] {$V$} (y);
\draw[->,dashed] (x) --node[midway,left]{$V_n$} (0,\vn);
\draw[->,dashed] (0,\vn) --node[midway,above]{$V_t$} (y);
\node[] at (3.35,-0.25) {$q(x)=0$};
\draw[|-|]  ({\raa*cos(\th+90)},{-\ra+\raa*sin(\th+90)}) -- node[midway,left] {$\sigma\bdy$} ({\rbb*cos(\th+90)},{-\ra+\rbb*sin(\th+90)});
\end{tikzpicture}
\end{minipage}\hfill
\begin{minipage}{0.47\textwidth}
\begin{tikzpicture}[scale=1]
\def\ra{15};    
\def\raa{15.07}; 
\def\rb{17.5};  
\def\rbb{17.45}; 
\def\th{7};   
\def\rr{0.8};  
\def\vv{0.7};  
\centerarc[name path = qq, ultra thick](0,{-\ra})(80:100:\ra);
\centerarc[dotted](0,{-\ra})(80:100:\rb);
\node[] at (3.35,-0.25) {$q(x)=0$};
\draw[|-|]  ({\raa*cos(\th+90)},{-\ra+\raa*sin(\th+90)}) -- node[midway,left] {$\sigma\bdy$} ({\rbb*cos(\th+90)},{-\ra+\rbb*sin(\th+90)});
\node[vertex,label=above:{x}] (x) at (0,{1+\vv}) {};
\draw[->,dashed] (x) -- node[midway,left]{$v\opt$} (0,\vv);
\draw[->,dashed] (0,\vv) -- node[midway,below]{$R$} (\rr,\vv);
\path[name path = losedir] (x) -- ({3*\rr/sqrt(\rr^2+1)},{1+\vv-3*1/sqrt(\rr^2+1)});
\draw[dotted,name intersections={of=qq and losedir, by=y}] (x) -- (y);
\node[vertex,label=below:{y}] at (y) {};
\draw[->,thick] (x) -- node[pos=0.3,right]{$V$} ({\rr/sqrt(\rr^2+1)},{1+\vv-1/sqrt(\rr^2+1)});
\end{tikzpicture}
\end{minipage}
\caption{Schematic illustrating the Gain and Lose tangent step proposals, for a one-sided boundary $q(x)=0$ in $\R^2$. Gain moves propose a step $V$ with components $V_n\sim\mbox{Unif}( [0,\sigma\bdy])$ in the normal direction to the boundary, and $V_t\sim \mathcal N(0,\sigtan^2V_n^2)$ in the tangential direction. 
Lose moves propose a unit vector $V$ by summing the vector $v\opt$ in the (estimated) normal direction to the boundary, and a vector with component $R\sim \mathcal N(0,\sigtan^2)$ in the tangential direction, and rescaling to obtain a unit vector. 
Lose moves only attempt to jump to boundaries that are no more than an (estimated) distance of $\sigma\bdy$. 
}\label{fig:LoseGain}
\end{figure}
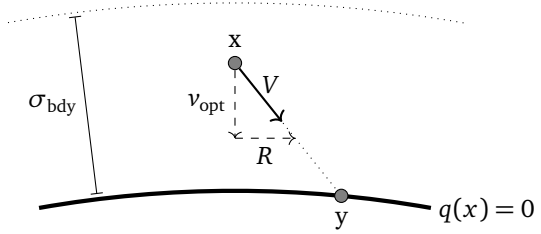

\subsubsection{Same Move: $x\in M_I$, $y\in M_I$}

To generate a proposal on $M_I$, we start by choosing a vector $v$ in the tangent space $\mathcal T_{x,I}$, and taking a step in this direction, as $x\to x+v$. This usually pushes us off the manifold, so we  project back onto $M_I$ by finding a vector $w$ such that 
\begin{equation}\label{qsolveSame}
q_i(x+v+w) = 0 \qquad \text{for all constraints $i\in\ieq$}.
\end{equation}
The proposal is $y=x+v+w$. Requiring $w$ to be in the normal space $\mathcal N_{x,I}$  gives this proposal several nice properties \cite{Zappa:2018jy}. 

Equation \eqref{qsolveSame} is solved numerically using a nonlinear equation solver. If this solver fails to find a solution, the proposal is $(x,I)$. If the solver finds a solution but it doesn't satisfy the inequalities, the proposal is also $(x,I)$. If the solver finds a solution $(y,J)$ that satisfies the inequalities and is accepted using the Metropolis-Hastings rule \eqref{acc}, the algorithm must still check that  starting from $(y,J)$ the nonlinear equation solver converges to $(x,I)$, a step we call the reverse check. This step is necessary to ensure the reverse move could have been proposed hence has nonzero probability density; see Section \ref{sec2:algorithm} and \cite{Zappa:2018jy}. 

We remark that in our implementation, we solve \eqref{qsolveSame} using Newton's method in dense arithmetic, using the LU decomposition to solve the linear equations at each iteration. This is fast for small systems, but will be slow for larger ones, where one should instead use sparse arithmetic and an iterative linear solver. For a discussion of the merits of different numerical methods for imposing constraints, see \cite{Eastman:2010ts}. 

We may choose the vector $v$ according to any density; in our examples we use an isotropic multidimensional Gaussian with standard deviation $\sigma$ in each independent direction. Specifically, from a vector of independent random variables $R\in \R^{d_I}$ we construct the step as
\begin{equation}\label{vSame}
v = T_{x,I} R, \qquad R_i \sim N(0,\sigma^2), \quad i=1,\ldots,d_I. 
\end{equation}


\subsubsection{Gain move: $x\in M_I$, $y\in M_J$ with $J\in \mathcal N\gain(I)$}

A Gain move is similar to a Same move, but with a different choice of tangent step $v$. 
Suppose that $q$ is the constraint dropped from $\ieq$ to obtain $\jeq$. We have $x\in M_I$, and we wish to propose $y\in M_J$ with $\jeq = \ieq - \{q\}$. 
We do this by dropping constraint $q$ and proposing a move on $M_J$, as if it were a Same move on $M_J$. 

The only difference with a Same move, is we use a different tangent step $v$. An isotropic Gaussian on $\mathcal T_{x,J}$ will work in theory, but in practice is inefficient, because it is not possible to make it symmetric with the Lose moves, as we discuss in Section \ref{ex:isotropic}. Instead, we choose $v$ to have different sizes in the directions normal to and tangential to the boundary. Let $u_n$ be a unit vector normal to the boundary with $M_I$ and in the tangent space $\mathcal T_{x,J}$. This vector, which points in the direction where the dropped constraint $q$ increases most quickly, is found by projecting $\grad q$ onto $\mathcal T_{x,J}$;  see Section \ref{sec2:algorithm}. 
Let the remaining, tangential directions be contained in the columns of $T_{x,I}$. 
Notice that together, $u_n$ and $\mathcal T_{x,I}$ span $\mathcal T_{x,J}$. 

We construct the tangent step $v$ by combining a random step in the direction of $u_n$, with a random step in the directions spanned by $T_{x,I}$, as 
\begin{equation}\label{GainStep}
v = u_nV_n + T_{x,I}V_t. 
\end{equation}
The random variables $V_n\in \R$, $V_t\in \R^{d_J-1}$ are constructed in different ways. In the direction normal to the boundary, the size of the step is  
\begin{equation}
V_n \sim \mbox{Unif}([0,\sigma\bdy])
\end{equation}
for a one-sided Gain move, and $ V_n \sim \mbox{Unif}([-\sigma\bdy,\sigma\bdy])$ for a two-sided Gain move. 
The choice for a one-sided move ensures it always moves directly away from the boundary to where the inequalities are not violated, and the upper bound of $\sigma\bdy$ ensures it does not move so far that $I$ would not be a nearby Lose neighbour of $J$. 

Given the value of $V_n$, we choose a step in the tangential directions as an isotropic Gaussian with standard deviation $\sigtan |V_n|$ in each independent direction, where $\sigtan>0$ is another parameter. 
That is, each component of $V_t$ is chosen as 
\begin{equation}
(V_t)_i \sim N(0,\sigtan^2V_n^2), \qquad 1=1,\ldots,d_J-1.
\end{equation}
For small $\sigtan$, the proposal step $v$ moves mostly away from the boundary. 
See Figure \ref{fig:LoseGain} for a sketch of a Gain move in two dimensions. 

We chose the variance of $V_t$ to increase with the magnitude of the proposed normal step $V_n$ so that Gain and Lose moves are more symmetric. 
The distributions of $V_n$ and $V_t$  balance, respectively, the overall probability $\lambda\lose$ of proposing a Lose move, and the distribution of the tangential components of the reversed Lose step. We discuss this further in Section \ref{sec2:flat}.

\subsubsection{Lose Move: $x\in M_I$, $y\in M_J$ with $J\in \nlosesig{I,x}$}

We have $x\in M_I$, and we wish to propose $y\in M_J$ with $\jeq = \ieq \cup \{q\}$, where $q$ is the new constraint.

We start by choosing a unit vector $v$ in the current tangent plane $\mathcal T_{x,I}$. 
This unit vector gives the direction to step in the tangent plane, but the magnitude of the tangent step is an unknown variable $\alpha\in \R$. Taking a step $x\to x+\alpha v$ usually does not put us on manifold $M_J$, so we project to $M_J$ by finding a vector $w$ in the normal space $\mathcal N_{x,I}$  such that 
\begin{equation}\label{qsolveLose}
q_i(x+\alpha v+w) = 0 \qquad \text{for all constraints $i\in\jeq$}.
\end{equation}
This equation must be solved for $\alpha,w$.
The proposal is $y=x+\alpha v+w$. As for Same and Gain moves, equation \eqref{qsolveLose} is solved numerically and the solver may not converge. If it does converge, one additional consideration is that the solution must have $\alpha>0$, for otherwise the step is not in direction $v$, but in $-v$. A solution with $\alpha <0$ must be rejected so the sampler is reversible. 

One possible way to choose $v$ is uniformly on the unit sphere in the tangent plane. This choice is not efficient, since it gives equal probability to all directions so it ignores the fact that there is a boundary nearby. 
A more efficient proposal comes from estimating the direction that points toward the boundary $M_J$ and moving mostly in this direction.

Specifically, let $v\opt$ be a unit vector in the tangent space $\mathcal T_{x,I}$ that is estimated to give the shortest distance to the boundary, i.e., it is normal to the boundary to a linear approximation. This vector is found by projecting $\grad q$ onto the tangent plane $\mathcal T_{x,I}$ as we explain in Section \ref{sec2:algorithm}. Let  $T_{x,I,v\opt}$ be the matrix whose columns span the remainder of $\mathcal T_{x,I}$, i.e. the part of the tangent space that is perpendicular to $v\opt$. 

We build $v$ by choosing a deterministic component in the direction of $v\opt$, and random isotropic components $R\in \R^{d_J-1}$ in the directions orthogonal to $v\opt$, and then normalizing to get a unit vector:
\begin{equation}\label{vlose}
v = \frac{v\opt + T_{x,I,v\opt}R}{|v\opt + T_{x,I,v\opt}R|}, \qquad R_i\sim N(0,\sigtan^2), \quad i=1,\ldots,d_J-1.
\end{equation}
The parameter $\sigtan$ is the same one used to construct a Gain move. 
If $\sigtan$ is small, the proposal direction is mostly in the direction that is estimated to give the shortest distance to the boundary. 
See Figure \ref{fig:LoseGain} for a sketch of a Lose move in two dimensions.


\section{Examples}\label{sec2:examples}

We now show several examples, that serve both to illustrate different kinds of stratifications, and how the Stratification Sampler works in practice. 
The examples in Sections \ref{ex:parline}-\ref{ex:ellipse} each introduce a different kind of stratification: the union of a parabola and a line and their interiors (\ref{ex:parline}), a trimer of sticky two-dimensional discs (\ref{ex:trimer}), a polymer of 6 sticky spheres whose sticky parameters can vary (\ref{ex:polymer}), a polymer adsorbing weakly to a surface (\ref{ex:polywall}), a high-dimensional shape whose surface area we wish to calculate (\ref{ex:ellipse}). Section \ref{ex:polymer} further shows that the sampler gives a good approximation to the probability distribution for a real system with an interaction potential that does not have infinitesimal range. 
Section \ref{ex:isotropic} shows what goes wrong when one uses chooses tangent steps isotropically in the tangent space. 
The codes to run all examples and reproduce the statistics and figures are available at \cite{gitfigures}.

\subsection{Parabola and line in two dimensions}\label{ex:parline}

\begin{figure}
\centering
\includegraphics[height=5cm]{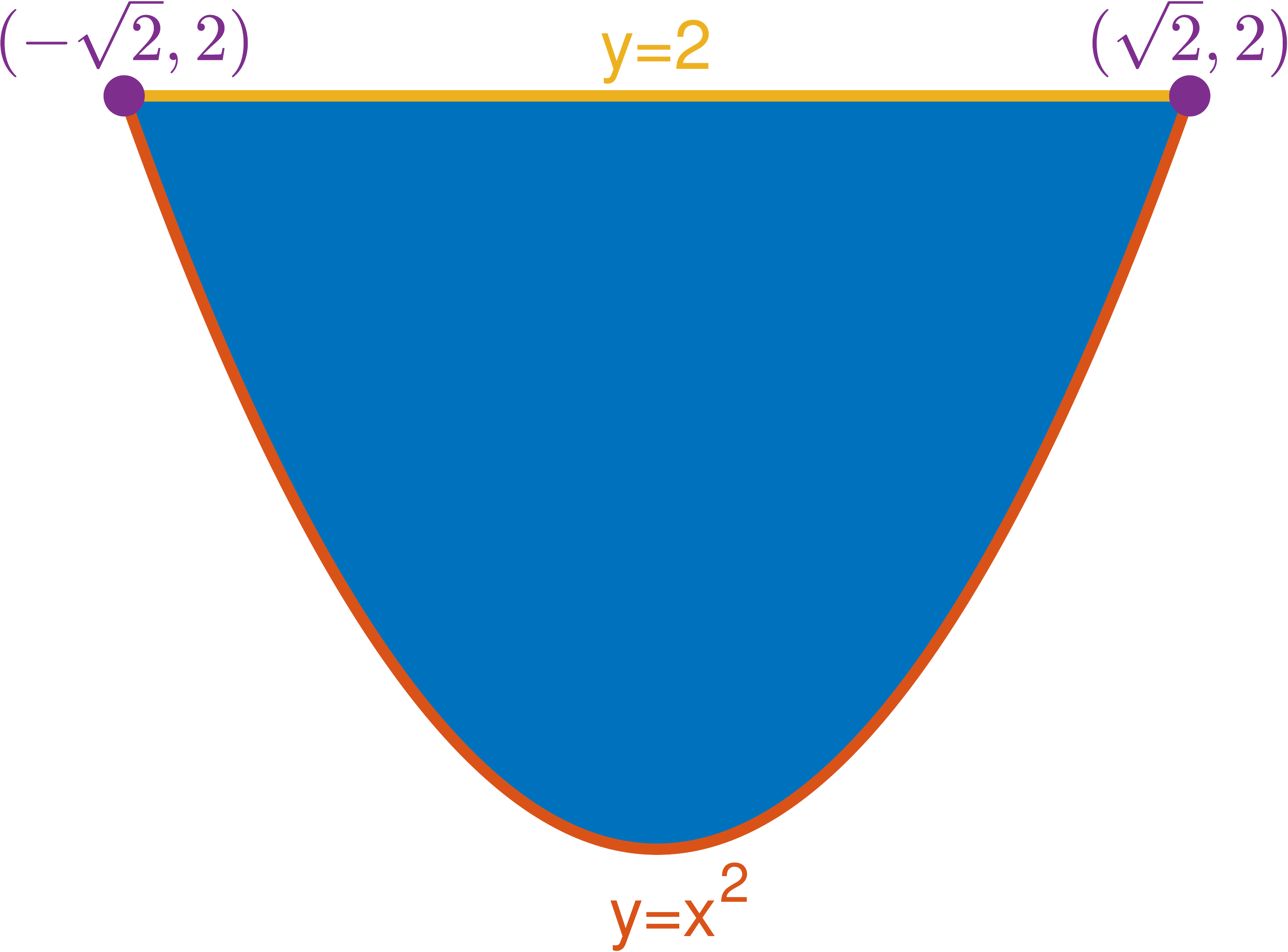}\hspace{1.75cm}
\raisebox{4mm}[0pt][0pt]{\includegraphics[height=4.25cm]{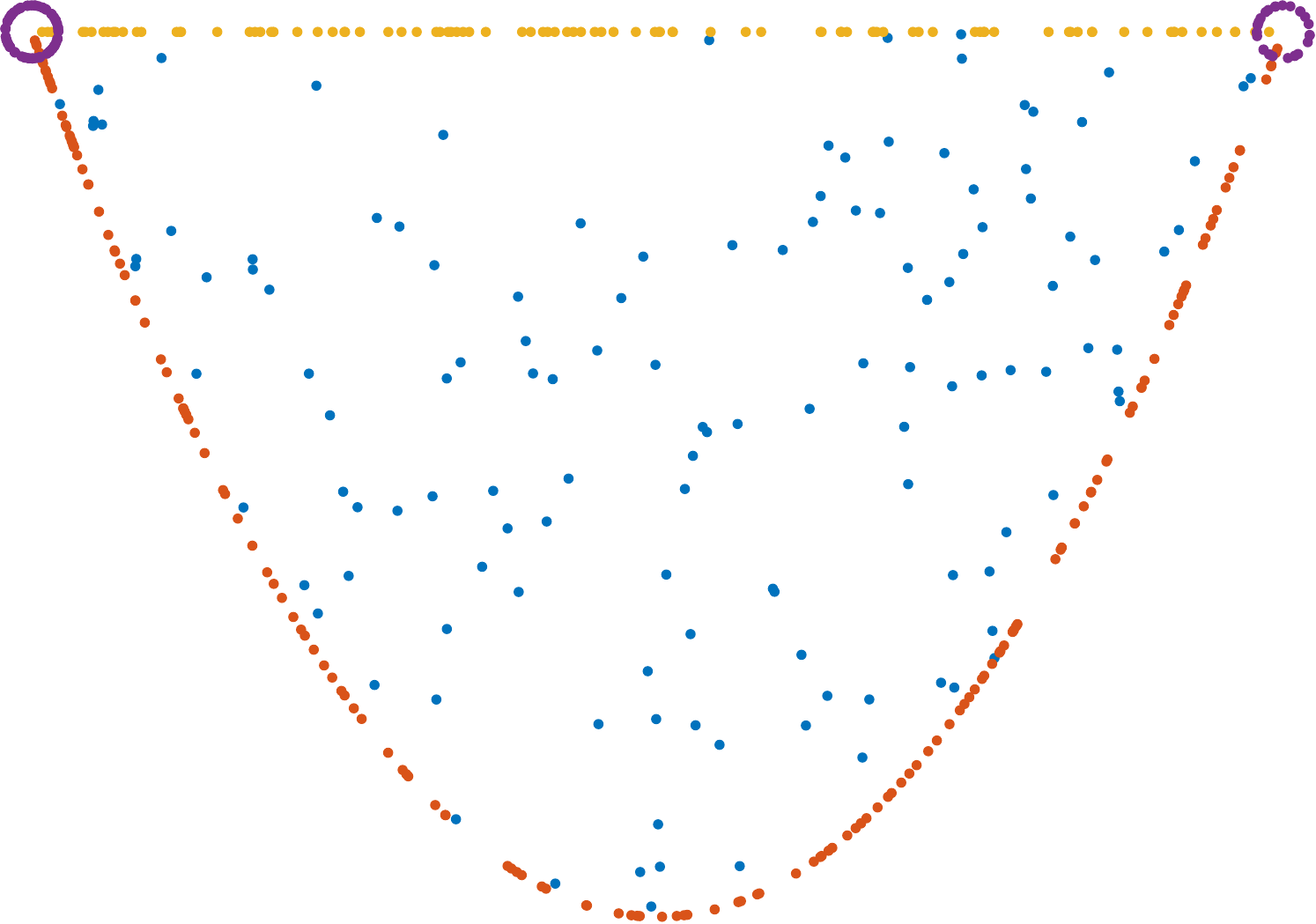}}
\caption{Left: Stratification from Example \ref{ex:parline}. $M_1$ (blue) is two-dimensional; $M_2$ (red) is one-dimensional; $M_3$ (yellow) is one-dimensional; $M_4$ (purple) is zero-dimensional. 
Right: Some points sampled from the stratification with probability distribution $\rho(dx)$ defined in \eqref{rhoparline}. Points sampled from the corners are plotted at equally-spaced angles along a small circle in the order they were generated to make them distinct. 
}\label{fig:parline}
\end{figure}

\begin{figure}
\centering
      \begin{minipage}{0.69\textwidth}
      \centering
\includegraphics[width=1\linewidth]{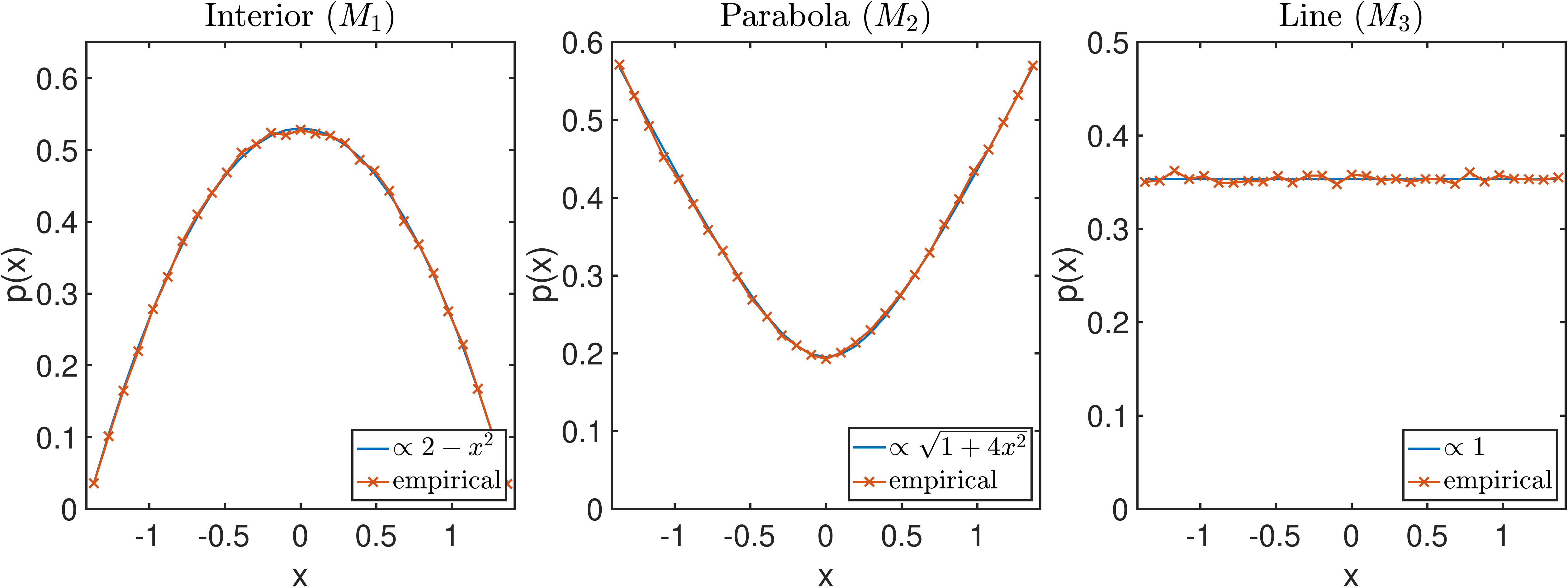}
\end{minipage}\hfill
\begin{minipage}{0.3\textwidth}
\centering
\begin{tabular}{ccc}
\multicolumn{3}{c}{Fraction of points}\\\hline
& Theory & Empirical \\\hline
$M_1$ & 0.2748 & 0.2751\\
$M_2$ & 0.3734 & 0.3726\\
$M_3$ & 0.2061 & 0.2066 \\
$M_4$ &0.1457  & 0.1458 \\\hline
\end{tabular}
\end{minipage}
\caption{Statistics for Example \ref{ex:parline} with probability distribution \eqref{rhoparline}. Plots show the marginal probability distributions in $x$ on each manifold, calculated theoretically (blue line) and empirically by sampling (red line.) The distribution in the left and right corners ($M_4$) was $\{0.50003,0.49997\}$. Table shows the fraction of points in each manifold. The sampling algorithm generated $10^7$ points and stored every 10$^{\rm th}$ point to calculate statistics, using parameters ${\sigma{=}0.9}$, ${\sigma\bdy{=}0.3}$, ${\sigtan{=}0.6}$, ${\lambda\lose{=}0.7}$, ${\lambda\gain{=}0.21}$. The sampling run took 110 seconds on a 3.6 GHz Intel Core i7 iMac.
}\label{fig:parline2}
\end{figure}

Our first example is not physically motivated, but rather designed to illustrate a stratification and its probability measure in a way that can be easily visualized. 
We found this example useful when developing our code, as it contains most of the special cases one must handle (e.g. manifolds with no constraints, zero-dimensional manifolds, manifolds with no Lose or Gain neighbours, manifolds with both Lose and Gain neighbours, flat manifolds, curved manifolds), yet it can be plotted in two dimensions and many statistics are known analytically. 

Let $(x,y)\in \R^2$, and let the constraint and inequality functions be $\mathcal Q = \{q_1(x,y),q_2(x,y)\}$ with
\[
q_1(x,y) = y-x^2,\qquad q_2(x,y) = 2-y\,.
\]
Consider the stratification formed by considering all possible labelings of $q_1,q_2$ as constraints or inequalities: $\mathcal I = I\e{1}\cup I\e{2}\cup I\e{3} \cup I\e{4}$ with 
$I\e{1} = (\emptyset,\{q_1,q_2\})$, $I\e{2} = (\{q_1\},\{q_2\})$, $I\e{3} = (\{q_2\},\{q_1\})$, $I\e{4} = (\{q_1,q_2\},\emptyset)$.   
That is, writing $M_{k} = M_{I\e{k}}$, the manifolds are defined mathematically and named for ease of reference as   
\begin{align*}
\text{``Interior'':} \qquad M_1 &= \{(x,y):y>x^2,y<2\} && \text{(two-dimensional, blue)}\\
\text{``Parabola'':} \qquad M_2 &= \{(x,y): y=x^2,y<2\}&& \text{(one-dimensional, red)}\\
\text{``Line'':} \qquad M_3 &= \{(x,y):y=2,y>x^2\} && \text{(one-dimensional, yellow)}\\
\text{``Corners'':} \qquad M_4 &= \{(x,y):y=x^2,y=2\} && \text{(zero-dimensional, purple.)}
\end{align*}
The stratification is visualized in Figure \ref{fig:parline}, with each manifold a different color as listed above. Notice that $M_1$ is a Gain neighbour of both $M_2$ and $M_3$, and each of $M_2,M_3$ in turn is a Gain neighbour of $M_4$. Since for each function $q_i$ that is not a constraint is an inequality, all Gain neighbours are one-sided Gain neighbours. 

Each manifold comes with a natural surface measure $\mu_{k}(dx) = \mu_{I\e{k}}(dx)$. The most familiar are the measures for the interior, $\mu_1(dx)$, and line, $\mu_3(dx)$, which are the two-dimensional and one-dimensional area measures on the corresponding manifolds respectively. The measure on the parabola,  $\mu_2(dx)$, is the arc-length measure: if $A\subset M_2$ is a connected segment of the parabola, then $\int_A \mu_2(dx)$ equals the arc length of $A$. The corner's measure $\mu_4(dx)$ is the counting measure, which is $\mu_4(A)\in\{0,1,2\}$ depending on how many points from the set $M_4=\{(-\sqrt{2},2), (\sqrt{2},2)\}$ are included in $A$. 

One example of a probability measure on $\mathcal S$ is
\begin{equation}\label{rhoparline}
\rho(dx)  \propto \mu_1(dx) + \mu_2(dx) + \mu_3(dx) + \mu_4(dx)\,.
\end{equation}
This weights points on each manifold according to their arc-length measures. 
Some points sampled from this measure using the Stratification Sampler are shown in Figure \ref{fig:parline}. 
Figure \ref{fig:parline2} shows that for the sampled points, the total fraction of time spent on each manifold, as well the marginal distributions in $x$ on each manifold, agree very well with the analytically calculated distributions. 

\subsection{A trimer of sticky discs in two dimensions}\label{ex:trimer}

\begin{figure}
\centering
\begin{minipage}{0.68\textwidth}\centering
\includegraphics[height=3.7cm]{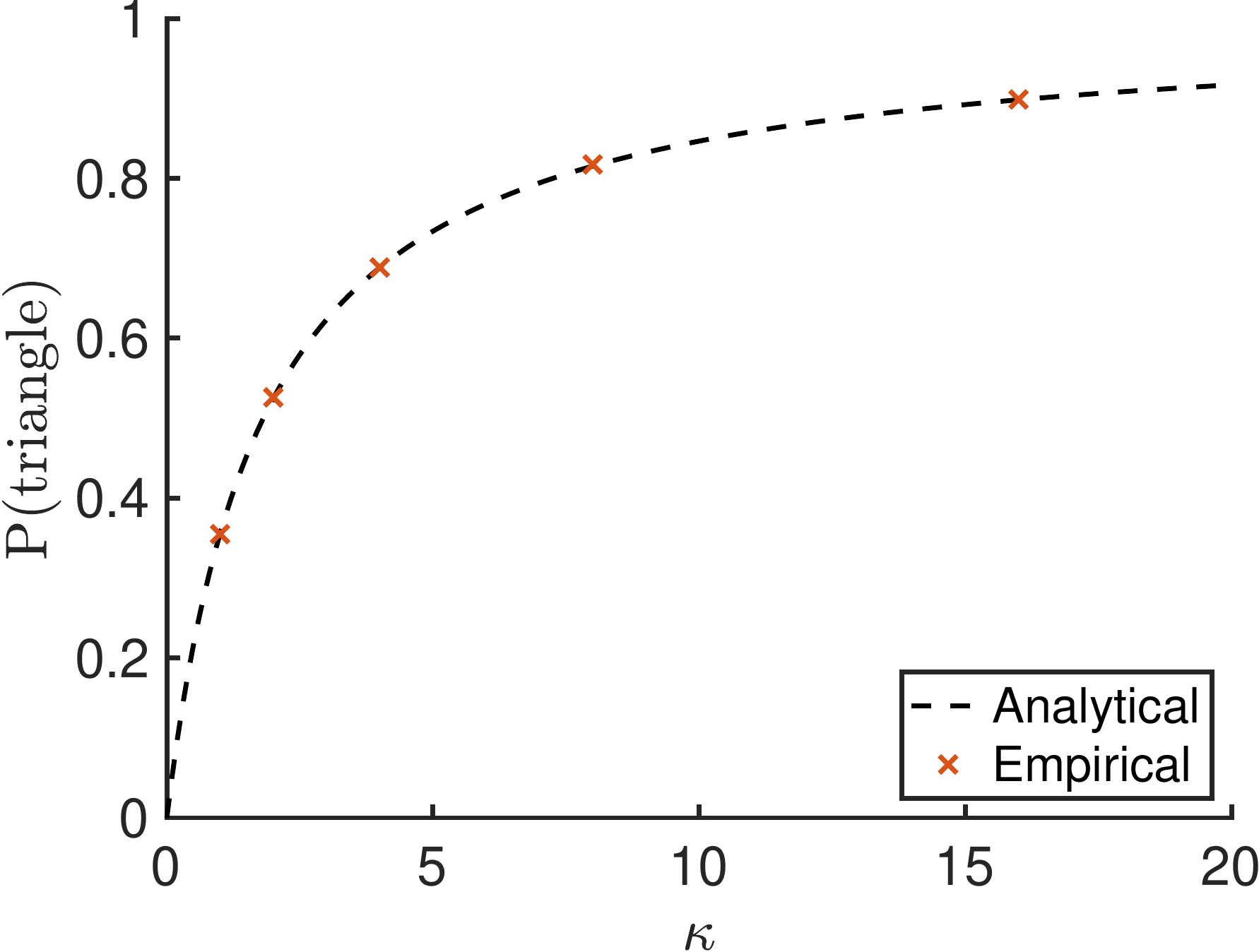}\hfill
\includegraphics[height=3.7cm]{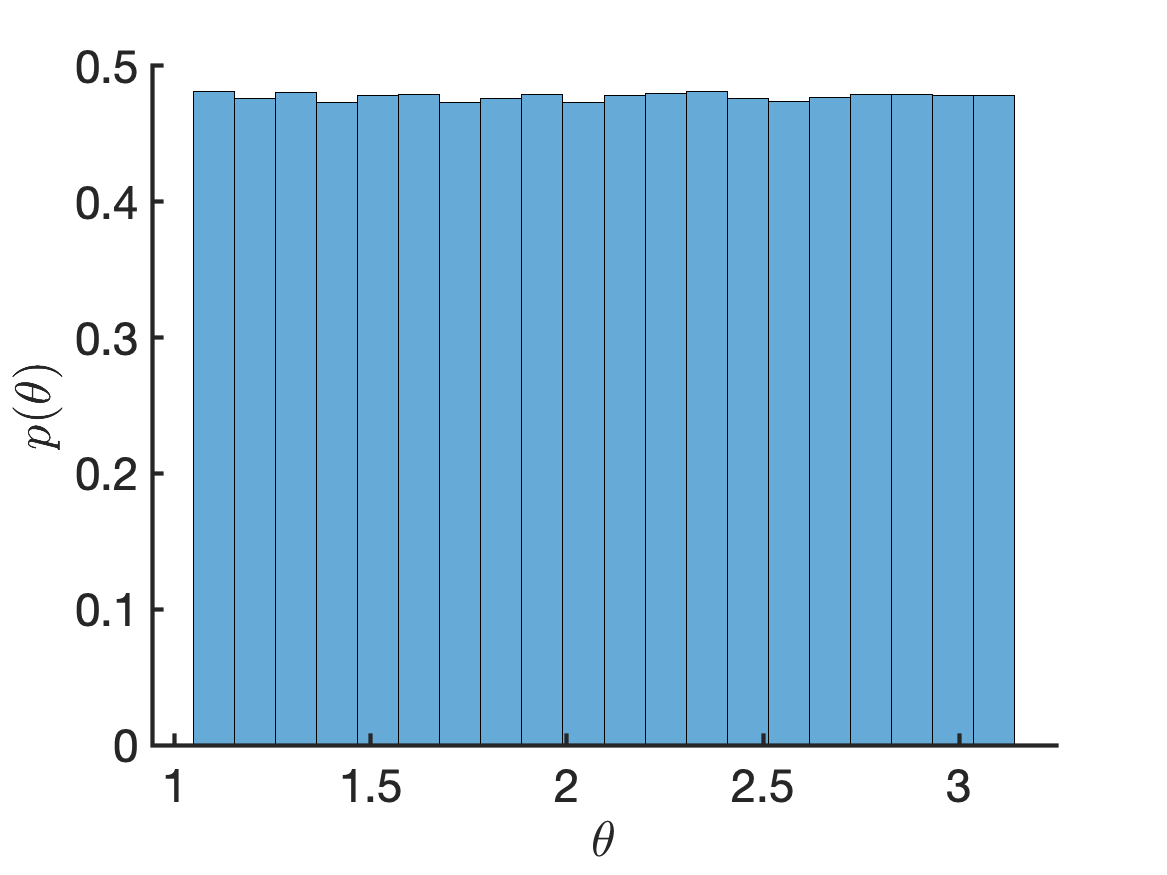}
\end{minipage}
\begin{minipage}{0.3\textwidth}
\centering
\begin{tabular}{ccc}
\multicolumn{3}{c}{Rejections}\\\hline
$\kappa$ & Gain & Lose\\\hline
1 &    $8\times 10^{-5}$ & 0.378  \\
2 &    0.28 & 0.10  \\
4 &   0.63 & 0.09    \\
8 &   0.82 & 0.08    \\
16 &   0.91 & 0.08  \\\hline
\end{tabular}
\end{minipage}
\caption{Statistics for Example \ref{ex:trimer} with probability measure \eqref{rhotrimer}. Left: probability of being in the triangle as a function of sticky parameter $\kappa$; analytically calculated in \eqref{ptriangle} and empirically calculated from a sampling run of $10^6$ points.  
Middle: distribution of internal angle $\theta$ conditional on being in the polymer, estimated from a sampling run with $10^7$ points and $\kappa=1$. The distribution is flat, as expected. It is zero for values of $\theta$ less than $\pi/3$ where the discs overlap. 
Right: Rejection rates for these sampling runs. For Gain moves (triangle $\to$ polymer) and Lose moves (polymer$\to$ triangle) all rejections came from the Metropolis step; we hypothesis that the rejection rate for Gain moves increases with $\kappa$ because we chose $\lambda\gain$ independently of $\kappa$. For Same moves, the rejection rates were uniformly 0.05 on the triangle (all from the NES failing), and 0.44 on the polymer (0.25 from NES; 0.12 from Metropolis step; 0.07 from violating inequalities.) 
Parameters for all sampling runs were  ${\sigma{=}0.5}$, ${\sigma\bdy{=}0.4}$, ${\sigtan{=}0.3}$, ${\lambda\lose{=}0.7}$, ${\lambda\gain{=}0.28}$, and runs of $10^6$ points took about 25 seconds on a 3.6 GHz Intel Core i7 iMac. 
}\label{fig:trimer2}
\end{figure}

We return to the example in Section \ref{sec2:stickysetup} of three sticky two-dimensional discs, with unit diameter. We find it convenient for numerical calculations to redefine the distance functions to be 
\[
q_{12}(x) = |x_1-x_2|^2-1,\quad q_{13}(x) = |x_1-x_3|^2-1,\quad q_{23}(x) = |x_2-x_3|^2-1\,.
\]
As before $\mathcal Q = \{q_{12},q_{23},q_{13}\}$. 
Instead of sampling the entire configuration space $S$ as in \eqref{Sdiscs1}, suppose we wish to sample only two configurations: 
a polymer whose internal angle can change, and a triangle, which is rigid. The polymer has labels $I$ with $\ieq = \{12,23\}$, $\iin= \{13\}$, and the triangle has labels $J$ with $\jeq = \{12,23,13\}$, $\jin = \emptyset$.  
The complete set of labels is $\mathcal I = \{I,J\}$  and the stratification is $S = M_I\cup M_J$. 
Notice that when constraint $q_{13}$ is dropped from $J$ to get $I$, it becomes an inequality, so $I$ is a one-sided Gain neighbour of $J$. 
One configuration from each of these manifolds is shown in Figure \ref{fig:setupexample}, where the polymer is $M_I = M_{12,23}$ and the triangle is $M_J = M_{12,23,31}$. 

We suppose that all bonds have the same sticky parameter $\kappa$, and sample the probability distribution constructed as in \eqref{rhos3},  
\begin{equation}\label{rhotrimer}
\rho_{\rm trimer}(dx) = Z^{-1}\left( \kappa^2|Q_I(x)|^{-1}\mu_I(dx) + \kappa^3 |Q_J(x)|^{-1}\mu_J(dx) \right).
\end{equation}
Recall that $Q_I,Q_J$ are defined in \eqref{Qsticky}, and  $Z$ as always is the normalizing constant. 

We keep track of the probability of finding the system in the triangle, for different sticky parameters $\kappa$, and the probability distribution for the internal angle of the polymer. Both probabilities can be calculated analytically so we may verify the sampler is working correctly. Figure \ref{fig:trimer2} verifies 
the empirical angle distribution for the polymer is flat, as expected from analytical calculations that consider a trimer with distance constraints but no momentum constraints \cite{vanKampen:1981ia}. This figure also shows that the empirical probability of finding the system in the triangle increases with $\kappa$, which is expected since the triangle has more bonds than the polymer. 

The probability of the triangle may be computed analytically by integrating $\rho_{\rm trimer}(dx)$ over the space of rotations and internal motions for each of the polymer and triangle. The triangle has no internal motions, so its probability is 
\[
\frac{P_{\rm Triangle}}{2\pi Z^{-1}} = \;\kappa^3\:2\:|Q_J(x)|^{-1}|x_{\rm cm}| = \frac{4\sqrt{3}}{9}\kappa^3 \,.
\]
Here $x_{\rm cm}$ is obtained from a configuration $x$ of the triangle by translating it so its center of mass lies at the origin, 
the factor $2\pi|x_{\rm cm}|$ comes from integrating over the triangle's rotations, and the factor 2 comes from counting the two different copies of the triangle obtained by permuting particle labels. For the unit triangle, $|x_{\rm cm}|=1$, $|Q_J| = \sqrt{27/4}$. 

The polymer has one internal degree of freedom, which may be parameterized by its internal angle $\theta$. Following \cite{HolmesCerfon:2016ji}, we adopt the parameterization 
$
x(\theta) = (-\sin \frac{\theta}{2}, -\frac{1}{3}\cos \frac{\theta}{2}, 0, \frac{2}{3}\cos \frac{\theta}{2}, \sin \frac{\theta}{2}, -\frac{1}{3}\cos \frac{\theta}{2}),
$
which keeps the center of mass of the polymer at the origin and doesn't rotate it. 
The surface measure in this parameterization is $\mu_J(x(d\theta)) = \big| \dd{x}{\theta}\big|d\theta$. 
The probability of the polymer is
\[
\frac{P_{\rm Polymer}}{2\pi Z^{-1}} = \kappa^2 \int_{\theta=\pi/3}^{5\pi/3} |Q_I(x(\theta))|^{-1}|x(\theta)| \Big| \dd{x}{\theta}\Big|d\theta  = \frac{4\pi}{9}\kappa^2.
\]
To evaluate the integral we calculated $|Q_I(x(\theta))| = \sqrt{4-\cos^2\theta}$, $|x(\theta)| = \sqrt{\frac{2}{3}+\frac{4}{3}\sin^2\frac{\theta}{2}}$, 
$\big| \dd{x}{\theta}\big| = \frac{1}{2}\sqrt{\frac{2}{3}+\frac{4}{3}\cos^2\frac{\theta}{2}}$; the integrand therefore is $\frac{1}{3}$.

Putting these calculations together shows the probability of the triangle is 
\begin{equation}\label{ptriangle}
P_{\rm Triangle} = \frac{\kappa }{\kappa + \pi/\sqrt{3}},
\end{equation}
which Figure \ref{fig:trimer2} shows agrees with our estimates from the sampled points. 

\subsection{A polymer of 6 sticky unit spheres in three dimensions}\label{ex:polymer}

\begin{figure}
\centering
\includegraphics[trim=0 3.5cm 0 0,clip,width=\textwidth]{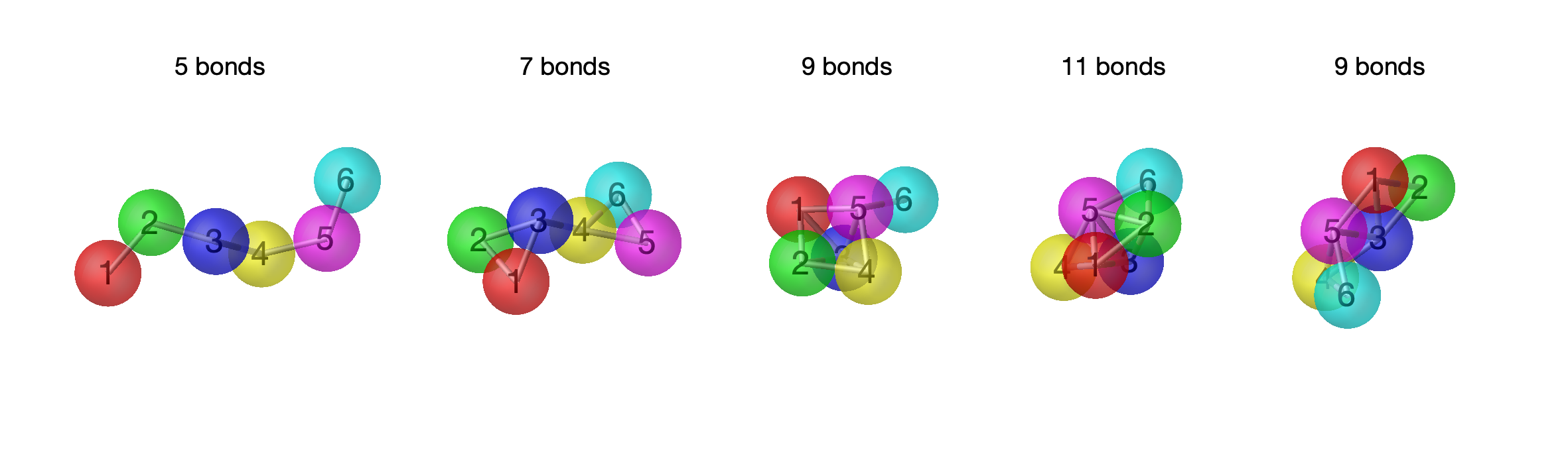}
\caption{Selected configurations of a polymer of $N=6$ unit spheres, from Example \ref{ex:polymer}. Bonds in the backbone of the polymer are fixed and cannot break, while all other bonds can break and form. 
}\label{fig:polymer1}
\end{figure}

\begin{figure}
\centering
\includegraphics[width=0.45\textwidth]{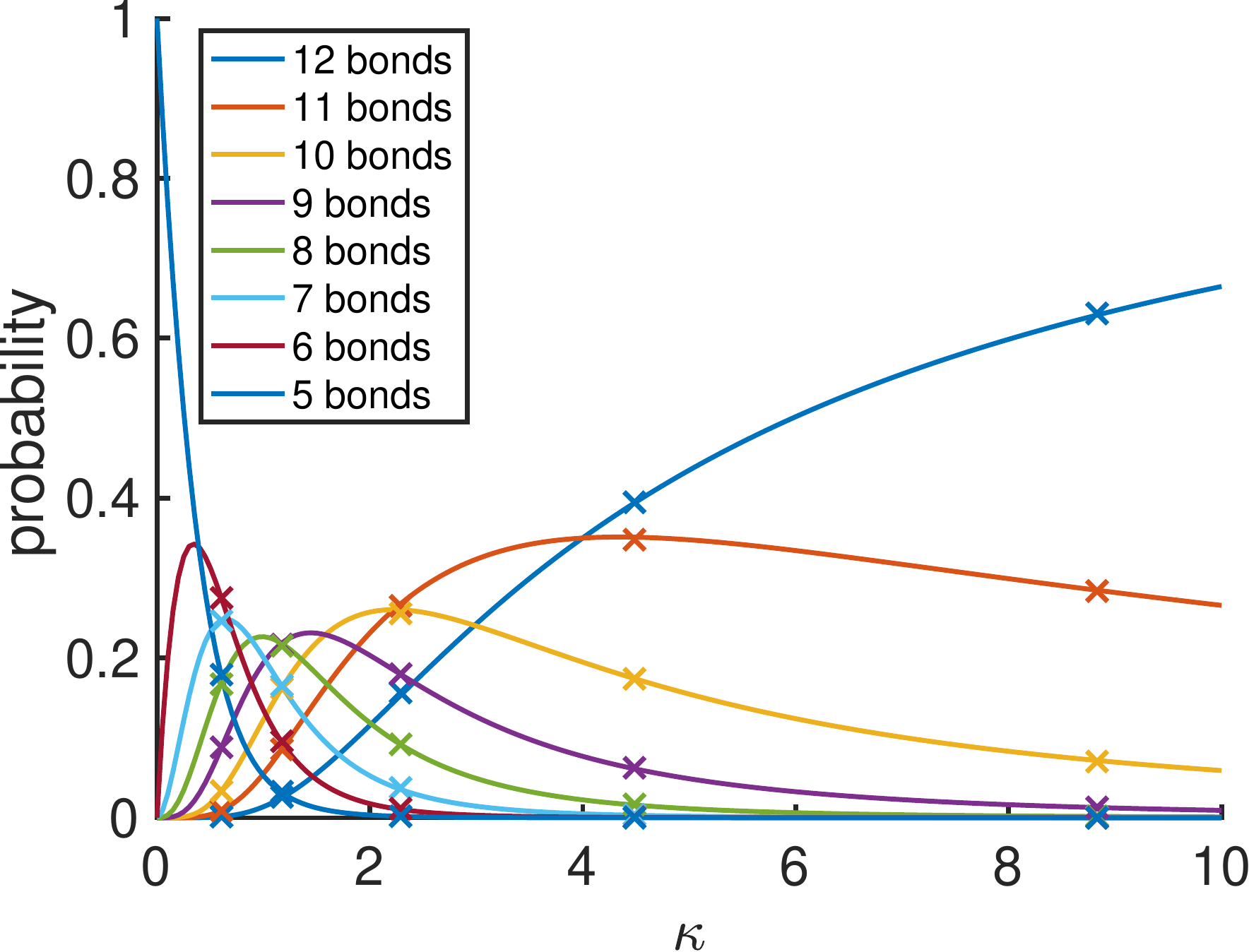}
\caption{
For Example \ref{ex:polymer}, the probability of finding a polymer of 6 unit spheres with each given number of bonds as a function of the sticky parameter $\kappa$,  the strength of the interactions between off-backbone spheres. Solid lines show the probabilities inferred by reweighting data from the Stratification Sampler at $\kappa_0=2$ with $4\times 10^7$ points (recording data every 4 points), and sampling parameters ${\sigma{=}0.4}$, ${\sigma\bdy{=}0.3}$, ${\sigtan{=}0.2}$, ${\lambda\lose{=}0.4}$, ${\lambda\gain{=}0.24}$. 
Markers show the probabilities estimated from Brownian dynamics simulations run to simulation time $10^4$ at different values of $\kappa$, using a Morse potential for the off-backbone interactions and a spring potential for the backbone interactions as in \eqref{Upolymer}, \eqref{BD}. 
The agreement with the Stratification sampler is excellent, verifying 
that the Stratification sampler gives good predictions even for a system whose interactions have small but nonzero range. 
}\label{fig:polymer2}
\end{figure}

\begin{figure}
\centering
\includegraphics[trim=2cm 2cm 2cm 2cm,clip,width=0.22\textwidth]{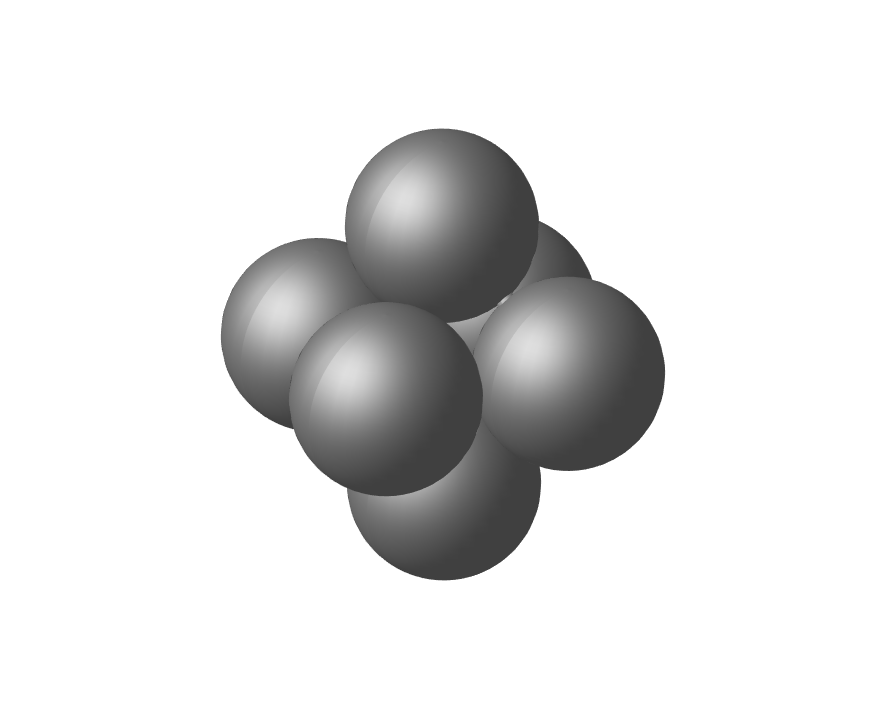}\qquad
\includegraphics[trim=2cm 1cm 2cm 2cm,clip,width=0.2\textwidth]{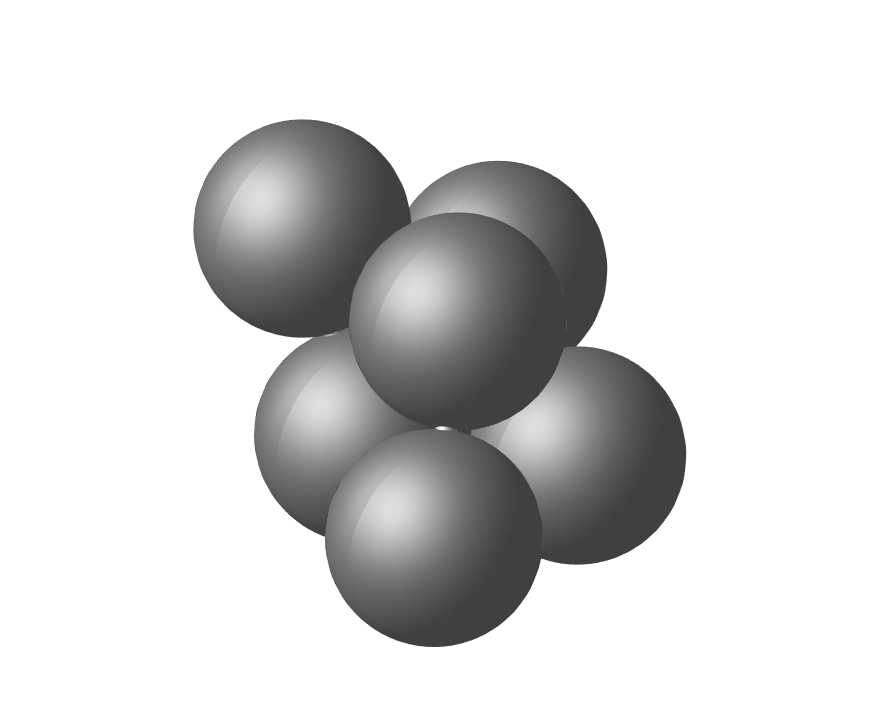}\qquad\quad
\includegraphics[trim=2cm 3.5cm 2cm 2cm,clip,width=0.3\textwidth]{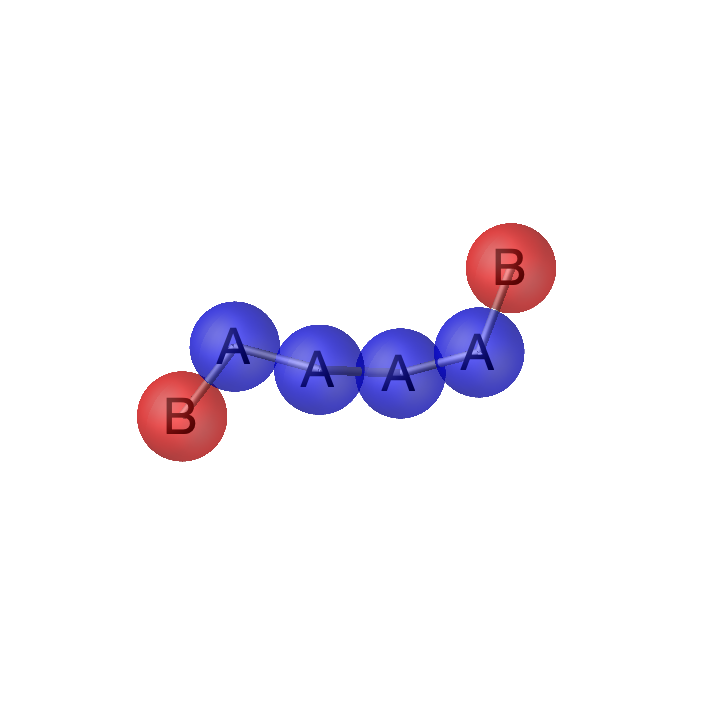}
\caption{
In Example \ref{ex:polymer} with identical sticky parameter for all bonds, the relative probabilities of forming the octahedron (left)  and polytetrahedron (middle) are 5\% and 95\% respectively, independent of sticky parameter.  For particles of types A and B as shown in the image on the right, we show that by choosing the sticky parameters $\kappa_{AA}, \kappa_{AB}, \kappa_{BB}$ for corresponding interacting pairs  appropriately, we can make the octahedron form with nearly 100\% probability. 
}\label{fig:polymerSA}
\end{figure}

\begin{figure}
\centering
\includegraphics[width=0.4\textwidth]{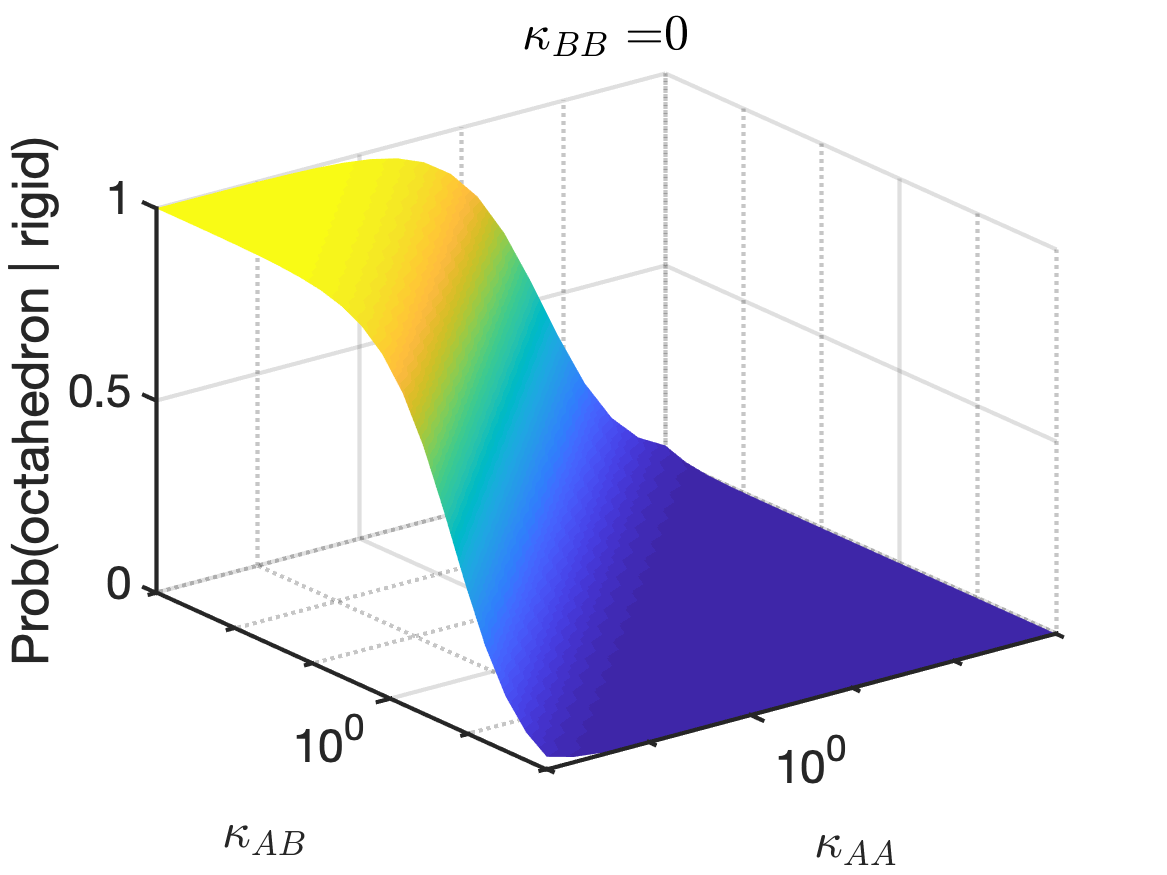}
\includegraphics[width=0.4\textwidth]{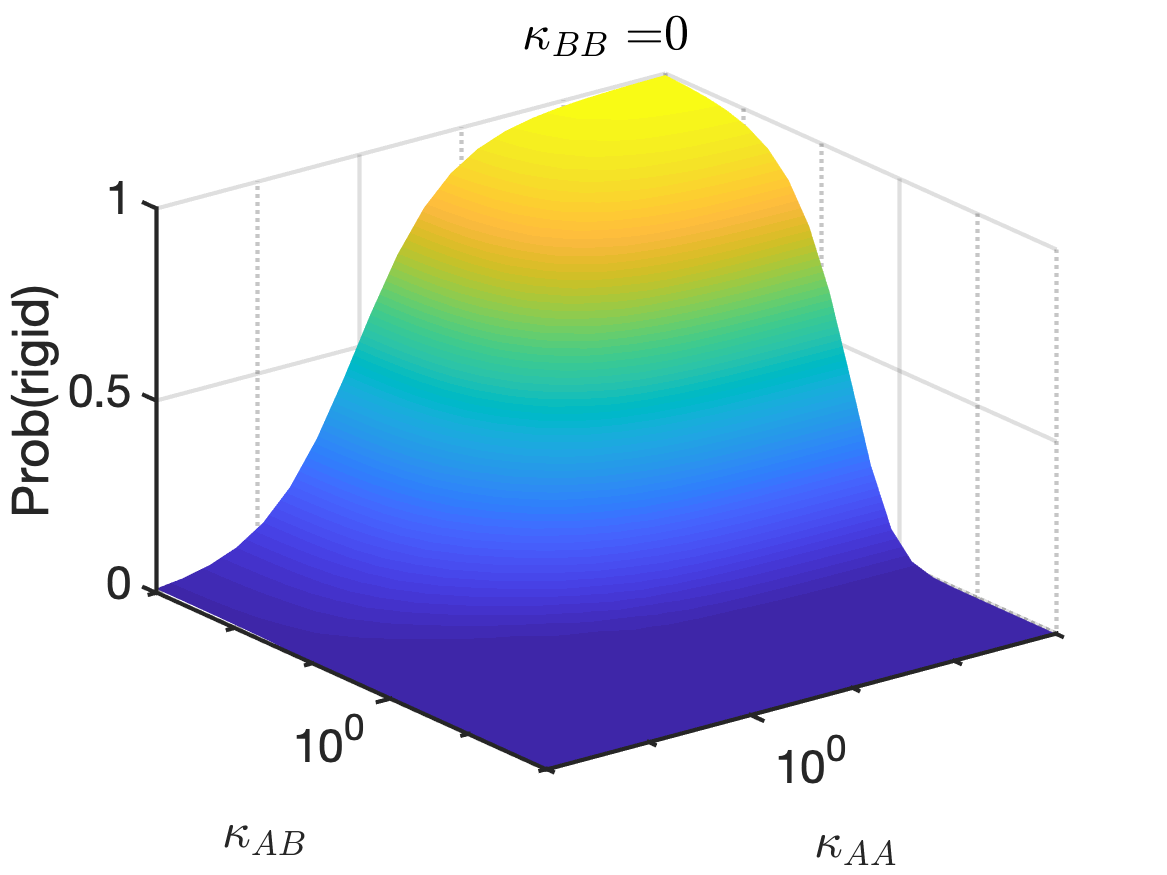}
\includegraphics[width=0.4\textwidth]{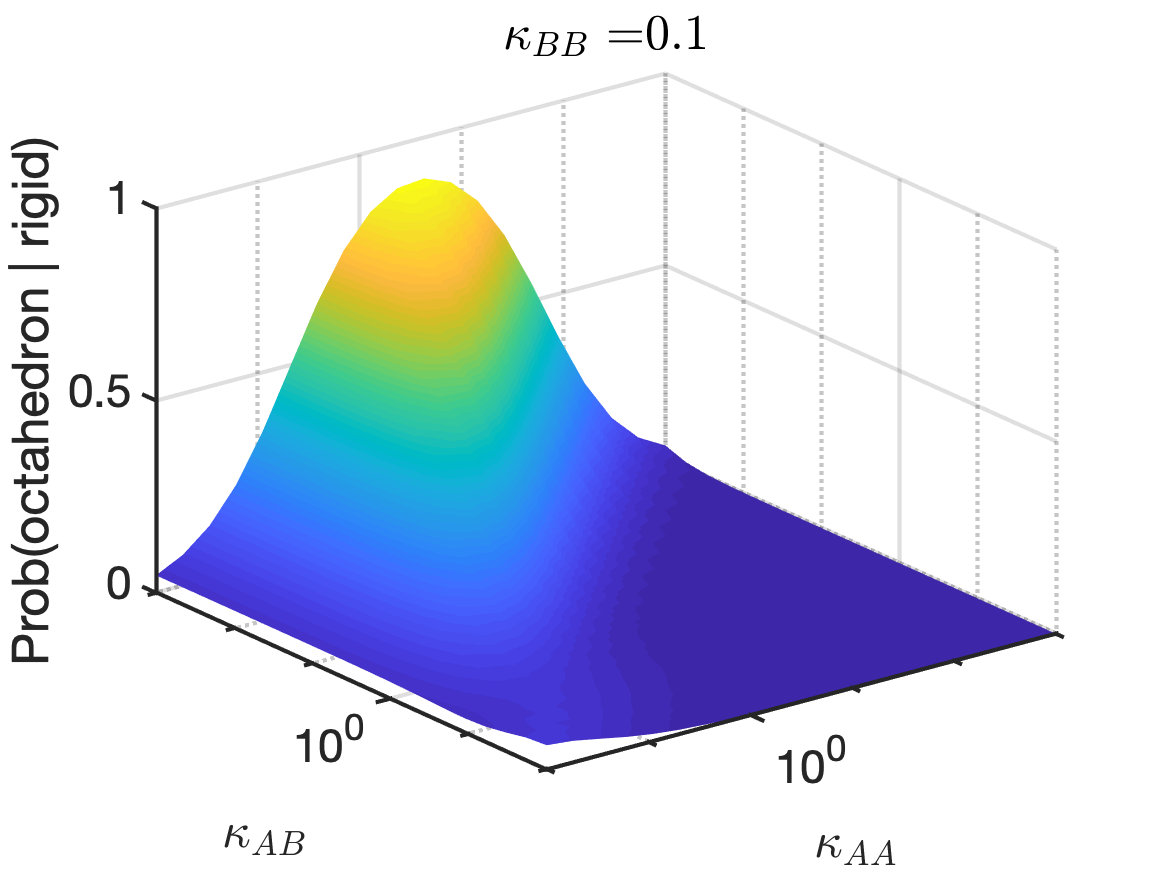}
\includegraphics[width=0.4\textwidth]{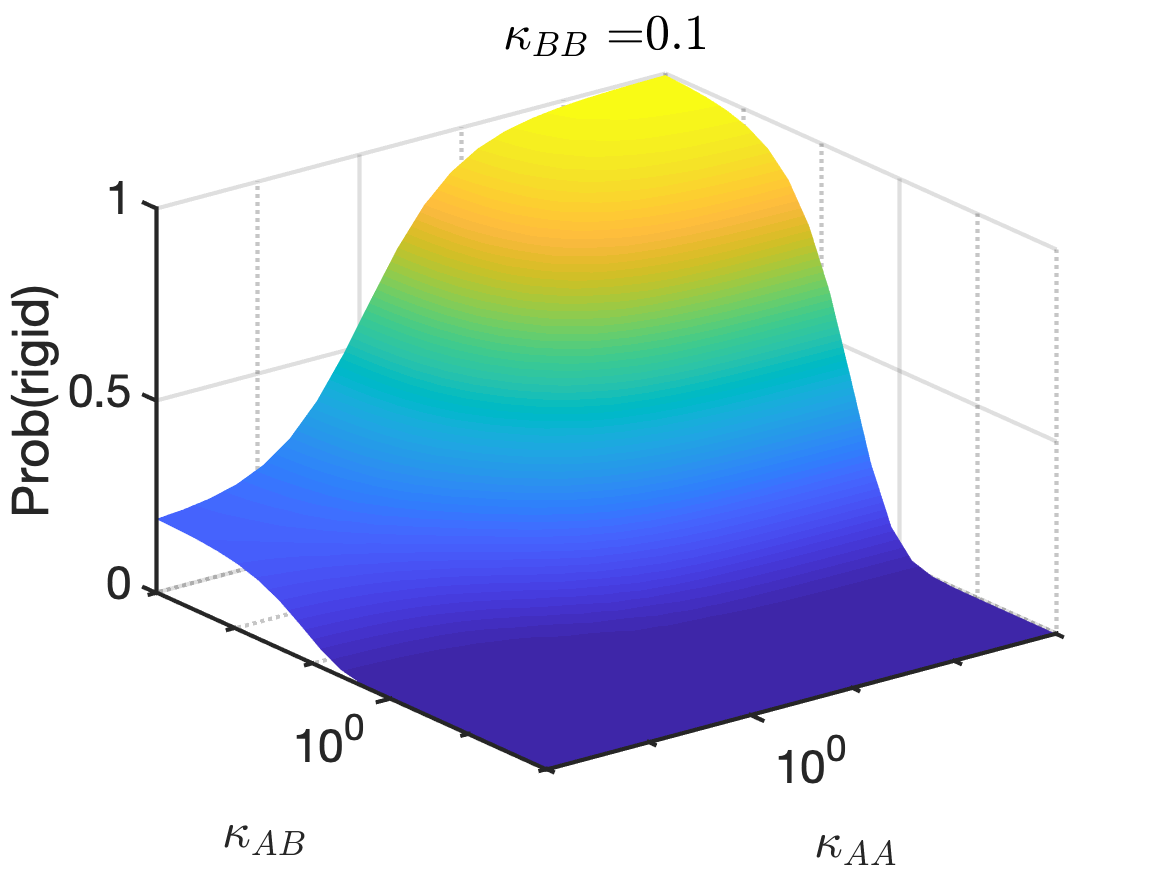}
\caption{
Studying inverse problems in interaction design, for $N=6$ unit spheres with two types of particles A, B as shown in Figure \ref{fig:polymerSA} and explained in Example \ref{ex:polymer}. Left column: equilibrium probability of observing an octahedron, given that the system is a rigid cluster (a cluster with $m=12$ bonds), for varying sticky parameters $\kappa_{AA}, \kappa_{AB}$, and for $\kappa_{BB}=0,0.1$. 
Right column: equilibrium probability of observing a rigid cluster as a function of $\kappa_{AA}, \kappa_{AB}$, for $\kappa_{BB}=0$ (top), $0.1$ (bottom.)  The top row shows the octahedron forms much more often than the polytetrahedron with small $\kappa_{AA},\kappa_{AB}$ and large $\kappa_{BB}$ (top left), however the probability of forming a rigid cluster is small in this limit (top right). The bottom row shows the octahedron  form with high probability when $\kappa_{BB}$ is small,  $\kappa_{AB} \sim O(1)$, and $\kappa_{AA}$ is large (bottom left), and furthermore, the probability of forming a rigid cluster is reasonably high in this limit (bottom right.) 
}\label{fig:SA}
\end{figure}

Next we consider a collection $N=6$ sticky three-dimensional spheres with unit diameter. We will study how to choose the sticky parameters to favour a desired ground state, and we will compare our sampler to Brownian dynamics simulation data, to show that the sampler can approximate a system with a noninfinitesimal potential. Designing interactions for small systems such as these has arisen as a challenge in materials science, since states with high symmetry can form the building blocks for meta-materials but entropy strongly disfavours symmetry \cite{Fan:2010jn, Meng:2010gsa,Miskin:2016bj}.

Suppose the spheres are connected in a chain, forming a polymer, and the bonds in the backbone of the polymer cannot break. We form a stratification by letting $\mathcal Q = \{q_{ij}(x)\}_{i\neq j}$ and considering all possible ways to let the non-backbone pairs come into contact: 
\[
\mathcal I = \{I: \{12,23,34,45,56\} \subset \ieq, \iin =\{ij:\;ij\neq \ieq\}\}.
\]
Some of the labels in $\mathcal I$ are infeasible, in the sense that $M_I = \emptyset$; for example, it is impossible to find a configuration where all pairs of spheres are in contact. 
The minimum number of bonds of a feasible configuration in $\mathcal I$ is 5, and the maximum number is 12. 
Some example configurations are shown in Figure \ref{fig:polymer1}. 

\subsubsection{Identical sticky parameters}
When all bonds have equal sticky parameter $\kappa$ the equilibrium distribution is 
\begin{equation}\label{rhoequal}
\rho_{\rm equal}(dx) = Z^{-1} \sum_{I\in \mathcal I}\kappa^{|\ieq|-5}|Q_I(x)|^{-1}\mu_I(dx). 
\end{equation}
We sampled $\rho_{\rm equal}$ with sticky parameter $\kappa_0=2$ for all bonds, recording the number of bonds at each time step. By reweighting this data (see \eqref{weights} below), we estimate the probability  $p_m(\kappa)$ of observing $m$ bonds at any sticky parameter $\kappa$. 
Figure \ref{fig:polymer2} shows that states with 5 bonds are most probable at low $\kappa$, while states with 12 bonds are most probable at large $\kappa$. 

The states with 12 bonds are the lowest-energy clusters, hence the ground states for this system. There are two distinct such clusters, after lumping together clusters with isomorphic adjacency matrices -- an octahedron and a polytetrahedron, illustrated in Figure \ref{fig:polymerSA}. These occurred with relative probabilities 5\%, 95\% respectively. It was shown in \cite{Meng:2010gsa} that the octahedron occurs much less frequently because it is has many more symmetries in its point group.\footnote{The numbers in our data are slightly different from those in Meng et al \cite{Meng:2010gsa}, because we consider the backbone to be fixed, while they allow all bonds to break and form, which changes the symmetry number for each cluster.}
Changing $\kappa$ does not change the relative probabilities of these clusters, because the equilibrium distribution only depends on the total number of bonds in a cluster.



\subsubsection{Designing the sticky parameters to self-assemble}
A problem that arises when studying the self-assembly of colloids is to form a particular ground state with high probability \cite{Hormoz:2011ir,Zeravcic:2014it}.
If we wish to form the octahedron, we must let the spheres have non-identical interactions. 
Suppose spheres 2,3,4,5 are type ``A'', and spheres 1,6 are type ``B'', with sticky parameters for AA, AB, BB interactions equal to $\kappa_{AA}, \kappa_{AB}, \kappa_{BB}$ respectively, as illustrated in Figure \ref{fig:polymerSA}. The equilibrium distribution is
\begin{equation}\label{rhoAB}
\rho_{AB}(dx) = Z^{-1} \sum_{I\in \mathcal I}\kappa_{AA}^{n_{AA}(I)}\kappa_{AB}^{n_{AB}(I)}\kappa_{BB}^{n_{BB}(I)}|Q_I(x)|^{-1}\mu(dx). 
\end{equation}
Here $n_{AA}(I),n_{AB}(I),n_{BB}(I)$ are the numbers of each kind of interaction present in $I$, respectively. 


We calculated the relative probability of the octahedron for different choices of parameters $\kappa_{AA}, \kappa_{AB},\kappa_{BB}$ by reweighting our data from a sampling run with identical sticky parameter $\kappa_0=2$. That is, instead of giving each sampled point $X_k,I_k$ equal weight when calculating probabilities or averages, we compute a weighted average using weights
\begin{equation}\label{weights}
w(X_k,I_k) = \left(\frac{\kappa_{AA}}{\kappa_0}\right)^{n_{AA}(I_k)}\left(\frac{\kappa_{AB}}{\kappa_0}\right)^{n_{AB}(I_k)}
\left(\frac{\kappa_{BB}}{\kappa_0}\right)^{n_{BB}(I_k)}. 
\end{equation}
Figure \ref{fig:SA} shows that when $\kappa_{BB} = 0$, the octahedron forms with 100\% probability as $\kappa_{AA}\to 0$ and $\kappa_{AB}\to \infty$. However, the probability of observing a rigid cluster goes to 0 in this limit, so the overall probability of observing the octahedron also goes to zero for these parameter choices. 
This observation is consistent with  \cite{Hormoz:2011ir}, which used a mean-field approach to argue that the best assembly occurs when sticky parameters are either zero or the same constant value; there are no intermediate maxima for the probability. 

In practice, it is hard to create sticky parameters that are exactly zero, especially for DNA-coated colloids where there is always some weak binding between non-complementary DNA strands \cite{Huntley:2016dn}. 
Therefore we also consider $\kappa_{BB} = 0.1$, and find the octahedron forms with 100\% probability for $\kappa_{AA} = O(1)$, and $\kappa_{AB}\to \infty$. In this limit the probability of observing a rigid cluster also goes to 100\%. Therefore, it is possible to form the octahedron with high probability, with the right choice of parameters. 

A tentative explanation for what makes the octahedron have high probability is that with 2 Bs, the polymer can be folded into an octahedron such that each B particle is in contact with 4 A particles, but neither B particle is in contact with itself. On the other hand, if the polymer folds into the polytetrahedron, then if the 2 Bs are each in contact with 4 As, then they must also be in contact with each other, forming a weak bond that wastes a precious A-A contact. 
We explored other labellings of the particles, and tentatively found that with any ordering of 2 nonadjacent Bs and 4 As, the octahedron formed with high probability under the same conditions, however with 3 Bs and 3 As we could not form the octahedron with high probability. 
It would be interesting to explore these observations further, and to find a general principle for the how the particle types and locations should be chosen for a polymer of $N$ spheres to fold into a given cluster.

\subsubsection{Comparing to Brownian dynamics simulations}
We now compare our data to Brownian dynamics simulations to verify that our algorithm is a good model for a real system of particles, which is never exactly at the sticky limit.  
We construct an energy $U(x)$ as 
\begin{equation}\label{Upolymer}
U(x) = \sum_{\substack{j=i+1}}U^{\rm spring}_{ij}(|x_i-x_j|) + \sum_{\substack{j\neq i+1}} U^{\rm morse}_{ij}(|x_i-x_j|)
\end{equation} 
where 
\[
U^{\rm spring}_{ij}(r) = \frac{1}{2}k_{\rm spring}(r-d_{ij})^2, \qquad 
U^{\rm morse}_{ij}(r) = E(1-e^{-\rho(r-d_{ij})})^2 - E
\]
are a spring potential, and Morse potential, respectively. The quantity $d_{ij}=1$ is the distance where spheres are exactly touching. The spring potential keeps the spheres in the backbone bonded, 
and the Morse potential creates a strong bond between non-backbone pairs when the distance between their surfaces is about $2.5/\rho$. 
The Boltzmann distribution  $Z^{-1}e^{-U(x)/k_BT}$ approaches $\rho_{\rm equal}$ in \eqref{rhoequal} as $k_{\rm spring}, \rho,E\to \infty$.

We numerically simulated the Brownian dynamics equations at temperature $T=1$,
\begin{equation}\label{BD}
\dd{X}{t} = -\grad U(X_t) + \sqrt{2}\:\eta(t),
\end{equation}
where $X_t\in\R^{6\cdot 3}$ is the configuration of the system and $\eta(t)\in \R^{6\cdot 3}$ is a white noise, using an Euler-Maruyama method with time step $\Delta t$.\footnote{The Euler-Maruyama method constructs a numerical approximation $X_0,X_1,X_2,\ldots$ to the solution to \eqref{BD} at times $0,\Delta t, 2\Delta t, \ldots$ as
\[
X_{k+1} = X_k + -\grad U(X_k)\Delta t + \sqrt{2 \Delta t} \:\xi_k, 
\]
where $\xi_1,\xi_2,\ldots \sim N(0,1)$ is a sequence of independent standard normal random variables. } 
We chose width parameters $\rho=60$, $k_{\rm spring} = 6\rho^2$ which gave us roughly the same width for the spring potential as for the Morse potential, and varied $E$. 
The choice of range is characteristic of certain DNA-mediated interactions (though some have smaller range \cite{Wang:2015ep}); it was small enough that the sticky limit gives a good description of the set of states on the energy landscape \cite{Trubiano:2019tl}, but large enough that is was not too prohibitive to resolve numerically. 
We needed a timestep of $\Delta t=10^{-6}$ to resolve the interactions.

We recorded the number of bonds in configuration $X_t$ every 0.05 units of time, identifying a pair $(i,j)$ as bonded if $|x_i-x_j| < 1+2.5/\rho$ \cite{Trubiano:2019tl}. Changing the cutoff (from about $2/\rho$ to $4/\rho$) did not significantly change the measured statistics. 
To compare to the Stratification sampler we determined the effective sticky parameter at each value of $E$ as
\begin{equation}
\kappa(E) = \int_0^{1+2.5/\rho}e^{-U(r;E)}dr.
\end{equation}
We calculated the integral by numerical integration, as Laplace asymptotics were not accurate enough for small values of $E$. 

The empirical probabilities $p_m(\kappa)$ estimated from the Brownian dynamics simulations at several values of $\kappa$ are shown in Figure \ref{fig:polymer2}. These agree with the probabilities predicted from the Stratification sampler. This not only verifies the correctness of the sampler, but it also shows that the sampler gives a good approximation to a real system with a non-delta-function potential.

\subsubsection{Efficiency compared to Brownian dynamics}

\begin{table}
\begin{tabular}{ccccc}
\hline
\# bonds & Stratification sampler & Brownian Dynamics & Ratio of Errors & Ratio of Physical Times \\
 $m$ & $\hat p_m^S\pm\sigma_m^S$  & $\hat p_m^{BD}\pm\sigma_m^{BD}$ & $\sigma_m^{BD}/\sigma_m^{S}$  & for fixed error \\\hline
12 & 0.1541 $\pm$ 0.00054 & 0.1568 $\pm$ 0.00147 & 2.7 & 162 \\
11 & 0.2675 $\pm$ 0.00048 & 0.2651 $\pm$ 0.00259 & 5.4 & 640 \\
10 & 0.2598 $\pm$ 0.00021 & 0.2556 $\pm$ 0.00152 & 7.2 & 1155 \\
9 & 0.1799 $\pm$ 0.00046 & 0.1798 $\pm$ 0.00184 & 4.0 & 358 \\
8 & 0.0916 $\pm$ 0.00039 & 0.0919 $\pm$ 0.00146 & 3.7 & 304 \\
7 & 0.0352 $\pm$ 0.00028 & 0.0374 $\pm$ 0.00130 & 4.6 & 467 \\
6 &0.0102 $\pm$ 0.00015 & 0.0115 $\pm$ 0.00038 & 2.6 & 147 \\
5 & 0.00177 $\pm$ 0.00006 &0.00197 $\pm$ 0.00017 & 2.7 & 158 \\\hline
\end{tabular}
\caption{Estimates of $p_m$ computed using the Stratification sampler (same parameters as in Figure \ref{fig:polymer2}) and by Brownian dynamics, both with $\kappa=2.885$. The Stratification sampler ran for $10^7$ steps and took a physical time of 13.6 minutes on a 2.40GHz Intel Xeon CPU E5-2680. The Brownian dynamics simulation ran for $10^4$ simulation time units and took a physical time of 5 hours on the same computer.  
The standard deviations $\sigma_m^S, \sigma_m^{BD}$ for each estimate $\hat p_m$ were calculated by binning the data into 8 bins; the standard deviation gives a measure of the error of the estimate. The ratio of errors $\sigma_m^{BD}/\sigma_m^S$  is shown in the second-last column, however it is a misleading comparison since the simulations took different amounts of physical time. By multiplying the squared ratio of errors by the ratio of physical simulation times, we obtain the ratio of physical simulation times that would be required for each method to produce the same given error, shown in the last column. 
}\label{tbl:polymer}
\end{table}


   
    
%

Our implementation of the stratification sampler, using dense arithmetic and standard matrix factorization, should be fast for small systems (not for large ones, which would require sparse arithmetic and more sophisticated matrix factorizations.) To see how much faster it is than Brownian dynamics for this example, we ran both methods with sticky parameter $\kappa = 2.885$, and computed the estimates $\hat p_m^S,\hat p_m^{BD}$ for $p_m(2.885)$, estimating the standard deviations $\sigma_m^S,\sigma_m^{BD}$ by binning the time series into 8 bins. The Stratification sampler was run for $10^7$ steps and the Brownian dynamics simulation was run for $10^4$ simulation time units. Table \ref{tbl:polymer} shows that under these conditions, the standard deviations (a measure of error) for the Stratification sampler are about 2.5-7 times smaller than those for Brownian dynamics. 

The standard deviation is not a good measure of efficiency, since the Stratification Sampler took much less physical time to run: it took 815 seconds (just under 14 minutes), compared to 17900 seconds (just under 5 hours) for Brownian Dynamics on the same processor in the same programming language and style -- about 22 times less time. 
Assuming that the simulation time scales with the inverse of the standard deviation squared, as for a typical Monte Carlo simulation, we can estimate the ratio of the physical times it would take each method to run to achieve the same value of standard deviation as $(\sigma^{BD}_m/\sigma^S_m)^2\cdot (17900/815)$. 
These predicted ratios of physical times are shown in the last column of Table \ref{tbl:polymer}; they range from about 150 to over 1000. That is, the Stratification sampler is \emph{two to three orders of magnitude faster} than Brownian dynamics for this example, to achieve a given level of accuracy. Of course, one of the sampling methods in the introduction would be faster than Brownian dynamics, but if one is interested in the dynamics of the particles, and not just their stationary distribution, then Brownian dynamics is the only method available. 


\subsection{Polymer adsorbing to a surface}\label{ex:polywall}

\begin{figure}
\centering
\includegraphics[trim=0cm 0cm 0cm 0cm,clip,width=0.48\textwidth]{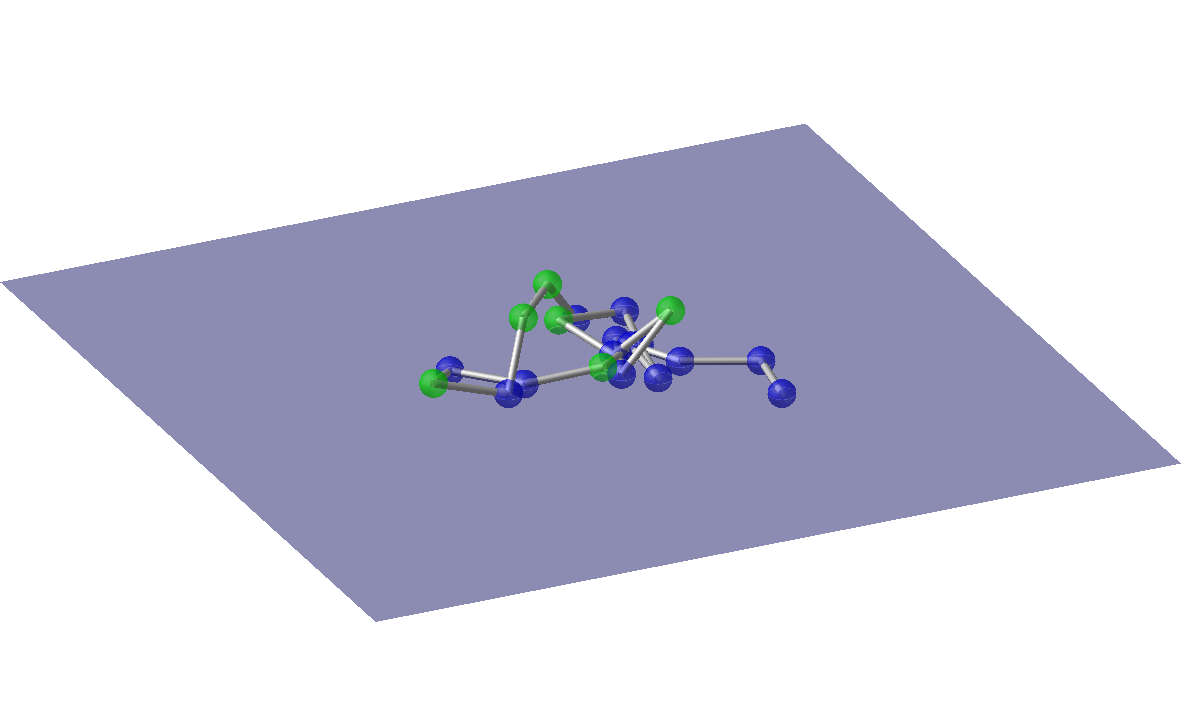}\quad
\includegraphics[trim=0cm 0cm 0cm 0cm,clip,width=0.48\textwidth]{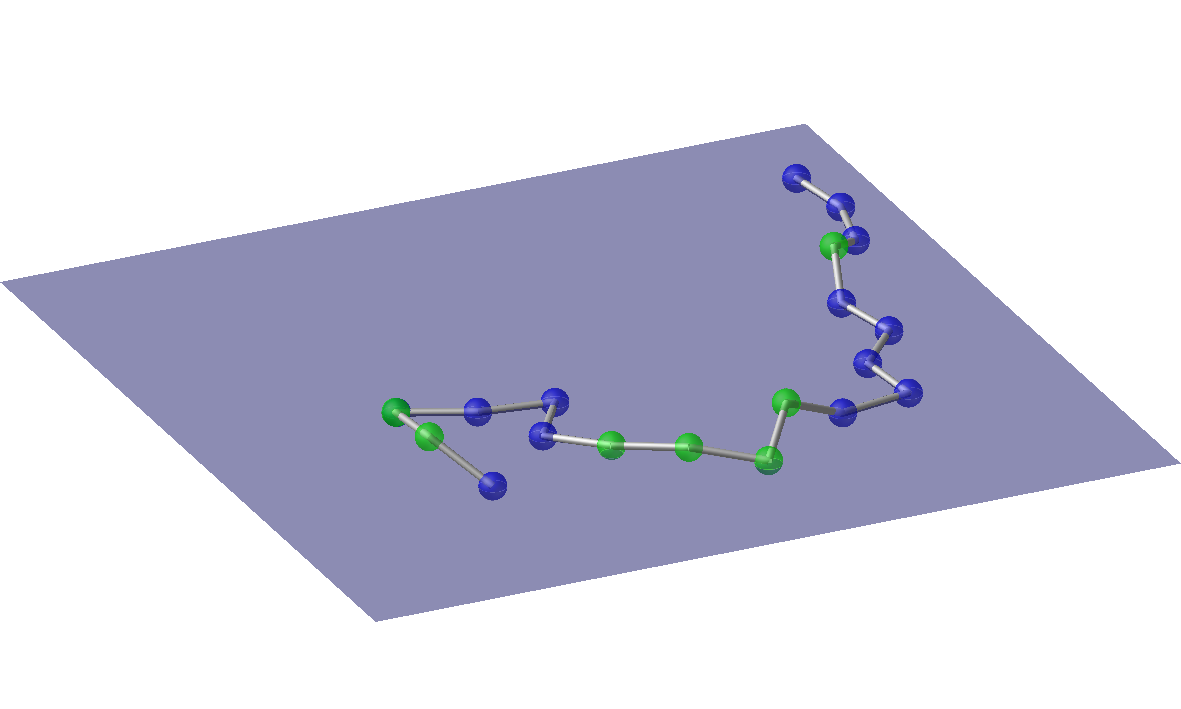}
\caption{Configurations drawn from the equilibrium distribution of a freely-jointed polymer (left) and a semiflexible polymer (right) with $N=20$ spheres, with sticky parameter $\kappa = 1$ with the wall and with bending stiffness parameter $k_{\rm bend} = 2$ in the semiflexible case. Blue spheres are in contact with the wall and green spheres are in free space. 
}\label{fig:polywall1}
\end{figure}

\begin{figure}
\centering
\includegraphics[width=0.32\textwidth]{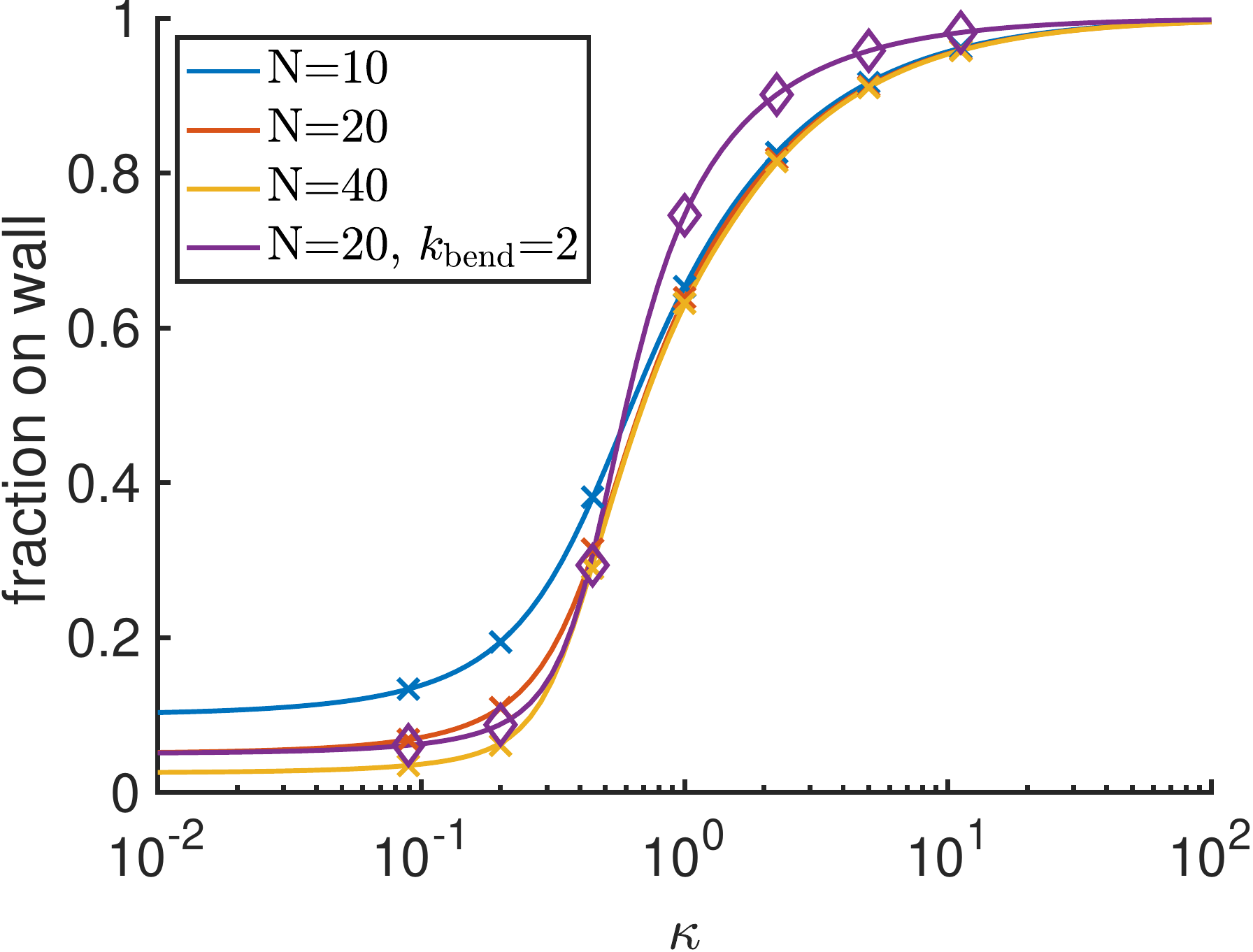}
\includegraphics[width=0.32\textwidth]{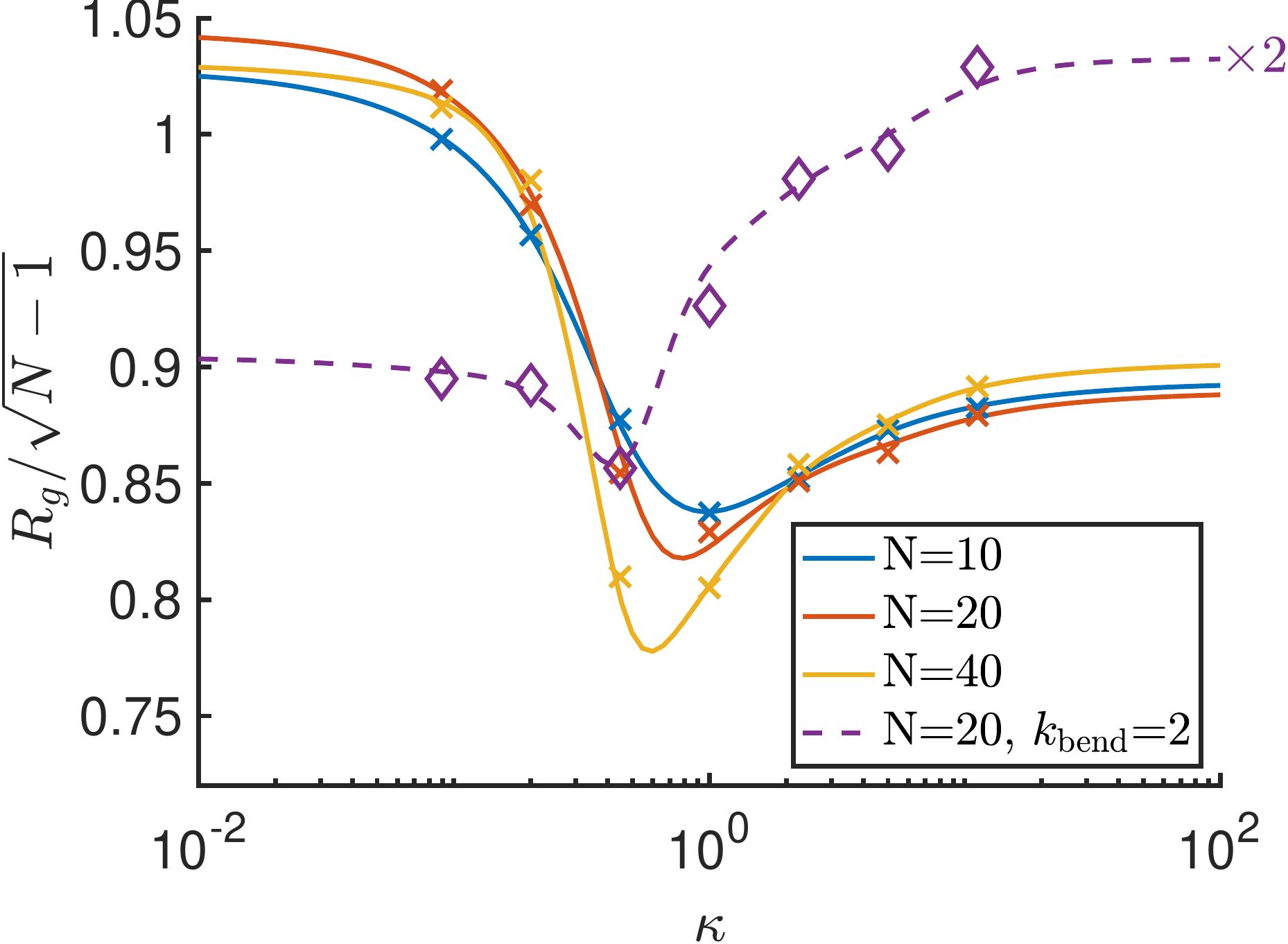}
\includegraphics[width=0.32\textwidth]{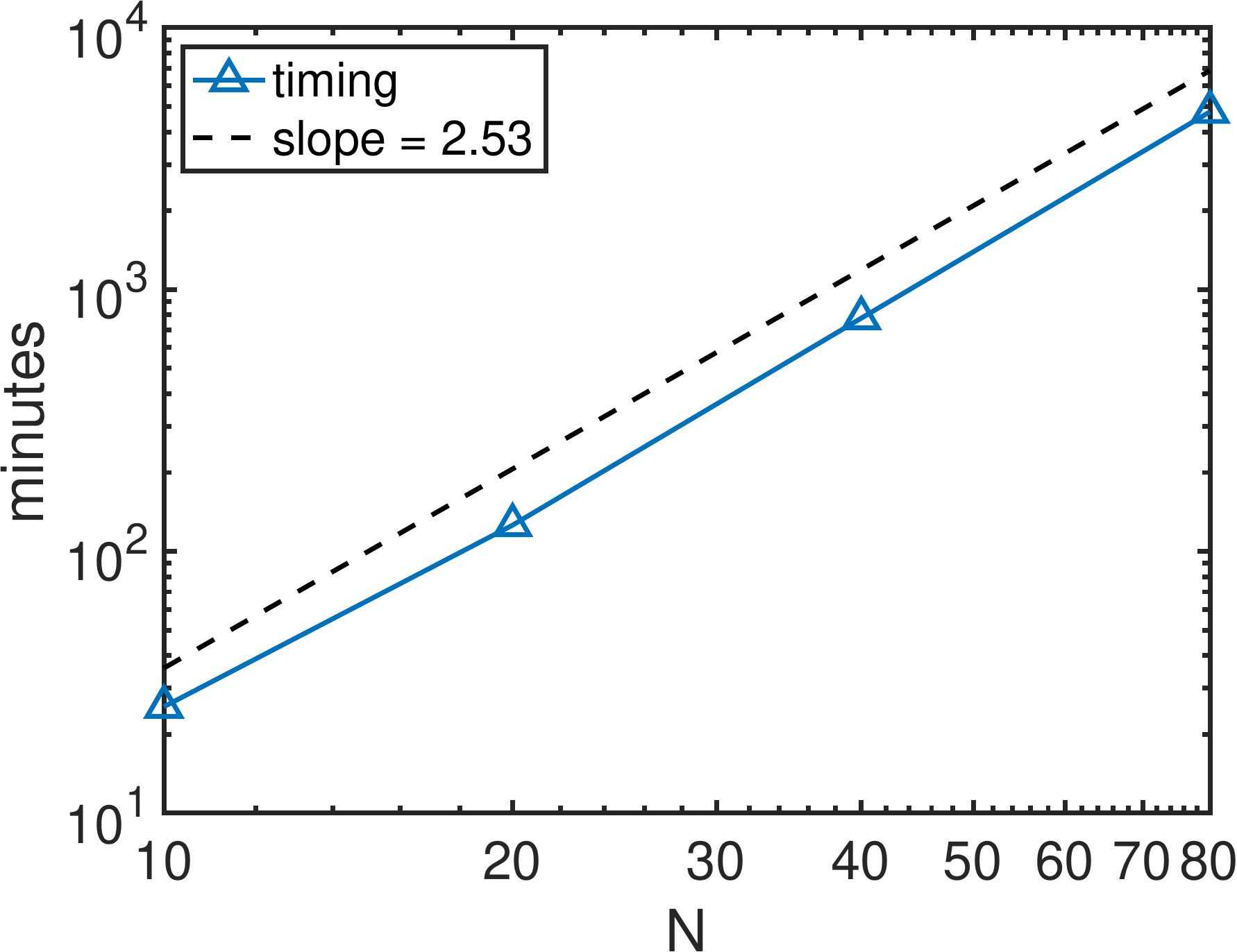}
\caption{
For Example \ref{ex:polywall}, the average fraction of particles on the wall (left) and average radius of gyration $R_g=|x_1-x_N|$ (middle) as a function of wall sticky parameter $\kappa$, for several values of $N$ as indicated in the legends, and for a semiflexible polymer with $N=20$ and $k_{\rm bend} = 2$ (see text.) For the semiflexible polymer we plot $R_g/2$ instead of $R_g$, to ease comparison. 
Markers indicate the values of $\kappa$ at which simulations were run, and curves indicate the averages estimated by reweighting data from nearby markers, in a similar way to \eqref{weights}. Each estimate was formed from $10^7$ simulation steps saved every 10 steps, with parameters $\sigma\bdy = 0.3$, $\sigtan = 0.2$, $\lambda\lose = 0.4$, $\lambda\gain = 2\sigma\bdy\lambda\lose$; for the freely jointed polymer we used $\sigma = 0.3, 0.2, 0.15, 0.12$ for $N=10,20,40,80$ respectively and for the semiflexible polymer we used $\sigma = 0.15$. 
The right plot shows the time in minutes taken by each simulation with $\kappa=1$ on a 2.40GHz Intel Xeon CPU E5-2680; on a log scale, the data is well-fit by a line with slope 2.5. 
}\label{fig:polywall2}
\end{figure}

A classic problem in polymer physics is to understand the statistical mechanics of a polymer that can adsorb weakly to a surface \cite{DeGennes:1981ho,Eisenriegler:1982il}. As the strength of interaction with the surface increases, the polymer behaves less like a three-dimensional polymer and more like a two-dimensional polymer (e.g. \cite{Milchev:2019db}, and references therein.) We can study this problem in the limit when the range of interaction with the surface is very short, so the particles are sticky on the surface. Consider a polymer formed from $N$ identical spheres with unit diameter that adsorb to a flat surface. Define a collection of functions $\mathcal Q = \mathcal Q_{\rm bonds} \cup \mathcal Q_{\rm surf}$, where the backbone distance functions are $\mathcal Q_{\rm bonds} = \{q_{i,i+1}(x)\}_{i=1\ldots N-1}$ with  $q_{i,i+1}(x) = |x_i-x_{i+1}|^2-1$, and the functions constraining a particle to the plane are $\mathcal Q_{\rm surf} = \{w_i(x)\}_{i=1\ldots N}$ with $w_i(x) = x_{i3}$ (the $z$-coordinate of particle $i$.) 
We form a stratification by considering all possible ways for the spheres to stick to the surface and never go through it, assuming the first sphere is always stuck:
\[
\mathcal I = \{I: \ieq=\mathcal Q_{\rm bonds}\cup I_{\rm surf} \; \text{where }w_1\in I_{\rm surf} \subset \mathcal Q_{\rm surf}, \; \; \iin = \mathcal Q_{\rm surf} - I_{\rm surf} \}. 
\]
For simplicity we don't consider inequalities between spheres in this example, though one could easily include them. 

We will study a both a freely-jointed polymer and a semiflexible polymer. A semiflexible polymer has a bending energy that penalizes deviations from a straight line: 
\begin{equation}\label{rhobend}
\rho_{\rm bend}(dx) = Z^{-1}e^{-\sum_{i=1}^{N-2}\frac{1}{2}k_{\rm bend}(1-\cos \theta_i)} \sum_{I\in\mathcal I} \kappa^{|I_{\rm surf}|}  |Q_I(x)|^{-1}\mu(dx). 
\end{equation}
Here $k_{\rm bend}$ is a parameter measuring the bending stiffness, and $\cos \theta_i = (x_{i+2}-x_{i+1})\cdot (x_{i+1}-x_i)$ is cosine of the $i$th internal angle $\theta_i$ in the polymer. 
A freely-jointed polymer has the same probability distribution but with $k_{\rm bend}=0$.


We sampled a freely-jointed polymer for $N=10,20,40,80$, and a semi-flexible polymer for $N=20$, each at seven different values of sticky parameter, $\kappa=5^{-3/2},5^{-2/2}, \ldots, 5^{3/2}$. For each $N$, we computed two statistics as a function of $\kappa$, the average fraction of spheres on the surface $f$, and the average radius of gyration, $R_g = |x_1-x_N|$, see Figure \ref{fig:polywall2}. For values of $\kappa$ not sampled, we estimated these averages by reweighting the data at the sampled values, constructing weights in a similar manner to \eqref{weights} and interpolating the estimates from each sampled value of $\kappa$. With the exception of $N\leq 10$, reweighting the data from a single simulation could not adequately capture the statistics over the full range of $\kappa$ values, because as $N$ increases, the distributions for each statistic become increasingly concentrated near their average values. 

As expected, the fraction $f$ of particles on the surface increases with $\kappa$, from close to $0$ at small $\kappa$, to nearly 1 at large $\kappa$. There is a sharp transition region whose width and location appear relatively independent of $N$, though is slightly sharper for the semiflexible polymer. The radius of gyration has interesting non-monotonic behaviour: for a freely jointed polymer, it appears to be slightly larger than $\sqrt{N-1}$ for small $\kappa$, and slightly smaller than $\sqrt{N-1}$ for larger $\kappa$, with a dip in between. That it is not always $\sqrt{N-1}$, its theoretical value for a polymer in free space of any dimension, must be due to the nonpenetration condition at the surface. 
The radius of gyration for a semiflexible polymer has the opposite behaviour, increasing with $\kappa$, but also nonmonotonically; overall it is nearly twice as large as for a freely jointed polymer. 

These observations are not new, but we present them here to show the variety of systems or constraints that can be studied by our sampler. An easy modification that could give new information would be to change the surface from a plane to a curved surface. 

We used this example to show that simulating larger systems will require implementing the stratification sampler using more sophisticated numerical linear algebra techniques. Figure \ref{fig:polywall2} compares the time it took to run each simulation for $10^7$ steps with $\kappa=1$. The increase in time is well-approximated by a power law $\propto N^{2.5}$. The reason for this poor scaling is that our implementation of the sampler performs several matrix factorizations in dense arithmetic at each step, such as the QR decomposition to find the tangent space, and the LU decomposition to solve systems of nonlinear equations using Newton's method, both of which scale as $O(N^3)$ \cite{Flannery:2007we,Trefethen:1997ee}.
These calculations can be done more efficiently, often with complexity $O(N)$, for problems with special structures, as we discuss in the conclusion, Section \ref{sec:conclusion}. Our implementation in dense arithmetic should be used only for systems with no more than a few dozen variables, where it does give a large speedup despite the poorer scaling with system size.



\subsection{Surface area of an $n$-dimensional ellipsoid}\label{ex:ellipse}

\begin{figure}
\centering
\includegraphics[trim=2.5cm 2.5cm 2.5cm 2.5cm,clip,width=0.5\textwidth]{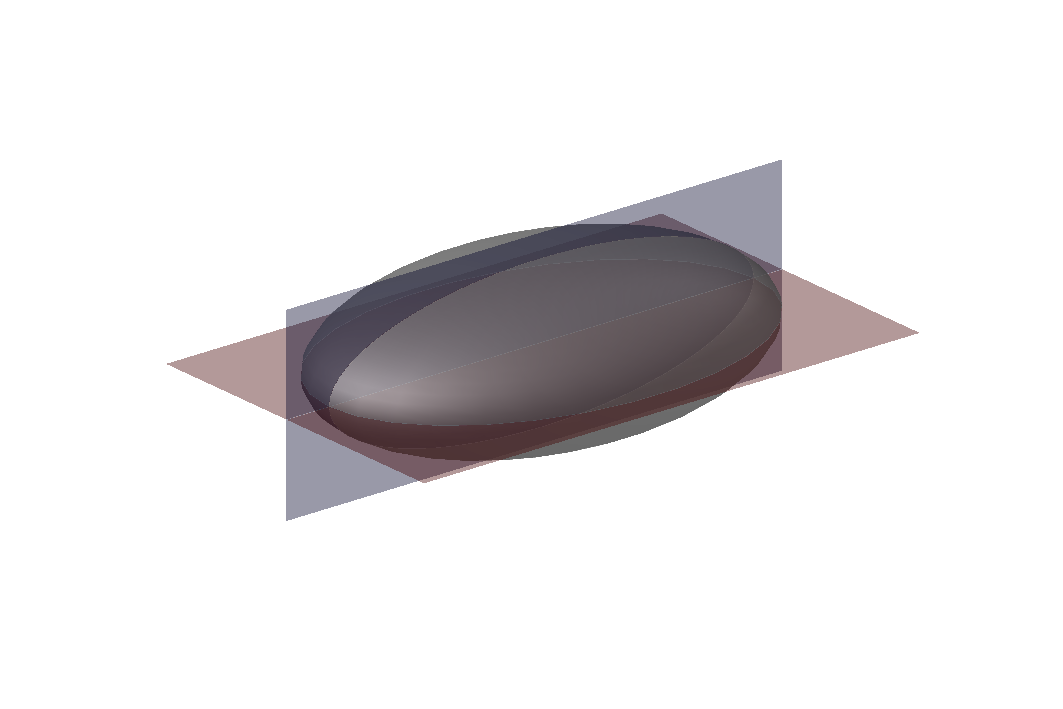}\quad
\caption{Example \ref{ex:ellipse}, with $n=3$. To estimate the surface area of the ellipse (grey), we form a stratification by intersecting it first with the blue plane, then with the red plane. The entire set of surfaces intersect in only two points, whose union is a 0-dimensional manifold. Since the volume of this 0-dimensional manifold is a known value (it is 2), we can infer the volumes of the other manifolds in the stratification from points sampled within the entire stratification. We apply this strategy in Example \ref{ex:ellipse} to calculate the surface area of a 10 dimensional ellipsoid. 
}\label{fig:ellipsoid}
\end{figure}

Our final example of a stratification illustrates how to use the Stratification Sampler to estimate the volume or surface area of a manifold. When the intrinsic dimension of the manifold is high enough (roughly larger than 4), deterministic methods cannot typically be used to compute its volume, because they require too many points to adequately cover the manifold  \cite{Simonovits:2003cq}. 
Consider the surface of an $n$-dimensional ellipsoid, defined as the solution to
\begin{equation}\label{qE}
E_n = \{ x\in \R^n: q_1(x) = \frac{x_1^2}{a_1^2} + \frac{x_2^2}{a_2^2} + \ldots + \frac{x_n^2}{a_n^2} -1 = 0\}. 
\end{equation}
Here $a_1,\ldots,a_n>0$ are the lengths of the semiaxes. 
Let $\text{Vol}(E_n)$ be the $n-1$-dimensional volume of $E_n$. 
While there exist expressions for the surface area of $E_n$, they are complicated (see Equation (10) in \cite{Rivin:2007gs}), so we show instead how one can estimate the surface area through sampling. 

The idea is to construct a stratification which contains at least one manifold whose volume we know, or can easily estimate. Since we can exactly calculate the volume of a 0-dimensional manifold, by counting the number of points in it, our strategy will be to add constraints until their intersection is a collection of isolated points.  The simplest constraint is a plane, so we intersect the ellipsoid with an ordered collection of planes going through the origin. 
To this end, let $\mathcal Q = \{q_1,q_2,\ldots,q_n\}$ where $q_1 = x_1^2/a_1^2 +  \ldots + x_n^2/a_n^2 -1$ is the constraint describing the ellipse, and $q_2(x) = x_2,q_3(x) = x_3,\ldots,q_n(x) = x_n$ are constraints describing the planes. 
Let $I\e{k} = (\ieq\e{k},\iin\e{k})$ with 
\[
\ieq\e{k} = \{q_1,q_2,\ldots,q_k\}, \qquad \iin\e{k} = \emptyset.
\]
Each $\ieq\e{k}$ is the ellipse $E_n$ intersected with the first $k-1$ planes.  
Our stratification is  
$
\mathcal I =  I\e{1} \cup I\e{2} \cup \cdots \cup I\e{n}. 
$
See Figure \ref{fig:ellipsoid} for an illustration of the stratification when $n=3$.
Notice that all Gain neighbours are two-sided since there are never any inequalities.  

We construct a probability distribution as 
\begin{equation}\label{rhoE}
\rho_E(dx) = Z^{-1}\sum_{k=1}^n c_k\mu_k(dx) \,,
\end{equation}
where $\mu_k = \mu_{I\e{k}}$ is the natural surface measure on each manifold $M_k = M_{I\e{k}}$, $c_k>0$ is a constant weight for each manifold, and $Z$ is the normalizing constant. 
A good strategy is to choose the weights $c_k$ so the sampler spends roughly the same amount of time in each manifold. Such a strategy has been shown to be optimal in certain systems when temperature is the variable that changes, rather than dimension \cite{Martinsson:2019df}. 

If we have a collection of points $X_1,X_2,\ldots, X_m \sim \rho_E$, then we can estimate the surface area of $E_n$ as 
\begin{equation}\label{VolEst}
\text{Vol}(E_n) \;\;\approx\;\; \frac{\text{\# of points in }M_1}{\text{\# of points in }M_n}\cdot \frac{2c_n}{c_1}
\end{equation}
This estimate comes from observing that, if the algorithm is ergodic, then 
\[ 
\frac{\sum_{i=1}^m1_{M_1}(X_i) }{ \sum_{i=1}^m1_{M_n}(X_i)} \to 
\frac{\int_{M_1}\rho(dx)}{ \int_{M_n}\rho(dx)} =
\frac{ \int_{M_1}c_1\mu_1(dx) }{ \int_{M_n}c_n\mu_n(dx)}
\] 
as $m\to \infty$. Since $\int_{M_1}c_1\mu_1(dx) = c_1\text{Vol}(E_n)$, and  $\int_{M_n}c_n\mu_n(dx) = 2c_n$ since the manifold $M_{n} = \{\pm a_1\}$ contains exactly two points, we solve for $\text{Vol}(E_n)$ to get the estimate above. 

We sampled the distribution \eqref{rhoE} for $n=10$, 
%
%
using semiaxes $(a_1,\ldots,a_{10}) = (2,2,2,2,3,3,3,1,1,1)$ 
and choosing the weights to be $c_k = e^{0.94 k}$. This choice was motivated by an initial sampling run with equal weights in which the probability to be in $I\e{k}$ was approximately proportional to $e^{-0.94 k}$. 
We used sampling parameters $\sigma = 0.6, \sigtan = 0.3, \sigma\bdy = 0.4, \lambda\lose = 0.4, \lambda\gain = \sigma\bdy\lambda\lose=0.16$. We ran the sampler for $10^7$ steps and used \eqref{VolEst} to estimate the surface area, dividing the data into 10 bins to estimate a one standard deviation error bar, and obtained 
\[
\text{Vol}(E_{10}) \approx 7155 \pm 162,
\] 
about a 2\% relative error. 

To check the accuracy of this estimate, we computed it using an alternative stratification formed from $\mathcal I = (\{-q_1\},\emptyset)\cup(\emptyset,\{-q_1\})$. This stratification consists of the surface of $E_{10}$, and its interior. The volume of the interior is $a_1a_2\cdots a_{10}\cdot \pi^5/120 $, so an estimate for the surface area with the same sampling parameters as above is
\[
\text{Vol}(E_{10}) \approx  \frac{\text{\# of points in $E_n$}}{\text{\# of points in interior($E_n$)}}\cdot \frac{\pi^5\cdot a_1a_2\cdots a_{10}}{120}
= 7138\pm 21,
\]
about a 0.2\% relative error.  While the second stratification clearly gives a more accurate estimate, because it is formed from fewer manifolds, one does not usually know the volume of a high-dimensional shape analytically, so this simpler method cannot usually be applied.

We remark that this strategy is a form of thermodynamic integration, but where \emph{dimension} is the variable that changes, rather than temperature \cite{Frenkel:1984ft}. For higher-dimensional volumes, where the sampler spends most of the time in the intermediate dimensions of the stratification, an alternative to choosing non-equal weights $c_k$ would be to break up the problem into several separate sampling problems: first estimate the ratio of the $d$-dimensional volumes in the stratification to the $d{-}1$-dimensional volumes, then the ratio of $d{-}2$-dimensional to $d{-}3$-dimensional, and so on down to the ratio of one-dimensional to 0-dimensional volumes. For a variant of this strategy
applied to shrinking sets, not changing dimensions, see \cite{Simonovits:2003cq, Lovasz:2006df}. 

\subsection{What goes wrong with isotropic Gain and Lose proposals}\label{ex:isotropic}

\begin{figure}
\centering
\begin{tikzpicture}[scale=1.6]
\def\ta{0.05};
\def\tb{0.25};
\def\tc{0.45};
\def\td{0.65};
\def\te{0.85};
\def\tf{1.05};
\def\tg{1.25};
\def\th{1.45};
\node[vertex,label=below:{x}] (x) at (0,0) {};
\draw[ultra thick] (0,1)-- (8.5,1);
\draw[|-|,dotted]  (-0.65,0) -- node[midway,left] {$k\sigma$} (-0.65,1);
\draw[->,dashed] (x) --({sin(\ta r)},{cos(\ta r)});
\draw[->,dashed] (x) --({sin(\tb r)},{cos(\tb r)});
\draw[->,dashed] (x) --({sin(\tc r)},{cos(\tc r)});
\draw[->,dashed] (x) --({sin(\td r)},{cos(\td r)});
\draw[->,dashed] (x) --({sin(\te r)},{cos(\te r)});
\draw[->,dashed] (x) --({sin(\tf r)},{cos(\tf r)});
\draw[->,dashed] (x) --({sin(\tg r)},{cos(\tg r)});
\draw[->,dashed] (x) --({sin(\th r)},{cos(\th r)});
\node[label={[align=center,font=\footnotesize]above:$a(y|x)$, k=1\\$a(y|x)$, k=4}] at (-0.65, 0.98) {};
\node[mystar,label={[align=center,font=\footnotesize]above:1\\1}] at ({tan(\ta r)}, 1) {};
\node[mystar,label={[align=center,font=\footnotesize]above:1\\1}] at ({tan(\tb r)}, 1) {};
\node[mystar,label={[align=center,font=\footnotesize]above:1\\1}] at ({tan(\tc r)}, 1) {};
\node[mystar,label={[align=center,font=\footnotesize]above:1\\0.67}] at ({tan(\td r)}, 1) {};
\node[mystar,label={[align=center,font=\footnotesize]above:1\\0.23}] at ({tan(\te r)}, 1) {};
\node[mystar,label={[align=center,font=\footnotesize]above:1\\0.012}] at ({tan(\tf r)}, 1) {};
\node[mystar,label={[align=center,font=\footnotesize]above:0.17\\$1.9\times 10^{-7}$}] at ({tan(\tg r)}, 1) {};
\node[mystar,label={[align=center,font=\footnotesize]above:$1.9\times 10^{-13}$\\0}] at ({tan(\th r)}, 1) {};
\end{tikzpicture}
\caption{Acceptance probabilities for isotropic proposals in Section \ref{ex:isotropic}, for a Lose move from $x=(0,-k\sigma)$ to the line $\{(x_1,x_2):x_2=0\}$. Here $\sigma$ is the standard deviation of the step size in the Gain proposal \eqref{pvsame}, and $k>0$ controls the distance to the boundary. 
We chose 8 directions for a Lose proposal with vertical angles equally spaced in $\theta\in[0.05,1.45]$ (dashed arrows), and for each proposed move (blue stars) we evaluated the acceptance probability in \eqref{acc} with $\lambda\lose=\lambda\gain$ and $k=1,4$, using the formula $a(y|x)=\min(1,e^{-\frac{k}{2\cos^2\theta}}\sqrt{2\pi}k/\cos^2\theta)$ to be derived in \eqref{a21}.
As $k$ increases, most directions give small acceptance probabilities, because they result in large displacements which are unlikely to have been proposed in a Gain move. Furthermore, such large moves are more likely to cause the nonlinear equation solver to fail, before even reaching the Metropolis step. Hence, Lose moves for this kind of proposal are accepted only if they jump to a nearby boundary, i.e. if $k$ is small and the proposed direction is nearly perpendicular to the boundary. However, for nearby boundaries the corresponding Gain move is small, so has only a small probability of being accepted during the Metropolis step. 
}\label{fig:BadMove}
\end{figure}
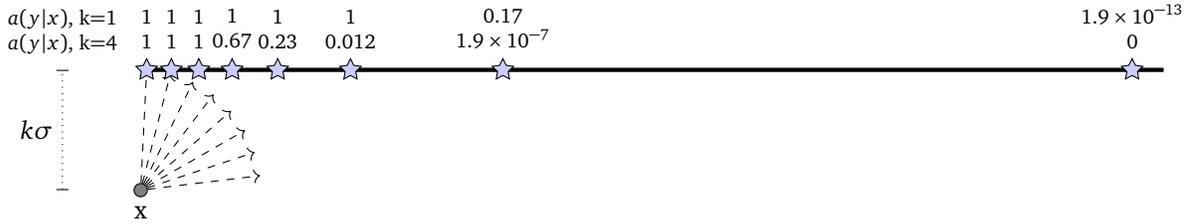

In this final example we show what happens when Gain and Lose tangent steps are chosen isotropically in the tangent space, instead of treating the directions normal to and tangential to the boundary separately. We modify our sampler to choose the tangent step $v$ for a Gain move the same way as for a Same move (recall \eqref{vSame}), and the tangent direction for a Lose move uniformly on the surface of the unit sphere. 

This method works in theory --  it samples the correct probability measure (once we adjust the proposal densities accordingly) --  but it is not efficient. 
The rejection rates for Gain and Lose moves are both high, and cannot be made small by an appropriate choice of parameters. 
For example, for Example \ref{ex:polymer} with $\kappa=2$, about 88\% of Lose moves and 93\% of Gain moves were rejected. These numbers were relatively insensitive to the sampling parameters (perturbing around $\lambda\gain{=}\lambda\lose{=}0.4$, $\sigma{=}0.3$), increasing somewhat for $\sigma\lessapprox 0.05$. 

The reason for the rejections gives a clue as to what is going wrong. 
Among the Lose proposals, a full 63\% were rejected while projecting back to the manifold: 26\% because the nonlinear solver failed to find a solution to \eqref{qsolveLose}, and 37\% because the solver found $\alpha<0$ so the step was not in direction $v$, but in $-v$. Another 20\% were rejected during the Metropolis step, 5\% during the inequality check and 0.05\% during the reverse projection check. 
Among the Gain proposals, 56\% were rejected during the inequality check, 27\% during the Metropolis step, 9\% because the nonlinear solver failed to find a solution to \eqref{qsolveSame};, and 0.5\% were rejected during the reverse projection check. 
Compare these statistics to those for Same moves, which are constructed in the same way as Gain moves, for which 19\%, 11\%, 6\% were rejected for the inequality, Metropolis, and nonlinear projection steps respectively.

The rejection rates for Gain and Lose moves are high because these proposals are  highly asymmetric with each other: by choosing a direction isotropically, Lose moves usually propose to jump to boundaries that are far away. Not only does this make the nonlinear solver more likely to fail (both because the boundary is far away, and because the direction proposed is not always consistent with the required sign of $\alpha$), but, even if the nonlinear solver is successful, such moves are unlikely to have been proposed in the Gain step, hence, lead to small Metropolis factors. This effect is illustrated in Figure \ref{fig:BadMove}. 
This is why most Lose moves are rejected in either the proposal step or the Metropolis step. 
Furthermore, about half of the Gain moves are rejected during the inequality check, because Gain moves don't distinguish between directions which move to the interior of the manifold, and directions which move away from it, where the inequality is violated. 
Proposing moves that are so unlikely to be accepted wastes computation; a sequence of smaller moves that are each rather likely to be accepted would give a more efficient sampler than one that proposes medium-sized moves that are unlikely to be accepted. 


\section{Details of the Stratification Sampler}\label{sec2:algorithm}

This section builds on Section \ref{sec2:overview}, giving the details of the Stratification Sampler algorithm that are necessary to implement it. 
The algorithm is summarized in pseudocode in Appendix \ref{sec:summary}, and codes are available on Github at \cite{gitfigures}. 

In what follows, the current point/label is $(x,I)$, and we wish to propose a new point/label  $(y,J)$. First we propose a label $J$ as described in Section \ref{sec2:label}, and then we propose a point $y$ on manifold $M_J$ as in Section \ref{sec2:point}. We give further details on the label proposals in Section \ref{sec2:label2} and on the point proposals in Section \ref{sec2:point2},  including the label and point probabilities $\lambda_{IJ}(x),h_{IJ}(y|x)$.

\subsection{Label proposal, revisited}\label{sec2:label2}

Recall from Section \ref{sec2:label} that we propose a Same, Gain, or Lose move with probabilities $\lambda\same, \lambda\gain, \lambda\lose$ respectively. A label $J$ is chosen uniformly from either the set of Gain neighbours $\ngain{I}$ or the set of nearby Lose neighbours $\nlosesig{I,x}$ (see \eqref{Nlosesig}.)
Letting $n\gain(I) = |\ngain{I}|$, $n\lose(I,x) = |\nlosesig{I,x}|$ be the number of each type of neighbour, the overall probabilities of each move type are  
\begin{gather}\label{LambdaAlg2}
\Lambda\gain(x) = \lambda\gain1_{n\gain>0}(x),\;\;
\Lambda\lose(x) = \lambda\lose1_{n\lose>0}(x), \;\;
 \Lambda\same(x) = 1-\Lambda\gain(x) - \Lambda\lose(x) . 
\end{gather}
The densities for the label proposals are
\begin{equation}
\text{Same:} \;\; \lambda_{II}(x) = \Lambda\same(x), \qquad
\text{Gain:} \;\; \lambda_{IJ}(x) = \frac{\Lambda\gain(x)}{n\gain(I)}, \qquad
\text{Lose:} \;\; \lambda_{IJ}(x) = \frac{\Lambda\lose(x)}{n\lose(I,x)}\,.
\end{equation}

It remains to explain how to calculate $\nlosesig{I,x}$. For this, we need to estimate the minimum distance $h_{\rm opt}=h_{\rm opt}(x,I,q)$ from a point $x$ on manifold $M_I$, to the boundary $q(x)=0$, if the distance is measured along manifold $M_I$. We do this by linearizing all constraints. Let $v$ be a unit vector in the  tangent space to $M_I$ at $x$ (i.e. any unit vector $v$ that satisfies $Q_I^T(x)v = 0$.) The distance $h$ in direction $v$ to $q(x)=0$ is estimated by linearizing the equation $q(x+hv)=0$ about $x$ to get\footnote{We should really linearize $q(x+hv+w)$, where $w$ is an (unknown) vector in the normal space to $M_I$ at $x$. The contribution from $w$ vanishes to linear order, since linearizing the equations $q_i(x+u)=0$, $i\in \ieq$ defining manifold $M_I$, gives $\sum_{i\in \ieq}u\cdot\grad q_i = 0$. Therefore any vector $u$ which maintains the constraints to linear order must lie in the tangent space to $M_I$ at $x$.} $q(x) + hv\cdot \grad q = 0$ , so  
\begin{equation}\label{hv}
h = -\frac{q(x)}{v\cdot \grad q(x)}\,.
\end{equation}
The minimum positive value of $h$ over all unit vectors $v\in \mathcal T_{x,I}$ occurs in direction 
\begin{equation}\label{vopt}
v_{\rm opt} = -\mbox{sgn}(q(x)) \frac{P_{x,I}\grad q(x)}{|P_{x,I}\grad q(x)|}, \qquad \text{where } \; P_{x,I} = T_{x,I}T_{x,I}^T
\end{equation}
is the orthogonal projection matrix onto $\mathcal T_{x,I}$. The minimum positive distance is then 
\begin{equation}\label{hopt}
h_{\rm opt} = \frac{|q(x)||T_{x,I}T_{x,I}^T\grad q|}{|T_{x,I}^T\grad q|^2} = \frac{|q(x)|}{|T^T_{x,I}\grad q|}\,.
\end{equation}
We simplified this expression using $|T_{x,I}T_{x,I}^T\grad q| = |T_{x,I}^T\grad q|$, since $T_{x,I}$ has orthonormal columns. When we need to clarify what the optimal direction and distance depend on we write $h=h_{\rm opt}(x,I,q)$, $v\opt=v\opt(x,I,q)$.

\subsection{Point proposal, revisited}\label{sec2:point2}

Given a proposal $J$, we then choose $y\in M_J$. 
We describe how to do this for each of the three move types, and we give the formula for proposal densities. We also explain how to calculate the proposal densities for the reverse move, $y\to x$, when all we are given is $x,y$; this calculation is required to evaluate the acceptance probability \eqref{acc}. 

\subsubsection{Same move: $x\in M_I$, $y\in M_I$}

\paragraph{Generating the move} 

Recall from Section \ref{sec2:point} that we choose a tangent step $v\in \mathcal T_{x,I}$; suppose its density is $p_{II}\esame(v;x)$. We take a step in direction $v$ then project back to $M_I$ by solving \eqref{qsolveSame} for the normal step $w \in \mathcal N_{x,I}$, and hence the proposal point $y=x+v+w$. By the regularity assumption, the columns of $Q_I(x)$ span $\mathcal N_{x,I}$, so we may write  $w=Q_I(x)a$  for $a\in \R^{m_I}$ and solve 
\begin{equation}\label{sameeqns}
q_i(x+v+Q_I(x)a) = 0, \qquad i\in \ieq\,,
\end{equation}
for $a$. 
Solving this system of nonlinear equations may be done using any nonlinear equation solver, provided it is deterministic with a  deterministic initial condition for reasons we explain momentarily. Let NES denote the particular choice of nonlinear equation solver used in the algorithm. Our implementation uses Newton's method as described the Appendix \ref{sec:summary}, Algorithm \ref{alg:takestep}, though there exist better choices. 

If a solution $a$ is successfully found by the NES, the proposal is constructed as $y = x+v+Q_I(x)a$. It must then be checked that $y$ satisfies the inequalities in $\ieq$; if not, the proposal is $y=x$. Additionally, if the solver fails to find a solution, the proposal is $y=x$. The solver could fail for a variety of reasons: there could be no solution for a particular $v$, or, even if there is a solution, the NES could simply fail to converge to a solution.

\paragraph{Density of the proposal } 
The density for a successful proposal that is found by the NES and satisfies the inequalities is $p_{II}\esame(v;x)\big|\pp{v}{y} \big|$, where $\big|\pp{v}{y} \big|$ is the inverse of the determinant of the Jacobian of the transformation from $v \to y$. This was shown in \cite{Zappa:2018jy} to be $\big|\pp{v}{y} \big| = |T_{x,I}^TT_{y,I}|$. 
For a plane, $\big|\pp{v}{y} \big| =1$, since $T_{x,I}=T_{y,I}$. 

Let $\mathcal A_{\rm same}(x)\subset M_I$ be set of values of $y$ that are \emph{accessible} from $x$ via a Same move: they can be found by solving \eqref{sameeqns} using our NES; by definition they satisfy the inequalities in $\iin$. In order for $\mathcal A_{\rm same}(x)$ to be well-defined, we must use a deterministic initial condition for the NES. In addition, there must be only one  $v$ for each move $x\to y$, a fact we verify shortly. Let $1_{\mathcal A_{\rm same}(x)}(y)$ be the characteristic function for the set $\mathcal A_{\rm same}(x)$, equal to 1 if $y\in\mathcal A_{\rm same}(x)$, and 0 otherwise. 
The proposal density for a successful proposal is
\begin{equation}\label{hsame}
h_{II}(y|x) = p_{II}\esame(v;x)|T_{x,I}^TT_{y,I}|1_{\mathcal A_{\rm same}(x)}(y)\,.
\end{equation}
If the tangent step is an isotropic Gaussian with step size $\sigma$ as in \eqref{vSame}, then 
\begin{equation}\label{pvsame}
p_{II}\esame(v;x) = \frac{1}{\sqrt{(2\pi\sigma^2)^{d_I}}}e^{-\frac{|v|^2}{2\sigma^2}}.
\end{equation}

The overall probability density for a Same move is 
\begin{equation}\label{SameDens}
\lambda_{IJ}(x)h_{IJ}(y|x)  = \Lambda\same(x)\frac{1}{\sqrt{(2\pi\sigma^2)^{d_I}}}e^{-\frac{|v|^2}{2\sigma^2}}|T_{x,I}^TT_{y,I}|1_{\mathcal A_{\rm same}(x)}(y)\,.
\end{equation}
Note that in practice we do \emph{not} have to calculate the Jacobian factor $|T_{x,I}^TT_{y,I}|$ for a Same move, since this factor cancels out in the acceptance ratio \eqref{acc}.

\paragraph{Evaluating the density for the reverse move}
Given $x,y\in M_I$, we now explain how to evaluate the proposal density \eqref{hsame}. We focus on evaluating the density from $x\to y$, to avoid rewriting \eqref{hsame}, but in practice one proposes a $y$ and then evaluates the density for the move from $y\to x$, and uses this density to calculate the acceptance probability in \eqref{acc}. 

The step $v$ can be found by projecting $y-x$ onto $\mathcal T_{x,I}$, as 
\begin{equation}\label{SameV}
v=T_{x,I}T_{x,I}^T(y-x)\,.
\end{equation}
From this we can evaluate $p_{II}(v;x)$. Notice that, as required, $v$ is uniquely defined from $x,y$. 

We also need to evaluate $1_{\mathcal A_{\rm same}(x)}(y)$. To do this, we must run the NES starting from $x$ with tangent step $v$, and check (i) the NES produces a solution, and (ii) this solution equals $y$. Even though we know that $y$ is a solution to \eqref{sameeqns}, because the pair $(x,y)$ was a successful proposal in a forward move, we still need to check that the NES actually produces this solution. If not, the density for $y$ must be set to zero. 

Therefore, the NES must be run twice for each move: once to produce a forward move, and once to compute the density of the reverse move. 
That this check is necessary for detailed balance was first pointed out in \cite{Zappa:2018jy}, and \cite{Lelievre:2019be} showed empirically that samplers can be more efficient for proposals where this check is critical to computing the correct acceptance probability.

\subsubsection{Gain move: $x\in M_I$, $y\in M_J$ with $J\in \mathcal N\gain(I)$} 

\paragraph{Generating the move} 

A Gain move is similar to a Same move but is on manifold $M_J$. 
Suppose that $q$ is the constraint dropped from $\ieq$ to obtain $\jeq$. 
We choose $v\in\mathcal T_{x,J}$ with density $p_{IJ}\egain(v;x)$. 
We take a step in direction $v$, then solve for $w\in\mathcal N_{x,J}$ such that $y=x+v+w$ lies on $M_J$. This is done by solving \eqref{sameeqns} but with constraints in $J$ instead of in $I$. 
If a proposal fails, either due to inequality constraints, or because the NES fails, then we set $y=x,J=I$. 

Recall we construct $v$ as in \eqref{GainStep} from a vector $u_n$ pointing away from the boundary, and the matrix $T_{x,I}$, whose columns are tangent to $M_I$ and normal to $u_n$. The unit vector $u_n\in \mathcal T_{x,J}$ that is normal to $\mathcal T_{x,I}$, i.e. normal to the ``boundary'' $M_I$, is 
\begin{equation}\label{un}
u_{n} = \frac{P_{x,J}\grad q(x)}{|P_{x,J}\grad q(x)|}, \qquad \text{where } \; P_{x,J} = T_{x,J}T_{x,J}^T
\end{equation}
 is the orthogonal projection matrix onto $\mathcal T_{x,J}$. Recall \eqref{vopt}, where a similar vector was defined to construct Lose neighbours.

\paragraph{Density of the proposal } 

Let $\mathcal A_{\rm gain}(x)\subset M_J$ be the set of values of $y$ that are accessible from $x$ via a Gain move using our NES; by definition they satisfy the inequalities in $\jin$. 
The proposal density is
 \begin{equation}\label{hgain}
h_{IJ}(y|x) = p_{IJ}\egain(v;x)|T_{x,J}^TT_{y,J}|1_{\mathcal A_{\rm gain}(x)}(y) \,.
\end{equation}

For a tangent step constructed as in \eqref{GainStep}, the tangent density factors as 
\begin{equation}\label{pgain}
p\egain_{IJ}(v;x) = p^\perp\gain(v_n)p^\parallel\gain(v_t|v_n)
\end{equation}
where the densities for the components of the step normal to and tangential to the boundary, are (in the one-sided case), respectively,
\begin{equation}\label{pgaincomp}
p^\perp\gain(v_n) = \frac{1}{\sigma\bdy}1_{[0,\sigma\bdy]}(v_n), \qquad
p^\parallel\gain(v_t|v_n) = \frac{1}{\sqrt{(2\pi\sigtan^2v_n^2)^{d_I}}}e^{-\frac{|v_t|^2}{2(\sigtan v_n)^2}}.
\end{equation}
For the two-sided case, the normal density is $p^\perp\gain(v_n) = \frac{1}{2\sigma\bdy}1_{[-\sigma\bdy,\sigma\bdy]}(v_n)$. 
We have omitted the densities' dependence on $I,J$ in the notation. 

The overall probability density (with respect to $\mu_J$) for a one-sided Gain move is 
\begin{equation}\label{GainDens}
\lambda_{IJ}(x)h_{IJ}(y|x) = \frac{\Lambda\gain(x)}{n\gain}\frac{1}{\sigma\bdy}1_{[0,\sigma\bdy]}(v_n)\frac{1}{\sqrt{(2\pi\sigtan^2v_n^2)^{d_I}}}e^{-\frac{|v_t|^2}{2(\sigtan v_n)^2}}|T_{x,J}^TT_{y,J}|1_{\mathcal A_{\rm gain}(x)}(y).
\end{equation}
For a two-sided move the modification is straightforward.

\paragraph{Evaluating the density for the reverse move}

Given $x\in M_I$, $y\in M_J$, the density \eqref{hgain} is evaluated in the same way as for a Same move: one first computes $v\in \mathcal T_{x,J}$, and then one runs the NES to check whether $y\in \mathcal A_{\rm gain}(x)$.

\subsubsection{Lose Move: $x\in M_I$, $y\in M_J$ with $J\in \nlosesig{I,x}$}

\paragraph{Generating the move} 

Suppose that $q$ is the constraint added to $\ieq$ to obtain $\jeq$. 
Recall from Section \ref{sec2:point} that we choose a tangent \emph{direction} $v\in \mathcal T_{x,I}$; the difference from Same and Gain moves is we now require $v$ to be a unit vector. Let the density for $v$ be $p_{IJ}\elose(v;x)$. 
We take a step in direction $v$ with unknown magnitude $\alpha$, and solve for the normal vector $w\in \mathcal N_{x,I}$ such that $y=x+\alpha v + w$ lies in $M_J$, as in \eqref{qsolveLose}. 
This is done by writing $w=Q_I(x)a$ for $a\in \R^{m_I}$ and solving the system of equations
\begin{equation}\label{loseeqns}
q_j(x+\alpha v + Q_I(x)a) = 0, \qquad j \in \jeq\,,
\end{equation}
for the unknown $(a,\alpha)\in \R^{m_J}$. 
There are $m_J$ equations and $m_J$ unknowns, so this system is well-posed in general. 
Call the method used to solve \eqref{loseeqns} NES-L. As before this method must be deterministic with a deterministic initial condition. 
Our specific method to solve using Newton's method is given in Appendix \ref{sec:summary}, Algorithm \ref{alg:takestep}.

If the solver NES-L finds a solution $(a,\alpha)$, \emph{and} if $\alpha > 0$, then the proposal is $y=x+\alpha v + Q_I(x)a$. If the reconstructed $y$ fails to satisfy the inequalities in $\jeq$, or if NES-L fails to find a solution, the proposal is rejected and set to $y=x,J=I$. 

We reject the proposal when $\alpha < 0$ to ensure the reverse move is reproducibile. If for some $v$ \eqref{loseeqns} has solution $(a,\alpha)$, then the solution for $-v$ is $(a,-\alpha)$. If we are given only $x,y$, then we don't know which of $v,-v$ was used in the proposal, so without a sign restriction on $\alpha$ we can't evaluate the density of that move, nor whether $y$ is accessible from $x$. 

This is also the reason why we choose a unit vector for $v$, not an arbitrary vector. If we chose an arbitrary vector, we could still solve \eqref{loseeqns}, but then given only $x,y$ we would not know which $v$ was actually used in the proposal.

\paragraph{Density of the proposal} 
Let $\mathcal A_{\rm lose}(x)\subset M_J$ be the set of values of $y$ that are accessible from $x$ via a Lose move using our NES-L; these necessarily satisfy the inequalities in $\jin$. 
The density of a successful proposal $y\in \mathcal A\lose(x)$ is $p_{IJ}\elose(v;x)\big|\pp{v}{y}\big|$, where as before $\big|\pp{v}{y}\big|$ is the determinant of the inverse of the Jacobian of the transformation from $v\to y$. The difference from the previous two move types is that now $v$ is constrained to be a unit vector. Therefore, to linear order $v$ can only vary in the subspace $\mathcal T_{x,I,v}=\mathcal T_{x,I}\cap (\mbox{span}\{v\})^\perp$, which is the part of the tangent space $\mathcal T_{x,I}$ that is orthogonal to vector $v$. Let $T_{x,I,v}\in \R^{n\times (d_I-1)}$ be a matrix whose columns form an orthonormal basis of $\mathcal T_{x,I,v}$.

The Jacobian factor can be calculated to be\footnote{
To calculate the Jacobian, suppose that $y,v$ vary as $y\to y+\Delta y$, $v\to v+\Delta v$, and simultaneously $\alpha\to\alpha + \Delta \alpha$. To linear order we must have $\Delta y = T_{y,J}b$ for some $b\in \R^{d_J}$, $\Delta v = T_{x,I,v}c$ for $c\in \R^{d_I-1}=\R^{d_J}$. The Jacobian $\pp{v}{y}$ will be the matrix $A$ such that to linear order $c = Ab$. 
Linearizing the equation $y = x+\alpha v + w$ with $w\perp \mathcal T_{x,I}$ gives
\[
\Delta y = \alpha \Delta v + v\Delta \alpha + \Delta w\,,
\]
where $\Delta w\perp \mathcal T_{x,I}$. Multiplying the above equation by $T_{x,I,v}^T$, using that $T_{x,I,v}^T w = 0$, $T_{x,I,v}^Tv = 0$, and substituting the expressions for $\Delta y,\Delta v$ gives
\[
T_{x,I,v}^TT_{y,J}b = \alpha c \qquad \Leftrightarrow \qquad c = \alpha^{-1}T_{x,I,v}^TT_{y,J}b\,.
\]
The Jacobian factor is therefore $|\alpha^{-1}T_{x,I,v}^TT_{y,J}|$. 
} 
\[
\Big|\pp{v}{y}\Big| = \alpha^{-d_J}|T_{x,I,v}^TT_{y,J}|,
\] 
where $\alpha = |y-x|$.
Therefore the proposal density is 
\begin{equation}\label{hlose}
h_{IJ}(y|x) = p_{IJ}\elose(v;x)|T^T_{x,I,v}T_{y,J}|y-x|^{-d_J}1_{A_{\rm lose}(x)}(y)\,.
\end{equation}

We now give the density for the tangent direction $v$ constructed from $v\opt$ and $T_{x,I,v\opt}$ as in \eqref{vlose}. 
First, note that given a step $v$, the tangential components $R$ are uniquely determined; call them $R=r(v)$. They can be found by scaling $v$ so it has unit length in the normal direction and then projecting onto the tangent space: 
\begin{equation}\label{rv}
r(v) = \frac{T^T_{x,I,v\opt}v}{v\cdot v\opt} \qquad \text{if } v\cdot v\opt > 0\,.
\end{equation} 
If $v\cdot v\opt < 0$, then $r(v)$ is undefined; such a step could never be proposed.

Because we can solve uniquely for $r(v)$, and because the density $p_R(r)$ for $R$ is known analytically, we obtain an analytic expression for the density in $v$, as $p\elose_{IJ}(v;x)=p_R(r)\big|\pp{r}{v}\big|1_{v\cdot v\opt > 0}$. Here $\big|\pp{r}{v}\big|$ is the determinant of the Jacobian of the transformation $v\to r$. Since there are no constraints on $r$, the Jacobian of the reverse transformation may be directly calculated from \eqref{vlose} to be
\begin{equation}\label{dvdr}
\pp{v}{r} = \frac{(I-v v^T)T_{x,I,v\opt}}{|v\opt + T_{x,I,v\opt}r|}.
\end{equation}
Therefore $|\partial v / \partial r| = |v\opt + T_{x,I,v\opt}r|^{-(d_I-1)}|(I-v v^T)T_{x,I,v\opt}|$, where $|\cdot|$ applied to a matrix is the pseudodeterminant. The power $(d_I-1)$ arises because this is the rank of the numerator. 

Substituting for $r$ using \eqref{rv} shows that $|v\opt + T_{x,I,v\opt}r| = |v / v\cdot v\opt| = |v\cdot v\opt|^{-1}$. For the numerator of \eqref{dvdr}, we claim its pseudodeterminant is 
\begin{equation}\label{losejac0}
|(I-v v^T)T_{x,I,v\opt}| = |v\cdot v\opt|. 
\end{equation}
To see this, note that the left-hand side has the form $|P_\mathcal A B|$, where  $P_{\mathcal A}$ is the orthogonal projector onto the subspace $\mathcal A = \mbox{span}\{v\}^\perp$, and $B = T_{x,I,v\opt}$ is a matrix with orthonormal columns spanning subspace $\mathcal B = \mathcal T_{x,I,v\opt}$. By Theorem 1  in \cite{Bjorck:1973jq} and the subsequent discussion, $|P_\mathcal A B| = \prod_{k}\cos\theta_k$, the product of the cosines of the principal angles $\theta_k$ between subspaces $\mathcal A, \mathcal B$. But the nonzero principal angles between $\mathcal A, \mathcal B$ equal the nonzero principal angles between their complements, $\mathcal A^\perp, \mathcal B^\perp$ \cite{Ipsen:1995ki,Knyazev:2007go}. Since $\mathcal A^\perp = \mbox{span}\{v\} \subset \mathcal B$, the largest principal angle between $\mathcal A^\perp, \mathcal B^\perp$ is the angle between $v,v\opt$, and the other principal angles are 0. This shows \eqref{losejac0}. 

Therefore $|\partial v / \partial r| = |v\cdot v\opt|^{d_I-1}|v\cdot v\opt|$, so 
\begin{equation}\label{rvjac}
\left|\pp{r}{v} \right| = |v\cdot v\opt|^{-d_I}, 
\end{equation}
so the tangent step density is 
\begin{equation}\label{plose}
p\elose_{IJ}(v;x)=\frac{1}{\sqrt{(2\pi\sigtan^2)^{d_J}}}e^{-\frac{|r(v)|^2}{2\sigtan^2}}|v\cdot v\opt|^{-d_I}1_{v\cdot v\opt > 0}(v)\,. 
\end{equation}

The overall probability density (with respect to $\mu_J(dx)$) of proposing a Lose move $(x,I)\to (y,J)$ is
\begin{multline}\label{LoseDens}
\lambda_{IJ}(x)h_{IJ}(y|x) = \\
\frac{\Lambda\lose(x)}{n\lose} \frac{1}{\sqrt{(2\pi\sigtan^2)^{d_J}}}e^{-\frac{|r(v)|^2}{2\sigtan^2}}|v\cdot v\opt|^{-d_I}1_{v\cdot v\opt > 0}(v)|T^T_{x,I,v}T_{y,J}|y-x|^{-d_J}1_{A_{\rm lose}(x)}(y).
\end{multline}

\paragraph{Evaluating the density for the reverse move}

Given $x\in M_I$, $y\in M_J$ we find $v$ by first calculating $\tilde v = T_{x,I}T^T_{x,I}(y-x)$, and then setting $v = \tilde v / |\tilde v|$, $\alpha = |\tilde v|$. We must also check that $y\in \mathcal A_{\rm lose}(x)$, by running the NES-L to solve \eqref{loseeqns} with the calculated $v$. In addition to checking that the solver converges, and gives $y$ as a solution, we must also check that the solution $(a',\alpha')$ has $\alpha'>0$.


\section{Acceptance probabilities for flat manifolds}\label{sec2:flat}

This section shows that our choice of Gain and Lose proposal moves give 100\% acceptance probabilities for a stratification consisting of two flat manifolds defined by affine constraints with no additional inequalities. We first motivate these proposals by considering in detail the simplest possible cases, a 0-dimensional point at the boundary of a 1-dimensional line, and then a 1-dimensional line at the boundary of a 2-dimensional plane. We explain how to construct proposals, including proposals different from the ones we introduced, to obtain a 100\% acceptance probability. Finally we show that the acceptance probability is 100\% for the general linear case, a $d$-dimensional hyperplane forming the boundary for a $(d+1)$-dimensional hyperplane.

\subsection{0-dimensional manifold $\leftrightarrow$ 1-dimensional manifold}\label{sec:01d}

  \definecolor{mypurple}{rgb}{0.3, 0.1, 0.5} 
  \definecolor{myblue}{rgb}{0.1, 0.4, 0.8}  
  \definecolor{mygreen}{rgb}{0.2, 0.7, 0.2} 
  \definecolor{myred}{rgb}{0.66, 0.06, 0.06} 
  
    \definecolor{myblue1}{rgb}{0.4, 0.6, 0.8}  
    \definecolor{myblue2}{rgb}{0.36, 0.54, 0.66} 
    \definecolor{myblue3}{rgb}{0.19, 0.55, 0.91}  
    \definecolor{myblue4}{rgb}{0.0, 0.53, 0.74}  
    \definecolor{myblue5}{rgb}{0.29, 0.59, 0.82} 
    \definecolor{myblue6}{rgb}{0.0, 0.48, 0.65} 
     
\begin{figure}
\centering
\begin{tikzpicture}[scale=2.5]
\draw[ultra thick,myblue] (0,0) --node[midway,label={[align=left,font=\small,black]${\color{myblue} \lambda\lose}=cst$
\\${\color{black}\lambda\same}=1-{\color{myblue}\lambda\lose}$
}] {} (1.5,0);   
\draw[thick] (1.5,0) --node[midway,label={[align=left,font=\small]${\color{black}\lambda\same}=1$}] {} (3.5,0); 
\draw[ultra thick,mypurple] (1.5,0.1) -- (1.5,-0.1) node[below] {{\color{mypurple}$\sigma\bdy$}};
\node[circle, draw, fill=myred, inner sep=0pt, minimum width=6pt,label={[align=left,font=\small]left:${\color{myred} \lambda\gain}={\color{mypurple}\sigma\bdy}{\color{myblue}\lambda\lose}$
\\${\color{black}\lambda\same}=1-{\color{myred} \lambda\gain}$
}] at (0,0) {};
\end{tikzpicture}

\caption{Schematic showing optimal choice of move types for a stratification with a 0-dimensional point and a 1-dimensional line, as in Section \ref{sec:01d}. 
}\label{fig:01}
\end{figure}
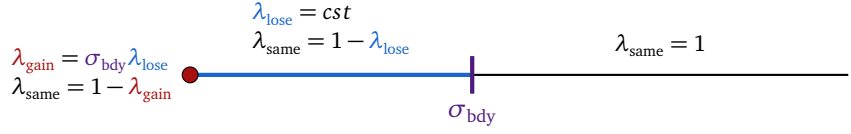

Let $x\in \R$, and define two manifolds $M_0$, $M_1$ by 
\begin{equation}
M_0 = \{0\}, \qquad M_1 = \{x\in \R: x>0\}.
\end{equation}
Let $\mathcal S = M_0\cup M_1$. One way to construct such a stratification via level sets is to define  a function $q(x) = x$. Let $f_1(x),f_0$ be functions defining the probability measure on each manifold.

Let $\lambda_{\rm gain},\lambda_{\rm lose}(x)$ be the probabilities of choosing a Gain and Lose move respectively from a given point (there is only one possible starting point for a Gain move, so we omit it in the notation.) 
Let the $v$-density for a Gain move be $p_{\rm gain}(v)$, and let the $v$-density for a Lose move be $p_{\rm lose}(v)$. Since there are only two unit vectors in the tangent plane to $M_1$, $\pm 1$, we make the decision that $p_{\rm lose}(v) = \delta_{-1}(v)$, a unit mass at $v=-1$. This ensures we always propose to move in the direction of the boundary, not away from it. For a general one-dimensional affine constraint $q(x)$, we can determine which direction points toward the boundary by evaluating $\dd{q}{x}$. 

Consider a Gain move from $x=0$ to $y\in M_1$. The tangent step is $v = y$ and the reverse step is $v_{\rm rev} = -y/|y| = -1$. The Jacobian factors are both 1 in this simple example, so the Metropolis ratio in the acceptance probability \eqref{acc} is
\begin{equation}
a(y,1|x,0) = \frac{f_1(y)\lambda_{\rm lose}(y)h_{01}(x|y)}{f_0\lambda_{\rm gain}h_{10}(y|x)} = \frac{f_1(y)}{f_0}\frac{\lambda_{\rm lose}(y)}{\lambda_{\rm gain}p_{\rm gain}(y)} 
\end{equation}
We wish to make the acceptance ratio as high as possible, to avoid wasting computation in proposing moves. 
For simplicity suppose that $f_1(y)=f_1$, a constant independent of $y$. For the acceptance ratio to be 1, we need
\begin{equation}\label{lam01}
\lambda\lose(y) = \frac{f_0}{f_1}\lambda\gain p\gain(y)\,.
\end{equation}
Since $p\gain(y)$ is a probability density, whose integral is 1, we also need that 
\begin{equation}\label{lamgain01}
\lambda\gain = \frac{f_1}{f_0}\int_0^\infty \lambda\lose(y)dy\,.
\end{equation}
Therefore once we choose $p\gain(y)$, this sets the ratio $\lambda\lose(y)/\lambda\gain$ via \eqref{lam01}, and conversely if we choose $\lambda\lose(y)$, this sets both $\lambda\gain$, via \eqref{lamgain01}, and then $p\gain(y)$ via \eqref{lam01}. 

A proposal that can be implemented in practice is to choose a Lose move with constant probability within a certain cutoff distance $\sigma\bdy$ from the boundary. This implies that $p\gain(y)$ must be uniform on $[0,\sigma\bdy]$, so
\begin{equation}\label{params01}
\lambda\lose(y) = \lambda\lose1_{[0,\sigma\bdy]}(y), \quad\;\;
p\gain(y) = \frac{1}{\sigma\bdy}1_{[0,\sigma\bdy]}(y), \quad\;\;
\lambda\gain = \frac{f_1}{f_0}\sigma\bdy\lambda\lose.
\end{equation}
There are two parameters defining these proposals, the cutoff $\sigma\bdy$, and the lose parameter $\lambda\lose$, which is subject to the constraint $0 < \lambda\lose <  \min(1, \sigma\bdy^{-1})$. The probability of proposing a Same move along $M_1$ is then $\lambda\same(y) = 1$ if $y > \sigma\bdy$, and $\lambda\same(y) =1-\lambda\lose$ if  $y < \sigma\bdy$. The probability of proposing a Same move on $M_0$ (repeating the point) is $\lambda\same(0) =1-\lambda\gain$.

\subsection{1-dimensional manifold $\leftrightarrow$ 2-dimensional manifold}\label{sec:12d}

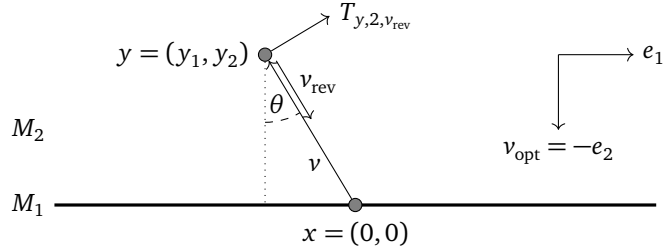
\begin{figure}
\centering
\begin{tikzpicture}
\def\yx{-1.2};
\def\yy{2};
\draw[very thick] (-4,0) -- (4,0);
\node[left,align=right] at (-4,0) {$M_1$};
\node[left,align=right] at (-4,1) {$M_2$};
\node[vertex,label=below:{$x=(0,0)$}] (x) at (0,0) {};  
\node[vertex,label=left:{$y=(y_1,y_2)$}] (y) at (\yx,\yy) {};  
\draw[dotted] (y) -- (\yx,0);  
\draw[->] (x) -- node[pos=0.25,left]{$v$} (y);
\draw[->,transform canvas={xshift=1mm}] (y) --node[pos=0.4,right]{$v\rev$} ($(y)+({cos(atan(-\yy/\yx))},{-sin(atan(-\yy/\yx))})$);
\centerarc[dashed](y)(-90:{-atan(-\yy/\yx)}:0.9);
\node at ({\yx+0.16},{\yy-0.65}) {$\theta$};
\draw[->] (y) -- ($(y) + ({cos(atan(-\yx/\yy))},{sin(atan(-\yx/\yy))})$) node[right] {$T_{y,2,v\rev}$};  
\begin{scope}[xshift=2.7cm,yshift = 1cm]
\draw[->] (0,1) -- (0,0) node[below]{$v\opt = -e_2$};
\draw[->] (0,1) -- (1,1) node[right] {$e_1$};
\end{scope}
\end{tikzpicture}
\caption{Setup for the example in Section \ref{sec:12d}, which constructs optimal proposal moves between the $x$-axis $M_1$ and the upper half plane $M_2$. }\label{fig:12}
\end{figure}

Next let $x=(x_1,x_2)\in\R^2$, and define two manifolds $M_1,M_2$ by 
\begin{equation}
M_1 = \{(x_1,x_2)\in \R^2: x_2=0\}, \qquad M_2 = \{ (x_1,x_2) \in \R^2 : x_2 > 0\}\,.
\end{equation}
See Figure \ref{fig:12}. For simplicity we assume that $f_1(x_1) = f_2(x_1,x_2) = cst$ everywhere. 
Let $e_1= (1,0)$, $e_2 = (0,1)$. 
We may assume that proposal moves are homogeneous in the $x_1$ direction, so we must choose functions $\lambda\gain$, $\lambda\lose(x_2)$, $p\gain(v)$, $p\lose(\theta;x_2)$. Here $\theta$ is the direction of a proposed Lose move $v$, measured counterclockwise from $v\opt = -e_2$; the corresponding unit vector is $v = (\sin\theta,-\cos\theta)$.

Consider a Gain move from $x=(x_1,0)\in M_1$ to $y=(y_1,y_2)\in M_2$. 
Since we assume homogeneity in $x_1$, we may assume without loss of generality that $x_1 = 0$. 
The forward tangent step is $v = y-x = (v_1,v_2) = (y_1,y_2)$. 
The Jacobian factor $\big|\pp{v}{y}\big|$ in \eqref{hgain} is $|T^T_{x,2}T^T_{y,2}| = 1$, since the step is along a plane. 

The reverse tangent step is $v_{\rm rev} = (v_{{\rm rev}, 1},v_{{\rm rev},2})=(y-x)/|y-x| = (\sin\theta,-\cos\theta)$ for some angle $\theta$. 
To calculate the jacobian factor $\big|\pp{v}{y}\big|$ in \eqref{hlose} for the Lose move from $y\to x$, we calculate $T_{y,2,v_{\rm rev}} = (\cos\theta,\sin\theta)$, the unit vector perpendicular to $v_{\rm rev}$. See Figure \ref{fig:12}. 
We have that $T_{x,1} = e_1$, the tangent space to the horizontal axis. Therefore $|T_{y,2,v_{\rm rev}}^TT_{x,1}| = |\cos \theta|$. By geometry, $\alpha = y_2/|\cos\theta|$, and it is raised to the power of $d_J = 1$ (where $J$ represents $M_1$ for the reverse move.) 
Therefore the Jacobian factor for the Lose move from $y\to x$ is 
\[
\frac{|T_{y,2,v}^TT_{x,1}|}{\alpha} = \frac{\cos^2 \theta}{y_2}\,.
\]
Plugging this into  \eqref{acc} gives a Metropolis ratio in the acceptance probability of
\begin{equation}\label{a21}
a(y,2|x,1) = \frac{p\lose(\theta;y_2)\lambda\lose(y_2)\cos^2\theta}{|y_2|p\gain(v)\lambda\gain}\,.
\end{equation}

We next explain how to choose the other proposal functions to make the acceptance probability equal 1. 
First, we will split $p\gain(v)$ into two parts, one which steps in the normal direction, and one which steps in the tangential direction. The aim is to use the normal component $v_2=y_2$ to balance the factor $\lambda\lose(y_2)/\lambda\gain$ as for the 1-dimensional problem in Section \ref{sec:01d}. These factors only depend on $y_2$, so $v_2$ can be chosen independently. The tangential component $v_1$ must also depend on $v_2$, since there is still a factor of $|y_2|$ in the denominator of \eqref{a21} that must be balanced. Therefore we choose the density to have the form
\begin{equation}
p\gain(v) = p\gain^\parallel(v_1|v_2)p\gain^\perp(v_2),
\end{equation}
where $p\gain^\perp(v_2)$ is the density for the normal component of the step, and $p\gain^\parallel(v_1|v_2)$ is the density for the tangential component, given the value of the normal component. 

Next, we use the choice \eqref{vlose} for $p\lose(\theta;y_2)$ and show this choice cancels the factor of $\cos^2\theta$ in the numerator. For this example, \eqref{vlose} constructs the Lose tangent step $V=(V_1,V_2)$ as
\begin{equation}
V = \frac{-e_2 + Re_1}{|-e_2+Re_1|}, \qquad R\sim p_R(r)\,,
\end{equation}
where $p_R(r)$ is the probability density of the random variable $R\in \R$, which we leave arbitrary for the moment. 
The inverse mapping from $V$ to $R$ takes the form 
\begin{equation}
r(v) = \frac{v\cdot e_1}{|-v\cdot e_2|} = \frac{\sin\theta}{|\cos\theta|} = \frac{v_1}{v_2}.
\end{equation}
The determinant of the jacobian of this transformation is given in \eqref{rvjac} to be $\big|\pp{r}{v}\big| = |v\cdot -e_2|^{-d_1} = |\cos\theta|^{-2}$, so the overall density for proposing the reverse step $v_{\rm rev}$ is 
\[
p\lose(\theta;y_2) = p_R\Big(\frac{v_{\rm rev,1}}{v_{\rm rev,2}}\Big)(\cos\theta)^{-2}\,.
\]
Substituting these specific Gain and Lose densities into \eqref{a21} gives an acceptance probability 
\begin{equation}\label{a21b}
a(y,2|x,1) = \frac{p_R\big(\frac{y_1}{y_2}\big)\lambda\lose(y_2)}{|y_2|p\gain^\parallel(y_1|y_2)p\gain^\perp(y_2)\lambda\gain}\,.
\end{equation}
We have written $v$ in terms of $y_1,y_2$, and used that $v_{\rm rev,1}/v_{\rm rev,2} = y_1/y_2$. 

Now, we show how to choose specific forms for these densities so the acceptance probability is 1. 
To balance the normal components, we can choose the densities so that 
\begin{equation}
\lambda\lose(y_2) = p^\perp\gain(y_2)\lambda\gain\,.
\end{equation}
This has the same form as \eqref{lam01} for the 1-dimensional problem, so the same considerations apply. In particular, our choice in \eqref{params01}, where $\lambda\lose(y_2)$ is constant within an interval $y_2 \in [0,\sigma\bdy]$,  $p^\perp\gain(y_2)$ is uniform on that same interval, and $\lambda\gain=\sigma\bdy\lambda\lose$, will work to satisfy the above equation. 

The remaining densities must be chosen so that 
\begin{equation}\label{pperp12}
p\gain^\parallel(y_1|y_2) = \frac{1}{y_2}p_R\Big(\frac{y_1}{y_2}\Big).
\end{equation}
A quick calculation shows that this is possible, since $\int_{-\infty}^\infty p\gain^\parallel(y_1|y_2)dy_1 = \int_{-\infty}^\infty p_R(r)dr = 1$, by changing variables to $r = y_1/y_2$. 
Therefore, once we choose a density for the tangential step of a Lose move, this sets the density for a tangential step for a Gain move via \eqref{pperp12}. For example, our choices from Section \ref{sec2:overview} would take the form 
\[
p_R(r) = \frac{1}{\sqrt{2\pi\sigtan^2}}e^{-\frac{r^2}{2\sigtan^2}}, \qquad 
p\gain^\parallel(y_1|y_2) = \frac{1}{\sqrt{2\pi(\sigtan y_2)^2}}e^{-\frac{y_1^2}{2(\sigtan y_2)^2}}.
\]

\medskip

We make a side remark that an alternative method could choose the tangential step for a Lose move in a way that depends on $y_2$, and let the tangential step for a Gain move be independent of distance. In this case the densities would have to solve 
\[
p_R\Big(\frac{y_1}{y_2}\Big|y_2\Big) = y_2 p\gain^\parallel(y_1).
\]
While we haven't tried it, we expect this  to have poorer performance compared to the choice \eqref{pperp12}, since it implies the tangential variance for a Lose move should \emph{increase} as the distance to the boundary increases, so the Lose move attempts to jump to a greater area of the boundary. For curved manifolds, this should lead to increased rejections than the method presented here, because jumps to further-away pieces of the boundary will cause the NES-L to fail more often, and will lower the acceptance probability since it will deviate even more from the flat case.

\subsection{General case: $d$-dimensional manifold $\leftrightarrow$ $d+1$-dimensional manifold}\label{sec:dgeneral}

Now we show that the choice of proposal densities in Section \ref{sec2:overview} leads to 100\% acceptance probability for two flat manifolds defined by affine constraints and no additional inequalities. 
Consider a stratification consisting of two manifolds: 
\begin{equation}
M_d = \{x\in \R^{d+1}: q(x) = 0\}, \quad M_{d+1} = \{x\in \R^{d+1}: q(x) > 0\}, \qquad
q(x) = a\cdot x + b
\end{equation}
where $a,b\in\R^{d+1}$ are constant vectors; we may assume by changing variables that $|a|=1$. Let $I,J$ represent the labels for manifolds $M_d, M_{d+1}$ respectively.  The dimensions of the manifolds are $d_I = d$, $d_{J} = d+1$. We could incorporate other affine constraints that are the same for both manifolds, but by reparameterizing would obtain a case equivalent to the above.  Assume the function we wish to sample is constant everywhere, $f(x) = cst$. 

Consider a Gain move from $x\in M_d$ to $y\in M_{d+1}$. This is a one-sided move. 
The tangent step for the Gain move is $v = y-x$ and the tangent step for the reversed Lose move is $v\rev = (x-y)/|x-y|$. The Metropolis ratio in the acceptance probability \eqref{acc} is calculated by substituting the densities \eqref{GainDens}, \eqref{LoseDens}:
\[
a(y,J|x,I) = \frac{\lambda_{JI}(y)h_{JI}(x|y)}{\lambda_{IJ}(x)h_{IJ}(y|x)}
= \frac{\lambda\lose }{\lambda\gain}\frac{p\elose_{JI}(v\rev)|T^T_{y,J,v\rev}T_{x,I}||x-y|^{-d_I}}{p_{IJ}\egain(v)|T_{x,J}^TT_{y,J}|}.
\]
We assume that NES, NES-L converge for any proposal when the constraints are affine.

Now we calculate the geometrical factors and substitute for the tangent step densities $p_{IJ}\egain(v)$, $p_{JI}\elose(v\rev)$. 
We have $|T^T_{y,I,v\rev}T_{x,J}| = 1$ since the Gain move is on a flat manifold. 
We have $|T^T_{y,J,v\rev}T_{x,I}| = |v\rev\cdot a|$, the cosine of the angle between the normal vectors of subspaces $\mathcal T_{y,J,v\rev}, \mathcal T_{x,I}$.
For calculating the densities, we have that $u_n = -v\opt = a$ since $q$ is affine (see \eqref{vopt},\eqref{un}.) 
Let $v_n = v\cdot u_n$, $v_t = T_{x,I,u_n}v$ be the normal and tangential components of the Gain tangent step. The proposal is such that $v_n < \sigma\bdy$, so  $\lambda\lose/\lambda\gain = 1/\sigma\bdy$. 
Substituting directly for all the factors above, 
\begin{align*}
a(y,J|x,I) &= \frac{1}{\sigma\bdy}\frac{ \frac{1}{\sqrt{(2\pi\sigtan^2)^{d}}}e^{-\frac{|r(v\rev)|^2}{2\sigtan^2}}|v\rev\cdot a|^{-(d+1)}|v\rev\cdot a||v|^{-d}}{ \frac{1}{\sigma\bdy} \frac{1}{\sqrt{(2\pi\sigtan^2v_n^2)^{d}}}e^{-\frac{|v_t|^2}{2(\sigtan v_n)^2}}}\\
&= \frac{ e^{-\frac{|r(v\rev)|^2}{2\sigtan^2}}|v\rev\cdot a|^{-d}|v|^{-d}\:v_n^d}{e^{-\frac{|v_t|^2}{2(\sigtan v_n)^2}}}
\end{align*}
Now we use that, by geometry, $|v_n|/|v| = |v\rev\cdot a|$, and $|v_t|/|v_n| = |r(v\rev)|$, to obtain that $a(y,J|x,I)=1$.


\section{Discussion and Conclusion}\label{sec:conclusion}

We introduced a Monte Carlo method to sample systems with constraints that can break and form, such as particles with variable bond distance constraints. Mathematically, the method generates samples from a probability distribution defined on a stratification, a union of manifolds of different dimensions. 
Our method  is a natural extension of a random walk Metropolis scheme on a manifold \cite{Zappa:2018jy}, and is perhaps the simplest possible sampler that works for a general stratification defined by level sets of functions. 
We illustrated the method with several pedagogical examples, ranging from sticky spheres, which are a model for DNA-coated colloids, to polymers adsorbing to a surface, to calculating volumes of high-dimensional manifolds. We showed the method gives an accurate description of a system with a stiff, strong, but reversible bonds, as is the case for DNA-coated colloids.

We have implemented the algorithm numerically using dense arithmetic and noniterative matrix factorizations, which makes it fast for small systems, of around a dozen particles, but slow for large ones, since the computational complexity of the factorizations grows as $O(n^3)$. To use this algorithm for larger systems will require addressing several issues in numerical linear algebra and optimization: (i) proposing a random step in a tangent space without calculating a basis for that tangent space; (ii) projecting a step back to the manifold of interest; (iii) calculating Jacobian factors required for the Metropolis ratios for Gain and Lose moves (the terms $|T_{x,J}^TT_{y,J}|$, $|T_{x,I,v}^TT_{y,J}|$, in \eqref{hgain}, \eqref{hlose}), and (iv) calculating the reverse move. Item (ii) is standard in constraint-based MC or MD simulations, and many papers have proposed iterative techniques to solve the necessary nonlinear equations, which have an $O(n)$ computational complexity if the constraints are sufficiently sparse \cite{Ryckaert:1977wr,Barth:1995vi,Weinbach:2005jt,Eastman:2010ts}. 
The other items are more specific to the stratification sampler, but we are confident that efficient, iterative solutions can be found. 


Beyond using more sophisticated numerical techniques to technically implement our algorithm, there are other ways to build upon the ideas presented here to make a method that samples a stratification more efficiently. 
For example, 
one could adapt other sampling schemes for manifolds, such as Hamiltonian \cite{Lelievre:2019be} or geodesic methods \cite{Leimkuhler:2016ia}, to work on stratifications; both of these methods have velocity-like terms that can help make successive samples less correlated. 
In addition, one can more carefully consider how to add and subtract constraints when there are several manifolds in the stratification. Our proposal moves were carefully constructed to achieve 100\% acceptance probability for two flat manifolds defined by affine constraints, but the strategy with  highest acceptance probability in more general stratifications might depend on the connectivity structure of the manifolds. 

An interesting extension of our method would be to generate a trajectory $X_1,X_2,\ldots$ that approximates the \emph{dynamics} of a particular system. 
Our method is purely a sampling method, which can generate approximate samples from a probability distribution containing singular measures, but is not designed to solve a particular set of dynamical equations. However, dynamics are often important, especially for systems such as DNA-coated colloids, which diffuse slowly hence don't always reach their stationary distribution on the timescale of observation. 
The challenge is in describing and simulating the dynamics of a diffusion process that change its intrinsic dimension, since it is not obvious how the dynamics should behave near changes in dimension.  
A step in this direction was taken in \cite{BouRabee:2020gg}, which introduced a numerical method to simulate a ``sticky Brownian motion'', a one-dimensional Brownian motion that can stick to the origin and spend finite time there. We are optimistic that the techniques from \cite{BouRabee:2020gg} can be built upon to handle higher-dimensional diffusion processes on stratifications, perhaps even by combining with tools from this paper: because our method achieves high acceptance probabilities, it already approximates some kind of dynamics. 

One limitation of our method is the regularity assumption that the constraint gradients are always linearly independent. For sticky spheres this assumption does not always hold. For example, a system of $N$ sticky unit spheres can assemble into a fragment of an fcc lattice, a packing which contains nearly 6 contacts per sphere, instead of the usual 3, creating redundant constraints. While redundant constraints can in principle be removed, a more subtle issue arises when constraints become linearly dependent without introducing redundant ones; for examples where this occurs with unit sphere packings see  \cite{HolmesCerfon:2016wa}. In these cases, our algorithm will fail because we will not be able to calculate the tangent space using only gradients of the constraints, and, even if we do take a step along the tangent space, projecting back to the manifold is not always well-defined because there are more variables than equations or vice versa. Sometimes it is possible to deal with such degenerate configurations by formulating additional optimization problems that can be solved by semidefinite programming \cite{HolmesCerfon:2019wm}, however it is not yet clear how to build these optimization methods into our sampler.

We expect the Stratification sampler to be useful in a variety of other applications. 
Constraints are frequently used in biology,  to model stiff bonds as well as flexible structures, such as networks of actin filaments or ligand-receptor systems, but the rate-limiting step in simulations is often set by the stiffest bond that can break \cite{Nedelec:2007bn}. 
Sampling methods are being explored in machine learning, where constraints on the weights of a neural network may offer efficient ways to explore parameter spaces of interest \cite{Leimkuhler:2019wr}. 
Our method may also generate interesting data for problems in discrete geometry; one example is to ask which kinds of graphs may be formed from packing spheres together \cite{Connelly:2019uq}. 
We hope the method can be adapted to these kinds of systems, and others yet to be imagined.

\subsection*{Acknowledgements}
Thank you to Jonathan Goodman, Jeff Cheeger, and Mark Goresky for helpful discussions regarding sampling and stratifications, to Nawaf Bou-Rabee and Glen Hocky for carefully reading early drafts of the paper, and to Anthony Trubiano for pointing out bugs in an earlier version of the code. 
I would also like to acknowledge support from the US Department of Energy DE-SC0012296, and the Alfred P. Sloan Foundation. 

\subsection*{Data availability statement} 
Codes to run each example, the data produced by each example, and Matlab scripts to generate the figures in this paper are available at \cite{gitfigures}. 


\appendix 

\section{Appendix: Verifying the invariant measure}\label{sec:DB}

In this section we show that the measure $\rho$ defined in \eqref{rho} is an invariant measure for the 
Markov chain $X_1,X_2,\ldots$ constructed in Section \ref{sec2:overview}: if $X_{k}\sim \rho$, then $X_{k+1}\sim \rho$. This argument is very similar to that in Green \cite{Green:1995dq} Section 3.2, which shows that RJMC satisfies detailed balance, the difference being that that Green \cite{Green:1995dq} considers a stratification formed from a union of Euclidean spaces, whereas we consider a  union of embedded manifolds. We include the argument nonetheless for completeness. 

To start, define the product space $\mathcal C =\mathcal S\times \mathcal I$.  
It is convenient to change notation slightly and let $\bar \rho$ be the measure in \eqref{rho}, and define 
 $\rho(dx,I) = Z^{-1}f_I(x)\mu_I(dx)$ to be the desired joint probability distribution to find the system at configuration $x\in \mathcal S$ and with labels $I\in  \mathcal I$. The measure in \eqref{rho} is obtained by summing over labels as $\bar \rho(dx) = \sum_{I\in\mathcal I}\rho(dx,I)$. 
We will show that (using our new notation)
\begin{equation}
(X_k,I_k) \sim \rho \qquad \Rightarrow \qquad (X_{k+1},I_{k+1}) \sim \rho \,,
\end{equation}
which implies the parallel result for $\bar\rho$. 

Let $x=X_k, I=I_k$. The probability to successfully propose $(y,J)$ starting at $(x,I)$ has distribution 
\begin{equation}
P(dy,J|x,I) = \lambda_{IJ}(x)h_{IJ}(y|x)a(y,J|x,I)\mu_J(dy)
\;+\; \xi_I(x)\delta_I(y-x)\mu_J(dy). 
\end{equation}
The first term in the sum is the probability distribution for proposing $(y,J)$ and then accepting it; the second term is the probability distribution for remaining at $(x,I)$. The latter is a product of $\xi_I(x)$, the probability of remaining at the starting point $(x,I)$ (either because of failing to produce a proposal point, or because the proposal point was rejected), and a measure $\delta_I(y-x)\mu_J(dy)$. 
We write $\delta_I(y-x)$ to mean the delta-function with respect to the measure $\mu_I$: it is defined so that, for any integrable function $g:M_J\to \R$, we have $\int_{y\in zM_J}g(y)\delta_I(y-x)\mu_J(dy) = g(x)$ if $I=J$, and $\int_{y\in M_J}g(y)\delta_I(y-x)\mu_J(dy) = 0 $ otherwise. 

Now suppose that $(X_k,I_k) \sim \rho$. Let $\rho'$ be the distribution of $(X_{k+1},I_{k+1})$. It is computed as 
\begin{align}
\rho'(dy,J) &= Z^{-1}\sum_{I\in\mathcal I}\int_{x\in\mathcal S}P(dy,J|x,I)\rho_I(dx)  \nonumber\\
&= Z^{-1}\sum_{I\in\mathcal I}\int_{x\in M_I} \Big( \lambda_{IJ}(x)h_{IJ}(y|x)a(y,J|x,I)\mu_J(dy)
\;+\; \xi_I(x)\delta_I(y-x)\mu_J(dy)\Big) f_I(x)\mu_I(dx)  \nonumber\\
&= Z^{-1}\mu_J(dy)\left( \sum_{I\in\mathcal I}\int_{x\in M_I} \!\! f_I(x)\lambda_{IJ}(x)h_{IJ}(y|x)a(y,J|x,I)\mu_I(dx)
\quad + \;\;  \xi_J(y)f_J(y) \right).\label{rhoprime}
\end{align}
In the last step we used $\int_{x\in M_I}f_I(x)\xi_I(x)\delta_I(y-x)\mu_I(dx) = f_I(y)\xi_I(y)$, and then we used that $\sum_{I\in \mathcal I}f_I(y)\xi_I(y)\mu_J(dy) = f_J(y)\xi_J(y)\mu_J(dy)$, since measure $\mu_J(y)$ only gives nonzero weight to $y\in M_J$, and if $y\in M_J$ and $I\neq J$ then $\xi_I(y) = 0$, by definition.

Now, suppose that the acceptance probability $a$ satisfies the following relation: 
\begin{equation}\label{accrel}
f(x)\lambda_{IJ}(x)h_{IJ}(y|x)a(y,J|x,I) = f(y)\lambda_{JI}(y)h_{JI}(x|y)a(x,I|y,J)\,.
\end{equation}
It is straightforward to verify that the particular choice of $a$ in \eqref{acc} satisfies this. 
Then, substituting for $a(y,J|x,I)$ in \eqref{rhoprime} gives
\begin{equation}\label{rhoprime2}
\rho'(dy,J) = Z^{-1}f_J(y)\mu_J(dy) \left( \sum_{I\in\mathcal I}\int_{x\in M_I} \!\! \lambda_{JI}(y)h_{JI}(x|y)a(x,I|y,J)\mu_I(dx)   \quad + \;\;  \xi_J(y)\right).
\end{equation}

As a final step we rewrite $\xi_J(y)$ in terms of the other ingredients of the proposal distribution. This function, which is the total probability of remaining at the starting point $(y,J)$, equals the probability of failing to produce a proposal, plus the probability of producing a proposal but rejecting it. Since the only other possibility is to to produce a proposal and accept it, we must have that $\xi_J(y)$ plus the probability of producing a successful proposal equals 1:
\begin{equation}\label{xiJ}
\xi_J(y) = 1 - \sum_{I\in\mathcal I}\int_{x\in M_I} \lambda_{JI}(y)h_{JI}(x|y)a(x,I|y,J)\mu_I(dx). 
\end{equation}
Substituting \eqref{xiJ} into \eqref{rhoprime2} gives 
\[
\rho'(dy,J) =  Z^{-1}f_J(y)\mu_J(dy)\,,
\]
which is the desired result.

\section{Appendix: Summary of the algorithm}\label{sec:summary}

This section provides pseudocode, written at a high level to show the main steps of the algorithm. 
The pseudocode depends on having first constructed a stratification, a topic we comment on in Section \ref{sec:stratconstruct}, before describing the main algorithms provided in the pseudocode in Section \ref{sec:pseudocode}. 
The code we used to run our examples is available at \cite{gitfigures}. 

\medskip

The algorithm depends on 5 parameters: $\sigma$, $\lambda\lose$, $\lambda\gain$ $\sigma\bdy$, $\sigtan$. These parameters could depend on the manifold, but in our implementation they are independent of manifold. 
In general $\sigma$ may be different for each manifold, and $\lambda\lose,\lambda\gain,\sigtan,\sigma\bdy$ may be different for each pair of manifolds.

\subsection{Constructing a Stratification}\label{sec:stratconstruct}

To form a stratification we need a set of functions $\mathcal Q$ and a set of labels $\mathcal I$. Here is how we represent these on a computer. 

We start from a list of functions $\mathcal Q = \{q_1,\ldots,q_{\nfcns}\}$ and their gradients $\{\grad q_1,\ldots, \grad q_{\nfcns}\}$, where $q_i:\R^n\to \R$, and $\nfcns = |\mathcal Q|$. These are implemented by providing functions $\texttt{eq}(i,x)$ and $\texttt{jac}(i,x)$ defined so that $q_i(x) = \texttt{eq}(i,x)$, $\grad q_i(x) = \texttt{jac}(i,x)$. Both $\texttt{eq},\texttt{jac}$ are application-dependent; to change the application one simply needs to rewrite these functions, and change the number of variables $n$ and number of functions $\nfcns$. 

The labels for a given manifold in the stratification are represented with vector 
\[L\in \{C_{\rm eq}, C_{\rm in}, C_{\rm none}\}^{\nfcns}\]
where $C_{\rm eq} = 1$, $C_{\rm in} =2$, $C_{\rm none}=0$ are constants labelling whether each function is a constraint, inequality, or neither, respectively. For example, the labels $L = (C_{\rm eq}, C_{\rm eq}, C_{\rm none}, C_{\rm in})$ represent the manifold $M_{L} = \{x\in \R^n: q_1(x) = 0, q_2(x) = 0, q_4(x) > 0\}$. We change notation slightly in this section, and use the vector $L$ instead of the corresponding subset $I$ to represent labels for manifolds and tangent spaces. 

Recall that for a given manifold $M_L$, the number of equations that define it is written as $m_L$ and its intrinsic dimension is $d_L$. The number of equations is computed as $m_L = \# i: L_i = C_{\rm eq}$, and $d_L = n-m_L$ by our assumption that the gradients of the constraints are linearly independent everywhere on $M_L$. 

To finish constructing a stratification we need a set of labels $\mathcal I$ indicating which manifolds are in the stratification. We do this in one of two ways: either by providing a list of labels in advance, or by forming $\mathcal I$ from all subsets of labels obtained by letting certain functions vary between constraint/inequality. We distinguish these cases in the code by the parameter $\nmanif$. When $\nmanif > 0$ we use the first method, and $\nmanif=|\mathcal I|$ represents the total number of manifolds in  the stratification. When $\nmanif = 0$, we use the second method, which doesn't require forming a complete list of manifolds. 

The overall goal of defining $\mathcal I$ is to access $\ngain{L}, \nlose{L}$, the set of Gain and Lose neighbours for a manifold $M_L$ defined by labels $L$, when requested. 
How we compute the neighbours depends on which method is used to construct $\mathcal I$. 

In the first case, when $\nmanif > 0$, we provide a list of labels in the form of a matrix, $\texttt{Llist}$, whose $i$th row $L\e{i} = \texttt{Llist}(i,:)$ gives the labels for the $i$th manifold in the stratification. The algorithm pre-computes $\ngain{L\e{i}}, \nlose{L\e{i}}$ by checking whether for each label $L\e{j} = \texttt{Llist}(j,:)$ with $j\neq i$, $L\e{j}$ differs from $L\e{i}$ by at most one equation. The sets $\ngain{L\e{i}}, \nlose{L\e{i}}$ are stored as lists of indices. 

In the second case, when $\nmanif = 0$, the algorithm computes $\ngain{L}, \nlose{L}$ on the fly. In the setup we must provide a set of flags $\texttt{FixFcns} \in \{C_{\rm fix}, C_{\rm vary}\}^{\nfcns}$ with $C_{\rm fix} = 1$, $C_{\rm vary}=0$, such that $\texttt{FixFcns}(i) = C_{\rm fix}$ if the label for function $q_i$ cannot change from what is provided initially, and $\texttt{FixFcns}(i) = C_{\rm vary}$  if the label for function $q_i$ can vary between constraint / inequality. Then, given labels $L$, the algorithm constructs $\ngain{L}$ by considering all possible ways to turn a single inequality into an equation (by performing a single flip $C_{\rm eq} \to C_{\rm in}$ for each function that is allowed to vary), and it constructs $\nlose{L}$ by considering all possible ways to turn a single equation into an inequality (by performing a single flip $C_{\rm in} \to C_{\rm eq}$ for each function that is allowed to vary.)

\subsection{Pseudocode} \label{sec:pseudocode}

The main algorithm is $\Call{SampleStrat}{}$ (Algorithm \ref{alg:samplestrat}).  
This takes as input a state $(X_n,L_n)$ representing the configuration $X_n\in \R^n$ and labels $L_n\in\{C_{\rm eq}, C_{\rm in}, C_{\rm none}\}^{\nfcns}$ of the $n$th step of the Markov chain, and outputs the $n+1$th state $X_{n+1},L_{n+1}$. This algorithm in turn calls $\Call{LamPropose}{}$ (Algorithm \ref{alg:lampropose}) to generate a proposal label $L$ and $\Call{VPropose}{}$ (Algorithm \ref{alg:vpropose})  to generate a proposal step $v$ in the tangent space. It then calls $\Call{TakeStep}{}$ (Algorithm \ref{alg:takestep}) to move in the direction $v$ and project back to the manifold, using projection methods \Call{NES}{} or \Call{NES-L}{} (also Algorithm \ref{alg:takestep}) to obtain a proposal $y\in M_L$.  $\Call{SampleStrat}{}$ then performs a variety of checks on the proposal $y,L$ to determine whether to accept it or reject it.

Most of these algorithms require information about the tangent space to a manifold, in the form of a matrix $T_{x,L}$, or sometimes $T_{x,L,v}$. Recall that $T_{x,L},T_{x,L,v}$ are matrices whose columns form an orthonormal basis of spaces $\mathcal T_{x,L},\mathcal T_{x,L,v}$ respectively, where $\mathcal T_{x,L}$ is the tangent space to $M_L$ at $x$, and $\mathcal T_{x,L,v}$ is the part of the tangent space that is orthogonal to vector $v\in \mathcal T_{x,L}$. Algorithm \ref{alg:tangent} shows how these matrices are calculated, using the QR decomposition of a matrix of gradients of the form \eqref{Q}. 
Once the matrices have been calculated for a state $x,L$, they are stored together with that state and reused wherever needed.

We remark that \Call{NES}{} and \Call{NES-L}{} perform essentially the same calculation, namely they solve a system of equations $q_k(x+Qa) = 0$ for $k\in I$ where $I$ is some subset of $\{1,\ldots,\nfcns\}$. The only differences between the algorithms are (i) they use a different initial condition for Newton's method, and (ii) the ordering of the constraints is different. We prefer to keep these functions separated to make it easier to interpret their output.

In the pseudocode that follows, we don't explicitly write out all the variables that must be passed to each function. Our implementation makes heavy use of objects, so we can efficiently pass all the information associated with a state $x,L$, such as $T_{x,L}$, $\ngain{L}$, $\nlose{L}$, $m_L$, $d_L$, and the values of the equations, inequalities, and gradients of the equations. We use the function argument $(x,L)$ somewhat loosely to mean all information associated with that particular state that has been calculated so far.



\algdef{SN}[propose]{StartPropose}{EndPropose}{\textbf{Proposal:}}{}
\algdef{SN}[project]{StartProject}{EndProject}{\textbf{Projection Check:}}{}
\algdef{SN}[inequality]{StartInequality}{EndInequality}{\textbf{Inequality Check:}}{}
\algdef{SN}[MH]{StartMH}{EndMH}{\textbf{Metropolis-Hastings Step:}}{}
\algdef{SN}[checkproject]{StartCheckProject}{EndCheckProject}{\textbf{Projection Check:}}{}
\algdef{SN}[reverse]{StartReverse}{EndReverse}{\textbf{Reverse Projection:}}{}

\begin{algorithm}[H]
  \caption{Given $x=X_n,L_n$, generate the next point $X_{n+1},L_{n+1}$ in the Markov chain}
  \label{alg:samplestrat}
\begin{algorithmic}[1]
\State {Parameters:} $\sigma$, $\lambda\lose$, $\lambda\gain$, $\sigma\bdy$, $\sigtan$
\Procedure{SampleStrat}{$x=X_n,L_n$} 
    \StartPropose \Comment{Generate a proposal $y,L$ and $\texttt{movetype}\in\{\texttt{Same,Gain,Lose}\}$}
        \State $\texttt{movetype},L =  \Call{LamPropose}{x,L_n}$  \Comment{Propose type of move and new labels}
       \State $v= \Call{VPropose}{x,L_0,\texttt{movetype},L}$ \Comment{Propose step in tangent space}
       \State $y,\texttt{newtonflag}, \alpha =$ \Call{TakeStep}{$x,L_n, \texttt{movetype}, L, v$}
       \State \Comment{\parbox[t]{0.83\linewidth}{Takes a step $v$ in the tangent space starting at $x$ and projects to manifold $M_L$.\\
       If $\texttt{movetype}\in\{\texttt{Same},\texttt{Gain}\}$, solves for $y\in M_L$ s.t. $y=x+v+w$, $w\perp \mathcal T_{x,L}$\\
       If $\texttt{movetype} == \texttt{Lose}$, solves for $\alpha\in\R,y\in M_L$ s.t. $y=x+\alpha v + w$, $w\perp \mathcal T_{x,L_n}$}}
    \EndPropose\Statex
    \StartProject
        \If{\texttt{newtonflag} == \texttt{fail}} \Comment{Projection step (\Call{NES}{} or \Call{NES-L}{}) failed to converge}
           \State Reject proposal: $X_{n+1}=X_n, L_{n+1}=L_n$. Return. 
        \EndIf
        \If{\texttt{movetype == Lose} and $\alpha < 0$}  \Comment{Lose move went in opposite $v$ direction}
            \State Reject proposal: $X_{n+1}=X_n, L_{n+1}=L_n$. Return. 
        \EndIf
    \EndProject\Statex
    \StartInequality
        \If{$q_i(y) < 0$ for some $i$ s.t. $L(i) = C_{\rm in}$}
            \State Reject proposal: $X_{n+1}=X_n, L_{n+1}=L_n$. Return. 
        \EndIf
    \EndInequality\Statex
    \StartMH
        \State Compute $\Lambda\same(x),\Lambda\gain(x),\Lambda\lose(x)$; $\Lambda\same(y),\Lambda\gain(y),\Lambda\lose(y)$ via \eqref{LambdaAlg2}
        \State Get densities $\lambda_{IJ}(x)h_{IJ}(y|x)$ for forward move via  \eqref{SameDens}, \eqref{GainDens}, or \eqref{LoseDens}
       \State Set \texttt{movetyperev} = \texttt{Same, Lose, Gain} for \texttt{movetype} =  \texttt{Same, Gain, Lose} respectively
        \If{\texttt{movetyperev} $==$ \texttt{Same} or \texttt{Gain}} \Comment{Construct steps for reverse move }
            \State  Find $v'\in \mathcal T_{y,L_n}$ such that $x=y+v'+w'$ with $w'\perp \mathcal T_{y,L_n}$
            \Comment{See \eqref{SameV}}
        \ElsIf{\texttt{movetyperev} == \texttt{Lose}}   
           \State  Find $v'\in \mathcal T_{y,L}$ such that $x=y+v'+w'$ with $w'\perp \mathcal T_{y,L}$ \Comment{Similar to \eqref{SameV}}
           \State $\alpha' = |v'|$, $v'\gets v'/|v'|$ 
       \EndIf
          \State Get densities $\lambda_{IJ}(y)h_{IJ}(x|y)$ for reverse move via  \eqref{SameDens},\eqref{GainDens}, or \eqref{LoseDens}, using $v',\alpha'$
          \State \Comment{Temporarily set accessibility factor $1_{A_{\texttt{movetype}}(y)}(x)=1$}
          \State Compute acceptance probability $a = a(y,L|x,L_n)$ via \eqref{acc} \Comment{Also requires $f(x),f(y)$}
          \If{$U\sim \text{Unif}([0,1]) > a$}
              \State Reject proposal: $X_{n+1}=X_n, L_{n+1}=L_n$. Return. 
          \EndIf
    \EndMH\Statex
    \StartReverse
        \State $x',\texttt{newtonflagrev}, \alpha' =$ \Call{TakeStep}{$y,L, \texttt{movetyperev}, L_n, v'$}
       \State \Comment{\parbox[t]{0.85\linewidth}{Takes a step $v'$ in the tangent space starting at $y$ and projects to manifold $M_{L_n}$.\\
       If $\texttt{movetyperev}\in\{\texttt{Same},\texttt{Gain}\}$, solves for $x'{\in}M_{L_n}$ s.t. $x'=y+v'+w'$, $w'\perp \mathcal T_{y,L_n}$ \\
       If $\texttt{movetyperev}{==}\texttt{Lose}$, solves for $\alpha'{\in}\R,x'{\in}M_{L_n}$ s.t. $x'=y'+\alpha' v' + w'$, $w'\perp \mathcal T_{y,L}$}}
       \If{$\texttt{newtonflagrev} == \texttt{fail}$}   \Comment{Reverse projection step failed to converge}
           \State Reject proposal: $X_{n+1}=X_n, L_{n+1}=L_n$. Return. 
        \EndIf
               \If{$x'\neq x$}  \Comment{Projection converged, but to wrong point}
          \State Reject proposal: $X_{n+1}=X_n, L_{n+1}=L_n$. Return. 
       \EndIf
        \If{\texttt{movetyperev == Lose} and $\alpha' < 0$}
            \State Reject proposal: $X_{n+1}=X_n, L_{n+1}=L_n$. Return. 
        \EndIf
    \EndReverse\Statex
    \State \textbf{Accept proposal:} $X_{n+1}=y, L_{n+1}=L$.  Copy $T_{y,L}$, $\ngain{L}$, $\nlose{L}$  to $X_{n+1}$. 
    Return. 
\EndProcedure
\State \textbf{end procedure}
  \end{algorithmic}
\end{algorithm}

\begin{algorithm}[H]
  \caption{Propose a type of move and new labels}
  \label{alg:lampropose}
  \begin{algorithmic}[1]
  \State {Parameters:} $\lambda\lose$, $\lambda\gain$, $\sigma\bdy$
  \Procedure{LamPropose}{$x,L_0$} 
  \State Get $\ngain{L_0}$, $\nlose{L_0}$ \Comment{Lists of labels, e.g. $\nlose{L_0} = \{L'_1,L'_2,\ldots,L'_{|\nlose{L_0}|}\}$}
  \State\Comment{\parbox[t]{0.85\linewidth}{If $\texttt{nmanif}{>}0$ these lists are pre-computed for each manifold; if $\texttt{nmanif}{=}0$ they are obtained by adding each inequality or subtracting each constraint to current labels}}
  \Statex
  \State \emph{Compute nearby Lose neighbours as in \eqref{Nlosesig}:}
  \State $\nlosesig{L_0,x} = \emptyset$ 
  \For{$i\gets 0$,  $i<|\nlose{L_0}|$}
      \If{\Call{DistMin}{$x,L_0,k\e{\rm change}_i$}} \Comment{$q_{k\e{\rm change}_i}$ is extra constraint in $L'_i$} 
          \State $\nlosesig{L_0,x} \gets \nlosesig{L_0,x}\cup \{L'_i\}$
      \EndIf
  \EndFor \State \textbf{end for}
  \State $n\gain = |\ngain{L_0}|$, $n\lose = |\nlosesig{L_0,x}|$ 
  \Statex
  \State \emph{Compute $\Lambda\same(x),\Lambda\gain(x),\Lambda\lose(x)$ via \eqref{LambdaAlg2}: }
  \State $\Lambda\gain(x) = \lambda\gain$, $\Lambda\lose(x) = \lambda\lose$ 
  \If{$n\gain == 0$}
       $\Lambda\gain(x) = 0$
  \EndIf
  \If{$n\lose == 0$}
       $\Lambda\lose(x) = 0$
  \EndIf
  \State $\Lambda\same(x) = 1 - \Lambda\gain(x) - \Lambda\lose(x)$
  \Statex
  \State \emph{Choose a move with the desired probabilities: }
  \State Generate $U\sim \text{Unif}([0,1])$  
  \If{$U  < \Lambda\same(x)$}  \Comment{Choose a Same move}
      \State \texttt{movetype} = \texttt{Same}
      \State $L_1 = L$
  \ElsIf{$U  < \Lambda\same(x) + \Lambda\gain(x)$}  \Comment{Choose a Gain move}
      \State \texttt{movetype} = \texttt{Gain}
       \State $L_1 \sim \text{Uniform}(\nlosesig{L_0,x})$ \Comment{Choose labels uniformly from Gain neighbours}
  \Else     \Comment{Choose a Lose move}
      \State \texttt{movetype} = \texttt{Lose}
      \State $L_1 \sim \text{Uniform}(\nlosesig{L_0,x})$ \Comment{Choose labels uniformly from nearby Lose neighbours}
  \EndIf
  \State \textbf{return} \texttt{movetype}, $L_1$
  \EndProcedure
  \State \textbf{end procedure}
    \Statex
   \State \emph{Estimate the minimum distance from $x$ to boundary $\{z:q_k(z)=0\}$ over all directions $v\in\mathcal T_{x,L}$.}
  \Function{DistMin}{$x,L,k$} 
      \State \textbf{return} $|q_k(x)|/|T_{x,L}\grad q_k(x)|$ \Comment{Estimated by linearizing $q_i$, as in \eqref{hopt}}
  \EndFunction
  \State \textbf{end function}
   \end{algorithmic}
\end{algorithm}

\begin{algorithm}[H]
  \caption{Propose a step $v$ in the tangent space}
  \label{alg:vpropose}
  \begin{algorithmic}[1]
  \State {Parameters:} $\sigma$, $\sigma\bdy$, $\sigtan$
  \Procedure{VPropose}{$x,L_0,\texttt{movetype},L_1$} 
      \If{\texttt{movetype} == \texttt{Same}}
          \If{$d_{L_0} == 0$} $v = (0,\ldots,0)\in \R^n$  \Comment{Dimension of $M_{L_0}$ is $0$; remain on current point}
          \Else
                 \State Generate $r_i \sim \sigma\cdot N(0,1)$ for $i=1,\ldots,d_{L_0}$  \Comment{$r$ is vector of i.i.d. normals}
              \State $v = T_{x,L_0}r$ \Comment{Step in tangent space $\mathcal T_{x,L_0}$}
          \EndIf
      \EndIf \State\textbf{end if}
      \If{\texttt{movetype} == \texttt{Gain}}
         \State Let $i_{\rm change}$ be  index of constraint removed from $L_0$ to $L_1$   
         \State Calculate $T_{x,L_1}$  \Comment{Basis for tangent space with $q_{i_{\rm change}}$ removed}
         \State $u_n = T_{x,L_1}T_{x,L_1}^T \grad q_{i_{\rm change}}$, \; $u_n\gets u_n/|u_n|$ \Comment{See \eqref{un}}
          \State Generate $U\sim \text{Unif}([0,1])$
           \If{move is one-sided} $v_n = \sigma\bdy U$  \Comment{Component of step in normal direction}
           \EndIf
          \If{move is two-sided} $v_n = \sigma\bdy (2U-1)$
          \EndIf
          \State $v = v_n u_n$  \Comment{Step in normal direction}
          \If{$d_{L_0} \geq 1$}  \Comment{There exist tangential directions}
              \State Generate $r_i \sim \sigtan u_n \cdot N(0,1)$ for $i=1,\ldots,d_{L_0}$ \Comment{Components of tangential steps}
              \State $v\gets v + T_{x,L_0}r$  \Comment{Add tangential steps to $v$; see \eqref{GainStep}}
          \EndIf
      \EndIf \State\textbf{end if}
      \If{\texttt{movetype} == \texttt{Lose}}
          \State Let $i_{\rm change}$ be  index of constraint added from $L_0$ to $L_1$   
          \State $u_n = T_{x,L_0}T_{x,L_0}^T \grad q_{i_{\rm change}}$, \; $u_n\gets u_n/|u_n|$  \Comment{$u_n = v_{\rm opt}$, as in \eqref{vopt}}
          \If{move is two-sided and $q_{i_{\rm change}}(x) < 0$}  $u_n\gets -u_n$ \EndIf
          \State $v = -u_n$  \Comment{Normal component of step}
          \If{$d_{L_0} \geq 2$}  
              \State Construct $T_{x,L_0,u_n}$
              \State Generate $r_i\sim \sigtan\cdot N(0,1)$ for $i=1,\ldots,d_{L_0}-1$  \Comment{Tangential components of step}
              \State $v\gets v + T_{x,L_0,u_n}r$   \Comment{See \eqref{vlose}}
              \State $v\gets v/|v|$ \Comment{Send back direction only}
             \EndIf
      \EndIf \State\textbf{end if}
        \State\textbf{return} $v$
  \EndProcedure
  \State \textbf{end procedure}
   \end{algorithmic}
\end{algorithm}

\begin{algorithm}[H]
  \caption{Take a step $v$ in the tangent space from $x$ and project back to manifold}
  \label{alg:takestep}
  \begin{algorithmic}[1]
    \Statex
  \State {Parameters:} $\texttt{tol}, \texttt{MaxIter}$ 
  \Procedure{TakeStep}{$x,L_0, \texttt{movetype}, L_1, v$} 
      \If{$\texttt{movetype}\in\{\texttt{Same},\texttt{Gain}\}$} 
          \State  $z = x+v$
          \If{$m_{L_1}>0$}    \Comment{$m_{L_1}$ = number of equations in $L_1$}
              \State Compute $Q_{L_1}$ via \eqref{Q}
              \State  $y,\texttt{newtonflag} = \Call{NES}{z,Q_{L_1},L_1}$  \Comment{Project back to manifold}
          \Else \Comment{No projection needed if $M_{L_1}$ is defined only by inequalities}
              \State $y=z$, \texttt{newtonflag}=\texttt{success}
          \EndIf  
       \ElsIf{\texttt{movetype} == \texttt{Lose}}  
          \State Construct $Q\e{v} = (Q_{L_0} \:|\: v)$, with $Q_{L_0}$ as in \eqref{Q}
          \State Let $i_{\rm change}$ be  index of constraint added from $L_0$ to $L_1$     \Comment{$q_{i_{\rm change}}$ is added equation}
          \State $y,\alpha,\texttt{newtonflag} = \Call{NES-L}{x,Q\e{v},L_0,i_{\rm change},v}$  \Comment{Modified projection for a Lose move}
       \EndIf  
       \State \textbf{return } $y,\texttt{newtonflag},(\alpha)$
  \EndProcedure
  \State \textbf{end procedure}
  \Statex
  
  \State \emph{Newton's method to solve $\{q_k(z+Qa)=0 \;\forall k{:}\:L(k){=}C_{\rm eq}\}$ for $a\in \R^{m_{L}}$}
  \Function{NES}{z,Q,L}  
      \State Set $a = (0,0,\ldots, 0)\in \R^{m_{L}}$ \Comment{Deterministic initial condition}
      \For{$i\gets 0, i < \texttt{MaxIter}$}
          \State$F\gets (q_k(z+Qa))_{k:L(k)=C_{\rm eq}}$  \Comment{Vector containing current values of equations}
          \If{$\max_j |F_j| < \texttt{tol}$} \Comment{Success! Solution converged}
              \State $y = z+Qa$ 
              \State \textbf{return} $y,\texttt{newtonflag} = \texttt{success}$
          \EndIf
          \State $\tilde Q = (\grad q_k(z+Qa))_{k:L(k)=C_{\rm eq}}$ \Comment{Current value of constraint gradients in $L$, as in \eqref{Q}}
          \State $J = \tilde Q^TQ$ \Comment{Jacobian of system to be solved}
          \State Solve $J\Delta a = -F$ for $\Delta a$ \Comment{We used LU decomposition with partial pivots}
          \State $a \gets a + \Delta a$
      \EndFor \State \textbf{end for}
      \State \textbf{return} $y=\texttt{NaN},\texttt{newtonflag} = \texttt{fail}$  \Comment{If we got this far, solution didn't converge}
  \EndFunction
  \State \textbf{end function} 
  \Statex
  
  \State \emph{Newton's method to solve $\{q_k(x+Qa)=0 \;\forall k{:}\:L_0(k){=}C_{\rm eq}, k=i_0\}$ for $a\in \R^{m_{L_0}+1}$}
  \Function{NES-L}{$x,Q,L_0,i_{0},v$} 
     \State  $h = -q_{i_0}(x)/\grad q_{i_0}(x)\cdot v$ \Comment{Estimated distance to $q_{i_0}(y){=}0$ in direction $v$; see \eqref{hv}}
     \State Set $a = (0,0,\ldots,0,h)\in \R^{m_{L_0}+1}$ 
     \Comment{Linearizing about 0 gives $q_{i_0}(hv)=0+O(h^2)$}
     \For{$i\gets 0, i < \texttt{MaxIter}$}
          \State $F\gets \left((q_k(x+Qa))_{k:L_0(k)=C_{\rm eq}}, q_{i_{0}}(x+Qa)\right)$  \Comment{Current values of equations in $L_0$ and $i_0$}
          \If{$\max_j |F_j| < \texttt{tol}$} \Comment{Success! Solution converged}
              \State $y = x+Qa$ 
              \State $\alpha = a_{m_{L_0}+1}$  \Comment{$\alpha$ is the last element of $a$}
              \State \textbf{return} $y,\alpha,\texttt{newtonflag} = \texttt{success}$
          \EndIf
           \State $\tilde Q = \left((\grad q_k(x+Qa))_{k:L_0(k)=C_{\rm eq}} \:|\: \grad q_{i_0}(x+Qa)\right)$ \Comment{\eqref{Q} for $L_0$ with additional column}
          \State $J = \tilde Q^TQ$ \Comment{Jacobian of system to be solved}
          \State Solve $J\Delta a = -F$ for $\Delta a$ \Comment{We used LU decomposition with partial pivots}
          \State $a \gets a + \Delta a$
     \EndFor \State \textbf{end for}
     \State \textbf{return} $y=\texttt{NaN},\texttt{newtonflag} = \texttt{fail}$  \Comment{If we got this far, solution didn't converge}
  \EndFunction
  \State \textbf{end function}
   \end{algorithmic}
\end{algorithm}

\begin{algorithm}[H]
\caption{Calculate an orthonormal basis of a tangent space}\label{alg:tangent}
  \begin{algorithmic}[1]
  \State \emph{Calculates $T = T_{x,L}$, a matrix whose columns form an orthonormal basis of $\mathcal T_{x,L}$}
  \Function{Tan}{$x,L$}
      \If{$0<d_L<n$} 
          \State Calculate $Q_L$ in \eqref{Q}
          \State $W = QR(Q_L)$ \Comment{QR decomposition}
          \State \Comment{\parbox[t]{0.85\linewidth}{$W$ has block form $W = (N \:|\: T)$, where $N\in \R^{n\times m_L}$ is an orthonormal basis of $\mathcal N_{x,L}$, and $T=T_{x,L}\in \R^{n\times d_L}$ is an orthonormal basis of $\mathcal T_{x,L}$}}
          \State Set $T$ to be the last $d_L$ columns of $W$
       \ElsIf{$d_L == n$} \Comment{No equations; moving in ambient Euclidean space}
           \State $T = I\in \R^{n\times n}$ \Comment{The $n\times n$ identity matrix}
       \ElsIf{$d_L == 0$} \Comment{We're on a 0-dimensional point}
           \State Don't need $T$, since there is no tangent space. 
       \EndIf\State\textbf{end if}
       \State \textbf{return } $T$
  \EndFunction
  \State \textbf{end function}
  \Statex
  \State \emph{Calculates $T_v = T_{x,L,v}$, a matrix whose columns form an orthonormal basis of $\mathcal T_{x,L,v}$. \\
  The function assumes that $v\in \mathcal T_{x,L}$; if not, the function should must be modified.}
  \Function{TanV}{$x,L,v$}
       \If{$d_L\geq 2$} 
           \State $T = \Call{Tan}{x,L}$
           \State $P = I - vv^T/|v|^2$, here $I\in \R^{n\times n}$ is identity matrix 
           \State \Comment{$P$ is the orthogonal projection matrix onto the space perpendicular to $v$}
           \State $W = QR(P\cdot T)$ \Comment{QR decomposition}
           \State Set $T_v$ to be the first $d_L-1$ columns of $W$. 
       \Else \Comment{We're on a point or a line}
           \State Don't need $T_{x,L,v}$ since there's at most 1 direction 
       \EndIf
       \State \textbf{return } $T_v$
  \EndFunction
  \State \textbf{end function}
 \end{algorithmic}
\end{algorithm}




\bibliographystyle{plainnat}
\bibliography{stratsample}

\end{document}